\documentclass[reqno]{amsart} 
\usepackage{amsmath} 
\usepackage{amsfonts}   
\usepackage{eufrak}      
\usepackage{amscd}      
\usepackage{amsthm}      
\usepackage{epsfig}      
\usepackage{amstext}      
\usepackage[all]{xy}      
\usepackage{amssymb}

\newcommand{\C}{{\mathbb C}}

\newcommand{\R}{{\mathbb R}}
\newcommand{\T}{{\mathbb T}}
\newcommand{\Z}{{\mathbb Z}}
\newcommand{\Q}{{\mathbb Q}}
   
\theoremstyle{plain}
\newtheorem{thm}{Theorem}[section]
\newtheorem{prop}[thm]{Proposition}
\newtheorem{lemma}[thm]{Lemma}
\newtheorem{cor}[thm]{Corollary}
\theoremstyle{definition}

\theoremstyle{remark}

\title{Classification of certain simple $C^*$-algebras with torsion in $K_1$} 
\author{Jesper Mygind}
\date{\today}
\address{The Fields Institute, 222 College Street, Toronto,
Ontario M5T 3J1, Canada}
\email{jmygind@hotmail.com}

\begin{document}

\begin{abstract}
We show that the Elliott invariant is a classifying invariant for the class
of $C^*$-algebras that are simple unital infinite dimensional inductive limits
of sequences of finite direct sums of building blocks of the form
$$
\{f\in C(\T)\otimes M_n: f(x_i)\in M_{d_i},\ i=1,2,\dots,N\},
$$
where $x_1,x_2,\dots,x_N\in\T$, $d_1,d_2,\dots,d_N$ are integers dividing $n$,
and $M_{d_i}$ is embedded unitally into $M_n$. Furthermore we prove existence
and uniqueness theorems for *-homomorphisms between such algebras and we
identify the range of the invariant.
\end{abstract}


\maketitle

\tableofcontents

\section{Introduction}

During the last decade the Elliott invariant has been used with amazing success
to classify simple unital $C^*$-algebras (see e.g \cite{EAIS}, \cite{EATS},
\cite{NT}, \cite{TATD}, \cite{JSSI}, \cite{JSSP}). This project is part of
Elliott's program which has the ambitious goal of a classification result
for all separable nuclear $C^*$-algebras by invariants of $K$-theoretical
nature.

The goal of the present paper is to unify and generalize classification
results due to Thomsen \cite{TATD} and Jiang and Su \cite{JSSP}. In order to
achieve this we will unfortunately have to consider the rather complicated
building blocks defined in the abstract. Our main result
(see Theorem \ref{i}) is the following:

\begin{thm}
The Elliott invariant is a classifying invariant for the class of unital simple
infinite dimensional inductive limits of sequences of finite direct sums of
building blocks.
\end{thm}

The main ideas of the proof are similar to those of Thomsen \cite{TATD}
who considers the simpler case $d_1=d_2=\dots=d_N$. The technical problems
are greater in our case, and in particular the possible lack of
projections in our building blocks (see Lemma \ref{proj}) means there is no
straightforward generalization of Thomsen's proof.

Let us introduce the notation used in this paper
before we describe our results in greater detail.
Recall that for a unital $C^*$-algebra $A$ the Elliott invariant consists of
the ordered group $K_0(A)$ with order unit, the group $K_1(A)$, the compact
convex set $T(A)$ of tracial states, and the restriction map
$r_A:T(A)\to SK_0(A)$, where $SK_0(A)$ denotes the state space of $K_0(A)$.

Let $A$ be a unital $C^*$-algebra. Let $\text{Aff}\, T(A)$ denote the order
unit space of all continuous real-valued affine functions on $T(A)$.
Let $\rho_A:K_0(A)\to\text{Aff}\, T(A)$ be the group homomorphism
$$
\rho_A(x)(\omega)=r_A(\omega)(x),\quad \omega\in T(A),\ x\in K_0(A).
$$
Let $U(A)$ denote the unitary group of $A$
and let $DU(A)$ denote its commutator subgroup, i.e the group generated by
all unitaries of the form $uvu^*v^*$, $u,v\in U(A)$.
If $A$ is a unital inductive limit of a sequence of finite direct sums of
building blocks then there is a natural short exact sequence of abelian groups
(see section 5)
$$
  \begin{CD}
      0 @>  >> \text{Aff}\, T(A)/\overline{\rho_A(K_0(A))}  @>\lambda_A  >>
U(A)/\overline{DU(A)} @>\pi_A  >> K_1(A) @>  >> 0 \\   
  \end{CD}
$$
that splits (unnaturally).
The group $U(\cdot)/\overline{DU(\cdot)}$ was introduced into the
classification program by Nielsen and Thomsen \cite{NT}.

Let $A$ and $B$ be unital $C^*$-algebras.
An affine continuous map $\varphi_T:T(B)\to T(A)$ gives rise to
a linear positive order unit preserving map
${\varphi_T}_*:\text{Aff}\, T(A)\to \text{Aff}\, T(B)$ by setting
${\varphi_T}_*(f)=f\circ\varphi_T$ for $f\in\text{Aff}\, T(A)$.
If furthermore ${\varphi_T}_*\circ\rho_A=\rho_B\circ\varphi_0$ for
some group homomorphism $\varphi_0:K_0(A)\to K_0(B)$ then $\varphi_T$
induces a group homomorphism
$$
\widetilde{\varphi_T}:\text{Aff}\, T(A)/\overline{\rho_A(K_0(A))}\to
\text{Aff}\, T(B)/\overline{\rho_B(K_0(B))}.
$$

Let $\psi:A\to B$ be a unital *-homomorphism. Let $\psi^*:T(B)\to T(A)$ be the
affine continuous map given by $\psi^*(\omega)=\omega\circ\psi$,
$\omega\in T(B)$. Define $\widehat{\psi}:\text{Aff}\, T(A)\to\text{Aff}\, T(B)$
by $\widehat{\psi}=(\psi^*)_*$. Note that
$\widehat{\psi}(f)(\omega)=f(\psi^*(\omega))$.
Since $\widehat{\psi}\circ\rho_A=\rho_B\circ\psi_*$ on $K_0(A)$, we see that
$\psi$ gives rise to a group homomorphism
$$
\widetilde{\psi}:\text{Aff}\, T(A)/\overline{\rho_A(K_0(A))}\to
\text{Aff}\, T(B)/\overline{\rho_B(K_0(B))}.
$$
Let $\psi^{\#}:U(A)/\overline{DU(A)}\to U(B)/\overline{DU(B)}$
be the homomorphism induced by $\psi$.

Besides the Elliott invariant, two other invariants will be crucial in
the proof of the classification theorem, namely $U(\cdot)/\overline{DU(\cdot)}$
and R\o rdam's $KL$-bifunctor \cite{RKL}. These invariants are both
determined by the Elliott invariant for the $C^*$-algebras under
consideration, and are therefore useless as additional isomorphism
invariants. They are, however, not determined canonically. This means that
*-homomorphisms
(or even automorphisms) between such $C^*$-algebras that agree on
the Elliott invariant may fail to be approximately
unitarily equivalent because they may act differently on these additional
invariants. This was demonstrated by Nielsen and Thomsen \cite[section 5]{NT}
for $U(\cdot)/\overline{DU(\cdot)}$ and by Dadarlat and Loring
\cite[p. 375-376]{DLUMCT} for $KL$.

It is therefore necessary to include these invariants in the following
uniqueness theorem (see Theorem \ref{aue}):
\begin{thm}
Let $A$ and $B$ be unital inductive limits of sequences of finite direct
sums of building blocks, with $A$ simple. Two unital *-homomorphisms
$\varphi,\psi:A\to B$ with $\varphi^*=\psi^*$ on $T(B)$,
$\varphi^{\#}=\psi^{\#}$ on $U(A)/\overline{DU(A)}$, and $[\varphi]=[\psi]$
in $KL(A,B)$ are approximately unitarily equivalent.
\end{thm}

Let $KL(A,B)_T$ denote the set of elements $\kappa\in KL(A,B)$ for which
the induced map $\kappa_*:K_0(A)\to K_0(B)$ preserves the order unit and
for which there exists an affine continuous map $\varphi_T:T(B)\to T(A)$
such that $r_B(\omega)(\kappa_*(x))=r_A(\varphi_T(\omega))(x)$ for
$x\in K_0(A),\ \omega\in T(B)$.

Let $A$ and $B$ be e.g simple unital inductive limits of sequences of finite
direct sums of building blocks. It turns out, perhaps surprisingly, that there
is a connection between $KL(A,B)$ and the torsion subgroups of
$U(A)/\overline{DU(A)}$ and $U(B)/\overline{DU(B)}$, see section 10. If
$\varphi,\psi:A\to B$ are unital *-homomorphisms with $[\varphi]=[\psi]$
in $KL(A,B)$ and if $x$ is an element of finite order in the group
$U(A)/\overline{DU(A)}$, then $\varphi^{\#}(x)=\psi^{\#}(x)$ in
$U(B)/\overline{DU(B)}$. More generally, an element $\kappa\in KL(A,B)_T$
gives rise to a group homomorphism
$$
s_{\kappa}:Tor(U(A)/\overline{DU(A)})\to Tor(U(B)/\overline{DU(B)}).
$$
The map
$$
KL(A,B)_T\to
Hom\bigl(Tor(U(A)/\overline{DU(A)}),Tor(U(B)/\overline{DU(B)})\bigr),
$$
where $\kappa\mapsto s_{\kappa}$, is natural with respect to the Kasparov
product and must be taken into account in the existence theorem:

\begin{thm}
Let $A$ and $B$ be simple unital inductive limits of sequences of finite
direct sums of building blocks, with $B$ infinite dimensional.
Let $\varphi_T : T(B)\to T(A)$ be an affine continuous map, let
$\kappa\in KL(A,B)_T$ be an element such that
$$
r_B(\omega)(\kappa_*(x))=r_A(\varphi_T(\omega))(x),
\quad x\in K_0(A),\ \omega\in T(B),
$$
and let $\Phi:U(A)/\overline{DU(A)}\to U(B)/\overline{DU(B)}$ be a
homomorphism such that the diagram
$$
  \begin{CD}
      \text{Aff}\, T(A)/\overline{\rho_A(K_0(A))} @> \lambda_A >>
      U(A)/\overline{DU(A)} @> \pi_A >> K_1(A)    \\   
      @V\widetilde{\varphi_T} VV  @V\Phi VV      @VV\kappa_* V \\
      \text{Aff}\, T(B)/\overline{\rho_B(K_0(B))}
      @>> \lambda_B > U(B)/\overline{DU(B)} @>> \pi_B > K_1(B)
  \end{CD}
$$
commutes. Assume finally that
$$
s_{\kappa}(y)=\Phi(y),\quad y\in Tor(U(A)/\overline{DU(A)}).
$$
There exists a unital *-homomorphism $\psi:A\to B$ such that
$\psi^*=\varphi_T$ on $T(B)$, such that $\psi^{\#}=\Phi$ on
$U(A)/\overline{DU(A)}$, and such that $[\psi]=\kappa$ in $KL(A,B)$.
\end{thm}

The above theorem follows by combining the slightly more general
Theorem \ref{e} with Lemma \ref{infty}, Lemma \ref{torK_0} and
Theorem \ref{inj}. It should be noted that it is possible to prove this
existence theorem (and our classification theorem) for $K_0(A)$
non-cyclic without using the map $s_{\kappa}$, see Corollary \ref{exnc}
(or \cite{TATD}).

Let us finally describe the range of the invariant for the $C^*$-algebras
in our class. By combining Theorem \ref{rannc} and Corollary \ref{ran} we have
the following:

\begin{thm}
Let $G$ be a countable simple dimension group with order unit, $H$
a countable abelian group, $\Delta$ a compact metrizable Choquet simplex, and
$\lambda:\Delta\to SG$ an affine continuous extreme point preserving
surjection. There exists a simple unital inductive limit of a sequence of
finite direct sums of building blocks $A$ together with an
isomorphism $\varphi_0:K_0(A)\to G$ of ordered groups with order unit, an
isomorphism $\varphi_1:K_1(A)\to H$, and an affine homeomorphism
$\varphi_T:\Delta\to T(A)$ such that
$$
r_A(\varphi_T(\omega))(x)=\lambda(\omega)(\varphi_0(x)),
\quad\omega\in\Delta,\ x\in K_0(A)
$$
if and only if $G$ is non-cyclic, or $G$ is cyclic and $H$ can
be realized as an inductive limit of a sequence of the form
$$
  \begin{CD}
      \Z\oplus H_1  @> >> \Z\oplus H_2 @> >> \Z\oplus H_3 @> >> \dots \\   
  \end{CD}
$$
where each $H_k$ is a finite abelian group.
\end{thm}

Let $A$ be a simple unital inductive limit of a sequence of
finite direct sums of building blocks.
It is easy to see that $A$ is unital projectionless if and only if
$(K_0(A),K_0(A)^+,[1])\cong (\Z,\Z^+,1)$. Hence our classification theorem
can be applied to a large class of simple unital projectionless $C^*$-algebras,
including the $C^*$-algebra $\mathcal Z$ constructed
by Jiang and Su \cite{JSSP}.

It would be interesting if one could extend our classification result to a
class that contains simple unital projectionless $C^*$-algebras
with arbitrary countable abelian $K_1$-groups. This could probably be obtained
by considering building blocks with $\T$ replaced by a general $1$-dimensional
compact Hausdorff space. It would also be interesting if one could include the class of
$C^*$-algebras considered by Jiang and Su in \cite{JSSI}.

Let $A$ be a unital $C^*$-algebra. If $a\in A_{\text{sa}}$ we define
$\widehat{a}\in\text{Aff}\, T(A)$ by
$\widehat{a}(\omega)=\omega(a)$, $\omega\in T(A)$. It is well-known that
$a\mapsto\widehat a$ is a surjective map from $A_{\text{sa}}$ to
$\text{Aff}\, T(A)$.
Let $q_A':U(A)\to U(A)/\overline{DU(A)}$ be the canonical map. We equip the
abelian group $U(A)/\overline{DU(A)}$ with the quotient metric
$$
D_A\bigl(q_A'(u),q_A'(v)\bigr)=\inf\{\|uv^*-x\|:x\in\overline{DU(A)}\}.
$$
Denote by $d_A'$ the quotient metric on the group
$\text{Aff}\, T(A)/\overline{\rho_A(K_0(A))}$. This group can be
equipped with another metric which
gives rise to the same topology, namely
$$
d_A(f,g)=\begin{cases}
2 & d'_A(f,g)\ge \frac 12, \\
|e^{2\pi id_A'(f,g)}-1| & d'_A(f,g)< \frac 12,
\end{cases}
$$
see \cite[Chapter 3]{NT}.
Let $q_A:\text{Aff}\, T(A)\to \text{Aff}\, T(A)/\overline{\rho_A(K_0(A))}$
be the quotient map.
Let finally $s(A)$ be the smallest positive integer
$n$ for which there exists a unital *-homomorphism $A\to M_n$ (we set
$s(A)=\infty$ if $A$ has no non-trivial finite dimensional representations).

Let $gcd$ denote the greatest common divisor and $lcm$
the least common multiple of a set of positive integers. Let $Tr$ denote the
(unnormalized) trace on a matrix algebra (i.e the number obtained by adding
the diagonal entries). If
$$
  \begin{CD}
      A_1  @> \alpha_1 >> A_2 @> \alpha_2 >>
A_3 @>\alpha_3  >> \dots \\   
  \end{CD}
$$
is a sequence of $C^*$-algebras and *-homomorphisms with inductive limit $A$,
we let
$\alpha_{n,m}=\alpha_{m-1}\circ\alpha_{m-2}\circ\dots
\circ\alpha_n:A_n\to A_m$
when $m>n$. We set $\alpha_{n,n}=id$ and let $\alpha_{n,\infty}:A_n\to A$
denote the canonical map.

{\it Acknowledgements:} Most of this work was carried out while I was a Ph.D
student at Aarhus University. I would like to thank the many people at the
Department of Mathematics that have helped me during my studies, both
mathematically and otherwise. In particular I would like to
express my gratitude towards my Ph.D advisor Klaus Thomsen for many inspiring
and helpful conversations. Special thanks to George Elliott for valuable
conversations and for his interest in my work. Finally, I would like to
thank Mikael R\o rdam for some useful suggestions.

\section{Building blocks}

Let $\T$ denote the unit circle of the complex plane. We will equip
$\T$ with the metric
$$
\rho(e^{2\pi is},e^{2\pi it})=\min_{k\in\Z} |s-t+k|
$$
which is easily seen to be equivalent to the usual metric on $\T$
inherited from $\C$.

As in \cite{NT} we say that a tuple $(a_1,a_2,\dots,a_L)$ of elements
from $\T$ is naturally numbered if there exist numbers
$s_1,s_2,\dots,s_L\in [0,1[$ such that $s_1\le s_2\le\dots\le s_L$ and
$a_j=e^{2\pi is_j}$, $j=1,2,\dots,L$.

We define a building block to be a $C^*$-algebra of the form
$$
A(n,d_1,d_2,\dots,d_N)=
\{f\in C(\T)\otimes M_n:f(x_i)\in M_{d_i},\ i=1,2,\dots,N\},
$$
where $(x_1,x_2,\dots,x_N)$ is a naturally numbered tuple of (different)
points in $\T$, $d_1,d_2,\dots,d_N$ are integers dividing $n$, and $M_{d_i}$
is embedded unitally into $M_n$, e.g via the *-homomorphism
$$
a\mapsto \text{diag}(\underbrace{a,a,\dots,a}_{\frac n{d_i}\ \text{times}}).
$$
The points $x_1,x_2,\dots,x_N$ will be called the exceptional points of $A$.
By allowing $d_i=n$ we may always assume that $N\ge 2$. It will also be
convenient to always assume that $1$ is not an exceptional point.

For every $i=1,2,\dots,N$, evaluation at $x_i$ gives rise to a unital
*-homomorphism from $A$ to $M_{d_i}$ which will be denoted by $\Lambda_i$, or
sometimes $\Lambda_i^A$. If $s$ is a non-negative integer we define
$\Lambda_i^s:A\to M_{sd_i}$ by
$$
\Lambda_i^s(f)=\text{diag}
(\underbrace{\Lambda_i(f),\Lambda_i(f),\dots,\Lambda_i(f)}_{s\ \text{times}}).
$$
Note that $\Lambda_i^{\frac n{d_i}}(f)=f(x_i)$ in $M_n$ for $f\in A$ and
$i=1,2,\dots,N$.

The following lemmas are left as exercises.

\begin{lemma}\label{repr}
Let $A=A(n,d_1,d_2,\dots,d_N)$ be a building block. The irreducible
representations (up to unitary equivalence) of $A$ are
$\Lambda_1,\Lambda_2,\dots,\Lambda_N$, together with point evaluations at
non-exceptional points.
\end{lemma}

\begin{lemma}\label{ideal}
Let $I$ be a closed two-sided ideal in $A$. There is a closed set
$F\subseteq \T$ such that
$$
I=\{f\in A : f(x)=0 \text{\ for all\ }x\in F\}.
$$
\end{lemma}

\begin{lemma}\label{trace}
Let $A=A(n,d_1,d_2,\dots,d_N)$ be a building block and let $\omega\in T(A)$.
There exists a Borel probability measure $\mu$ on $\T$ such that
$$
\omega(f)=\frac 1n \int_{\T} Tr(f(x))\, d\mu(x).
$$
It follows that $C_{\R}(\T)$ and $\text{Aff}\, T(A)$ are isomorphic as
order unit spaces via the map $f\mapsto \widehat{f\otimes 1}$,
$f\in C_{\R}(\T)$.
\end{lemma}

\begin{thm}\label{stable}
Let $A$ be a finite direct sum of building blocks. Then
$A$ is finitely generated and semiprojective.

\begin{proof}
First note that $A$ is a one-dimensional non-commutative  $CW$ complex, as
defined in \cite{ELPNCCW}. Hence $A$ is semiprojective by
\cite[Theorem 6.2.2]{ELPNCCW} and finitely generated by
\cite[Lemma 2.4.3]{ELPNCCW}.
\end{proof}
\end{thm}

Note that if $A=A(n,d_1,d_2,\dots,d_N)$ then $s(A)=\min(d_1,d_2,\dots,d_N)$.

Building blocks will sometimes be called circle building blocks in order to
distinguish them from interval building blocks. An interval building block
is a $C^*$-algebra $A$ of the form
$$
I(n,d_1,d_2,\dots,d_N)=
\{f\in C[0,1]\otimes M_n:f(x_i)\in M_{d_i},\ i=1,2,\dots,N\},
$$
where $0=x_1<x_2<\dots<x_N=1$ and $d_1,d_2,\dots,d_N$ are integers dividing
$n$. We will call $x_1,x_2,\dots,x_N$ the exceptional points of $A$.

\section{$K$-theory}
The purpose of this section is to calculate and interpret the $K$-theory
of a building block. We start out with the following lemma, which will be
used to calculate the $K_1$-group.

\begin{lemma}\label{coker}
Let $N\ge 2$ and let $a_1,a_2,\dots,a_N$ be positive integers. Define a
group homomorphism $\varphi:\Z^N\to \Z^N$ to be multiplication with the
$N\times N$ matrix
$$
C=
\begin{pmatrix}
\ \, \, a_1 & -a_2 &           \\
     &\ \, \,  a_2 & -a_3 &    \\
     &      &\ \, \,  a_3 & \ddots & \\
     &      &      & \ddots & -a_N \\
-a_1 &      &      &        &\ \, \,  a_N \\
\end{pmatrix}.
$$
For $k=1,2,\dots,N-1$, set
$$
s_k=lcm(a_1,a_2,\dots,a_k)
$$
and
$$
r_k=gcd(s_k,a_{k+1})=gcd(lcm(a_1,a_2,\dots,a_k),a_{k+1}).
$$
Choose integers $\alpha_k$ and $\beta_k$ such that
$$
r_k=\alpha_ks_k+\beta_ka_{k+1},\quad k=1,2,\dots,N-1.
$$
Then
$$
\text{coker}(\varphi)\cong
\Z\oplus\Z_{r_1}\oplus\Z_{r_2}\oplus\dots\oplus\Z_{r_{N-1}}.
$$
This isomorphism can be chosen such that for $k=1,2,\dots,N-2$, a generator
of the direct summand $\Z_{r_k}$ is mapped to the coset
$$
(\underbrace{0,0,\dots,0}_{k-1\ \text{times}},
1,-\frac{\beta_ka_{k+1}}{r_k},\underbrace{0,0,\dots,0}_{N-k-2\ \text{times}},
-\frac{\alpha_ks_k}{r_k})+ im(\varphi),
$$
such that a generator of the direct summand $\Z_{r_{N-1}}$ is mapped to the
coset
$$
(0,0,\dots,0,1,-1)+ im(\varphi),
$$
and such that a generator of the direct summand $\Z$ is mapped to the coset
$$
(0,0,\dots,0,1) + im(\varphi).
$$

\begin{proof}
Let $I_j$ denote the $j\times j$ identity matrix for any non-negative
integer $j$. For each $k=1,2,\dots,N-2$, define an integer matrix of size
$N\times N$ by
$$
A_k=\begin{pmatrix}
I_{k-1} \\
& 1 \\
& -\frac{\alpha_ks_k}{r_k} & 1 \\
& -\frac{\alpha_ks_k}{r_k} & & 1 \\
& \vdots & & & \ddots \\
& -\frac{\alpha_ks_k}{r_k} & & & & 1 \\
& 0 & & & & & 1 \\
\end{pmatrix}.
$$
Let $D_k$ denote the $2\times 2$ matrix
$$
\begin{pmatrix}
\ \, \, \alpha_k & \frac{a_{k+1}}{r_k} \\
-\beta_k & \frac{s_k}{r_k} \\
\end{pmatrix},
$$
and define for $k=1,2,\dots,N-1$, an integer matrix of size $N\times N$ by
$$
B_k=
\begin{pmatrix}
I_{k-1} \\
& D_k \\
& & I_{N-k-1} \\
\end{pmatrix}.
$$
For $k=0,1,2,\dots,N-2$, define yet another $N\times N$ matrix by
$$
X_k=
\begin{pmatrix}
r_1 \\
& r_2 \\
& & \ddots \\
& & & r_k \\
& & & & s_{k+1} & -a_{k+2} \\
& & & & s_{k+1} & & -a_{k+3} \\
& & & & \vdots & & &\ \ddots \\
& & & & s_{k+1} & & & & -a_N \\
& & & & 0 & & & & 0 \\
\end{pmatrix}.
$$
Finally, let $P$ be the $N\times N$ matrix
$$
\begin{pmatrix}
1 \\
1 & 1 \\
1 & 1 & 1 \\
\vdots & \vdots & \vdots & \ddots \\
1 & 1 & 1 & \hdots & 1 \\
\end{pmatrix}.
$$
Note that for $k=1,2,\dots,N-2$,
$$
s_{k+1}=lcm(s_k,a_{k+1})=\frac{s_ka_{k+1}}{r_k}.
$$
Using this, it is easily seen by induction that
$$
A_kA_{k-1}\cdots A_1 PCB_1B_2\cdots B_k = X_k,\quad k=0,1,2,\dots,N-2.
$$
It follows that
$$
A_{N-2}A_{N-1}\cdots A_1 PCB_1B_2\cdots B_{N-1} =
\begin{pmatrix}
r_1 \\
& r_2 \\
& & \ddots \\
& & & r_{N-1} \\
& & & & 0
\end{pmatrix}.
$$
Since all the matrices on the left-hand side, except $C$, are invertible in
$M_N(\Z)$, we obtain the desired calculation of $\text{coker}(\varphi)$.
Finally, it is easily verified that
$$
(A_{N-2}A_{N-1}\cdots A_1 P)^{-1} =
\begin{pmatrix}
1 \\
-\frac{\beta_1a_2}{r_1} & 1 \\
0 & -\frac{\beta_2 a_3}{r_2} & 1 \\
0 & 0 & -\frac{\beta_3 a_4}{r_3} & 1 \\
\vdots & \vdots & & \ddots & \ddots \\
0 & 0 & & & -\frac{\beta_{N-2} a_{N-1}}{r_{N-2}} & 1 \\
-\frac{\alpha_1s_1}{r_1} & -\frac{\alpha_2s_2}{r_2} & \hdots & \hdots &
-\frac{\alpha_{N-2}s_{N-2}}{r_{N-2}} & -1 & 1 \\
\end{pmatrix}.
$$
The last part of the lemma follows from this.
\end{proof}
\end{lemma}

Let $A=A(n,d_1,d_2,\dots,d_N)$ be a building block with exceptional
points $e^{2\pi i t_k}$, $k=1,2,\dots,N$, where $0<t_1<t_2<\dots<t_N<1$.
Set $t_{N+1}=t_1+1$. Define continuous functions $\omega_k:\T\to\T$ for
$k=1,2,\dots,N$, by
\begin{equation*}
\omega_k(e^{2\pi it})=
\begin{cases}
\exp(2\pi i\frac{t-t_k}{t_{k+1}-t_k}) &\ t_k\le t\le t_{k+1}, \\
1&\ t_{k+1}\le t\le t_k+1.
\end{cases}
\end{equation*}
Let $U_k^A$ be the unitary in $A$ defined by
$$
U_k^A(z)=\text{diag}(\omega_k(z),1,1,\dots,1),\quad z\in\T.
$$

\begin{thm}\label{K_1c}
Let $A=A(n,d_1,d_2,\dots,d_N)$ be a building block.
Set for $k=1,2,\dots,N-1$,
$$
s_k=lcm(\frac n{d_1},\frac n{d_2},\dots,\frac n{d_k}),
$$
and
$$
r_k=gcd(s_k,\frac n{d_{k+1}})=
gcd(lcm(\frac n{d_1},\frac n{d_2},\dots,\frac n{d_k}),\frac n{d_{k+1}}).
$$
Choose integers $\alpha_k$ and $\beta_k$ such that
$$
r_k=\alpha_ks_k+\beta_k\frac n{d_{k+1}},\quad k=1,2,\dots,N-1.
$$
Then
$$
K_1(A)\cong
\Z\oplus\Z_{r_1}\oplus \Z_{r_2}\oplus\dots\oplus\Z_{r_{N-1}}.
$$
This isomorphism can be chosen such that for $k=1,2,\dots,N-1$, a generator of
the direct summand $\Z_{r_k}$ is mapped to
$$
[U_k^A]-\frac{\beta_kn}{r_kd_{k+1}}[U_{k+1}^A]-\frac{\alpha_ks_k}{r_k}[U_N^A],
$$
and such that a generator of the direct summand $\Z$ is mapped to $[U_N^A]$.

\begin{proof}
Define a *-homomorphism $\pi : A\to M_{d_1}\oplus
M_{d_2}\oplus\dots\oplus M_{d_N}$ by
$$
\pi(f)=(\Lambda_1(f),\Lambda_2(f),\dots,\Lambda_N(f)).
$$
Via the identification $SM_n\cong \{f\in C[0,1]\otimes M_n:f(0)=f(1)=0\}$
we define a *-homomorphism $\iota : (SM_n)^N \to A$ by
$$
\iota(f_1,f_2,\dots,f_N)(e^{2\pi it})=f_k(\frac{t-t_k}{t_{k+1}-t_k}),
\qquad t_k\le t\le t_{k+1}.
$$
The short exact sequence
$$ 
  \begin{CD}
      0 @>  >> (SM_n)^N  @> \iota >> A @> \pi >>
M_{d_1}\oplus M_{d_2}\oplus\dots\oplus M_{d_N} @>  >> 0 \\   
  \end{CD}
$$
gives rise to a six-term exact sequence
$$
\begin{CD}
K_0((SM_n)^N) @> \iota_* >> K_0(A)
@> \pi_* >> K_0(M_{d_1}\oplus\dots\oplus M_{d_N})  \\
@AA A  @.  @VV \delta V \\
K_1(M_{d_1}\oplus\dots\oplus M_{d_N}) @<< \pi_* < K_1(A)
@<< \iota_* < K_1((SM_n)^N)
\end{CD}
$$
where $\delta$ denotes the exponential map.

By Bott periodicity $K_1((SM_n)^N)\cong\Z^N$ is generated by
$[V_1],[V_2],\dots,[V_N]$, where
$$
V_k(t)=(1,1,\dots,1,
\underbrace{\text{diag}(e^{2\pi it},1,\dots,1)}_{\text{coordinate } k},
1,1,\dots,1),\qquad t\in [0,1],
$$
is a unitary in $\widetilde{(SM_n)^N}$.
Note that ${\iota}(V_k-1)=U_k^A-1$ and hence $\iota_*([V_k])=[U_k^A]$ in
$K_1(A)$. Since the map $\iota_*:K_1((SM_n)^N)\to K_1(A)$ is surjective
it follows that $K_1(A)$ is generated by $[U_1^A],[U_2^A],\dots,[U_N^A]$,
and that $\iota_*$ gives rise to an isomorphism between the cokernel of
$\delta$ and $K_1(A)$.

Let $\{ e_{ij}^k\}$ denote the standard matrix units in
$M_{d_1}\oplus\dots\oplus M_{d_N}$.
Recall that $K_0(M_{d_1}\oplus\dots\oplus M_{d_N})\cong\Z^N$ is generated by
$[e_{11}^1],[e_{11}^2],\dots,[e_{11}^N]$. We leave it with the reader to check
that
$$
\delta([e_{11}^1])=-\frac n{d_1}[V_N]+\frac n{d_1}[V_1],
$$
and for $k=2,3,\dots,N$,
$$
\delta([e_{11}^k])=-\frac n{d_k}[V_{k-1}]+\frac n{d_k}[V_k].
$$
The conclusion follows from Lemma \ref{coker}.
\end{proof}
\end{thm}

Choose a continuous function $\gamma:\T\to\R$ such that
$$
Det(U_N^A(z))=z\exp(2\pi i\gamma(z)),\quad z\in\T.
$$
Define a unitary $v^A$ in $A$ by
$$
v^A(z)=U_N^A(z)\exp(-2\pi i\frac{\gamma(z)}n),\quad z\in\T.
$$
Note that $Det(v^A(z))=z$, $z\in\T$.

\begin{lemma}\label{K1det}
Let $A=A(n,d_1,d_2,\dots,d_N)$ be a building block and let $u\in A$
be a unitary. If
\begin{gather*}
Det(\Lambda_k(u))=1,\quad k=1,2,\dots,N, \\
Det(u(z))=1,\quad z\in\T,
\end{gather*}
then $u$ can be connected to $1$ via a continuous path of unitaries in $A$.

\begin{proof}
Let us start with a simple and well-known observation. Let $v$ be a unitary
in the $C^*$-algebra $B=\{f\in C[0,1]\otimes M_n:f(0)=f(1)\}$ such that the
winding number of $Det(v(\cdot))$ is $0$. Then $v$ can be connected to $1$
via a continuous path $(v_t)_{t\in [0,1]}$ in $U(B)$. If $v(0)=1$ we may
assume that $v_t(0)=1$ for every $t\in [0,1]$.

Let $e^{2\pi it_1},\dots, e^{2\pi it_N}$ be the exceptional
points of $A$, where $t_1<t_2<\dots<t_N$ are numbers in $]0,1[$.
Set $t_0=t_N-1$, $t_{N+1}=t_1+1$ and let $\iota_k:M_{d_k}\to M_n$ be the
inclusion, $k=1,2,\dots,N$. Since the group of unitaries in $M_{d_k}$ with
determinant $1$ is path-connected there exists a continuous function
$\gamma_k:[t_{k-1},t_{k+1}]\to U(M_{d_k})$ such that
$\gamma_k(t_{k-1})=\gamma_k(t_{k+1})=1$, $\gamma_k(t_k)=\Lambda_k(u)$, and
$Det(\gamma_k(\cdot))=1$. Set
$$
w_k(e^{2\pi it})=
\begin{cases}
\iota_k(\gamma_k(t))\quad &t\in [t_{k-1},t_{k+1}], \\
1\quad &t\in [t_{k+1},t_{k-1}+1].
\end{cases}
$$
It follows from the above observation that $w_k$ can be connected to $1$
via a continuous path of unitaries in $A$. Upon replacing $u$
with $uw_1^*w_2^*\dots w_N^*$ we may
thus assume that $u(e^{2\pi it_k})=1$ for $k=1,2,\dots,N$. Set
$$
y_k(e^{2\pi it})=
\begin{cases}
u(e^{2\pi it}) &\ t\in [t_k,t_{k+1}], \\
1    &\ t\in [t_{k+1},t_k+1].
\end{cases}
$$
Then $u=y_1y_2\dots y_N$. Again by the above observation, $y_k$ can be
connected to $1$ within $U(A)$ for $k=1,2,\dots,N$.
\end{proof}
\end{lemma}

Let $A=A(n,d_1,d_2,\dots,d_N)$ be a building block and set
$d=gcd(d_1,d_2,\dots,d_N)$. Since $d$ divides $d_i$ for every $i=1,2,\dots,N$,
there exists a unital and injective *-homomorphism $M_d\to A$ given by
$f\mapsto\text{diag}(f,f,\dots,f)$.

\begin{lemma}\label{cd}
Let $p$ be a projection in $A=A(n,d_1,d_2,\dots,d_N)$.  Then $p$ is
unitarily equivalent to a projection in $M_d\subseteq A$.

\begin{proof}
Let $r\in\Z$ denote the rank of $p$ and let
$e^{2\pi it_1},e^{2\pi it_2},\dots, e^{2\pi it_N}$ be the exceptional
points of $A$, where $0<t_1<t_2<\dots<t_N<1$. Since
$\frac n{d_k}$ divides $r$ for $k=1,2,\dots,N$, it follows
that $\frac nd$ also divides $r$. Hence there is a projection
$e\in M_d\subseteq A$ with the same trace as $p$.

For each $t\in [0,1]$ there is a unitary $u_t\in M_n$ such that
$$
e=u_tp(e^{2\pi it}){u_t}^*.
$$
We may assume that $u_{t_k}\in M_{d_k}$, $k=1,2,\dots,N$, and that $u_0=u_1$.
By compactness
$$
[0,1]=\bigcup_{j=1}^{L-1} [s_j,s_{j+1}],
$$
where $0=s_1<s_2<\dots<s_L=1$,
$\{t_1,t_2,\dots,t_N\}\subseteq \{s_1,s_2,\dots,s_L\}$, and
$$
t\in [s_j,s_{j+1}]\ \Rightarrow\ \|u_{s_j}\, p(e^{2\pi it})\, u_{s_j}^*-e\|<1.
$$
Set $z_j(t)=v_j(t)|v_j(t)|^{-1}$ for $t\in [s_j,s_{j+1}]$,
$j=1,2,\dots,L-1$, where
$$
v_j(t)=1-u_{s_j}\, p(e^{2\pi it})\, u_{s_j}^*-e+
2\, e\, u_{s_j}\, p(e^{2\pi it})\, u_{s_j}^*.
$$
Then $t\mapsto z_j(t)$, $t\in [s_j,s_{j+1}]$, is a continuous path of
unitaries in $M_n$, and by \cite[Lemma 6.2.1]{MURPHY}
$$
e=z_j(t)\, u_{s_j}\, p(e^{2\pi it})\, u_{s_j}^*\, z_j(t)^*,
\qquad t\in [s_j,s_{j+1}].
$$
As $U(M_n)\cap \{ e \} '$ is path-connected there is for each
$k=1,2,\dots,L-1$ a continuous map
$\gamma_j : [s_j,s_{j+1}]\to U(M_n)\cap \{ e \} '$ such that
$$
\gamma_j(s_j)=1,\qquad
\gamma_j(s_{j+1})=u_{s_{j+1}}\, u_{s_j}^*\, z_j(s_{j+1})^*.
$$
Since $z_j(s_j)=1$ for $j=1,2,\dots,L-1$, we can define a unitary $u\in A$ by
$$
u(e^{2\pi it})=\gamma_j(t)z_j(t)u_{s_j},\qquad t\in [s_j,s_{j+1}].
$$
Then $upu^*=e$.
\end{proof}
\end{lemma}

\begin{cor}\label{cutdown}
If $p\in A=A(n,d_1,d_2,\dots,d_N)$ is a projection of rank $r\neq 0$ then
$$
pAp\cong A(r,\frac rnd_1,\frac rnd_2,\dots,\frac rnd_N).
$$
\end{cor}

\begin{cor}\label{K_0c}
The embedding $M_d\subseteq A$ gives rise to an isomorphism of ordered groups
with order units between $K_0(M_d)$ and $K_0(A)$. In other words,
$$
\bigl(K_0(A),K_0(A)^+,[1]\bigr)\cong\bigl(\Z,\Z^+,d\bigr).
$$
\end{cor}

By Lemma \ref{trace} we have the following:
\begin{cor}\label{imK_0}
If $A=A(n,d_1,d_2,\dots,d_N)$ then $\rho_A(K_0(A))=\Z\frac 1d\widehat{1}$
in $\text{Aff}\, T(A)$.
\end{cor}

\begin{lemma}\label{proj}
$A=A(n,d_1,d_2,\dots,d_N)$ is unital projectionless if and only if $d=1$.
\begin{proof}
As in the proof of Lemma \ref{cd} we see that there exists a projection
$p\in A$ of rank $r\le n$ if and only if $\frac nd$ divides $r$. The conclusion
follows.
\end{proof}
\end{lemma}

\begin{lemma}\label{K_1ex}
Let $K$ be a positive integer and let $H$ be a finite abelian
group. There exists a unital projectionless building block $A$ with
$s(A)\ge K$ such that $K_1(A)\cong \Z\oplus H$.
\begin{proof}
Let
$$
H\cong \Z_{p_1^{k_1}}\oplus\Z_{p_2^{k_2}}\oplus\dots\oplus\Z_{p_m^{k_m}},
$$
where $m$ is a positive integer, $k_1,\dots,k_m$ are non-negative integers,
and $p_1,\dots,p_m$ are prime numbers. Let $q_1,q_2,\dots,q_{m+1}\ge K$
be prime numbers, mutually different as well as different from
$p_1,p_2,\dots,p_m$. Define integers $n$ and $d_1,d_2,\dots,d_{m+1}$ by
\begin{align*}
n&=p_1^{k_1}\, p_2^{k_2}\dots p_m^{k_m}\, q_1\, q_2\dots q_{m+1}, \\
d_1&=q_2\, q_3\dots q_{m+1}, \\
d_i&=\frac{p_1^{k_1}p_2^{k_2}\dots p_m^{k_m}}{p_{i-1}^{k_{i-1}}}\,
\frac{q_1q_2\dots q_{m+1}}{q_i},\quad 2\le i\le m+1.
\end{align*}
Set $A=A(n,d_1,d_2,\dots,d_{m+1})$. Then $K_1(A)\cong \Z\oplus H$ by
Theorem \ref{K_1c}. $A$ is unital projectionless by Lemma \ref{proj}.
\end{proof}
\end{lemma}

\section{$KK$-theory}

Recall a few facts about $KK$-theory that can be found in e.g \cite{BKOA}.
$KK$ is a homotopy invariant bifunctor
from the category of $C^*$-algebras to the category of abelian groups that is
contravariant in the first variable and covariant in the second.
A *-homomorphism $\varphi:A\to M_n(B)$ defines an element
$[\varphi]\in KK(A,B)$. We have an associative map $KK(B,C)\times KK(A,B)\to
KK(A,C)$, the Kasparov product, that generalizes composition of
*-homomorphisms.

The purpose of this section is to analyze the $KK$-theory of our building
blocks. Inspired by the work of Jiang and Su \cite[section 3]{JSSP},
we will consider the $K$-homology groups $K^0(A)=KK(A,\C)$.
A *-homomorphism $\varphi:A\to M_n(B)$ induces a group
homomorphism $\varphi^*:K^0(B)\to K^0(A)$ via the Kasparov product.
$K^0(M_n)\cong\Z$ is generated by the class of the identity map on $M_n$.

If $A$ and $B$ are unital $C^*$-algebras we let $KK(A,B)_e$ be the set of
elements $\kappa\in KK(A,B)$ such that $\kappa_*:K_0(A)\to K_0(B)$ preserves
the order unit.

\begin{lemma}\label{K^0i}
Let 
$$
A=\{f\in C[0,1]\otimes M_n:f(t_i)\in M_{d_i},\ i=1,2,\dots,N\}
$$
where $N\ge 2$, $0\le t_1<t_2<\dots<t_N\le 1$, and $d_1,d_2,\dots,d_N$
are integers dividing $n$. Let $\Lambda_i:A\to M_{d_i}$ be the *-homomorphism
induced by evaluation at $t_i$, $i=1,2,\dots,N$.
Then $K^0(A)$ is generated by $[\Lambda_1],[\Lambda_2],\dots,[\Lambda_N]$.
Furthermore, for $a_1,a_2,\dots,a_N\in\Z$ we have that
$$
a_1[\Lambda_1]+a_2[\Lambda_2]+\dots +a_N[\Lambda_N]=0
$$
if and only if there exist $b_1,b_2,\dots,b_N\in\Z$ such that
$\sum_{i=1}^N b_i=0$ and
$$
a_i=b_i\frac{n}{d_i},\quad i=1,2,\dots,N.
$$
\begin{proof}
Choose $y\in ]0,1[$ such that $t_1<y<t_2$. Set
\begin{gather*}
B=\{f\in C[0,y]\otimes M_n : f(t_1)\in M_{d_1}\}, \\
C=\{f\in C[y,1]\otimes M_n : f(t_i)\in M_{d_i}, i=2,3,\dots,N\}.
\end{gather*}
We have a pull-back diagram
$$ 
  \begin{CD}
      A @> g_1 >> B    \\   
      @Vg_2VV      @VVf_1V \\
      C   @>> f_2 > M_n
  \end{CD}
$$
where $g_1,g_2$ are the restriction maps and $f_1,f_2$ evaluation at $y$.
Apply the Mayer-Vietoris sequence \cite[Theorem 21.5.1]{BKOA} to get a
six-term exact sequence
$$
\begin{CD}
K^0(M_n) @> (-f_1^*,f_2^*) >> K^0(B)\oplus K^0(C)
@> g_1^*+g_2^* >> K^0(A)  \\
@AA A  @.  @VV V \\
K^1(A) @<< g_1^*+g_2^*< K^1(B)\oplus K^1(C)
@<<(-f_1^*,f_2^*)< K^1(M_n).
\end{CD}
$$
Note that $K^1(M_n)=0$ and $K^0(M_n)\cong\Z$. Thus the exact sequence becomes
$$
\begin{CD}
\Z@> \varphi >> K^0(B)\oplus K^0(C) @> \psi >> K^0(A)  \\
@AA A  @.  @VV V \\
K^1(A) @<<  <K^1(B)\oplus K^1(C)  @<< < 0.
\end{CD}
$$
Since $f_1$ is homotopic to evaluation at $x_1$ in $B$ and $f_2$ is homotopic
to evaluation at $x_2$ in $C$ we see that
$$
\varphi(k)=(-k\frac{n}{d_1}[\Lambda_1|_B],k\frac{n}{d_N}[\Lambda_N|_C]),
\quad k\in\Z.
$$
$B$ is homotopic to $M_{d_1}$ via $\Lambda_1|_B$ and hence
$K^0(B)\cong\Z$ is generated by $[\Lambda_1|_B]$.

For $N=2$ we have that $K^0(C)$ is generated by $[\Lambda_2|_C]$ and that
$$
\psi(a_1[\Lambda_1|_B],a_2[\Lambda_2|_C])=a_1[\Lambda_1]+a_2[\Lambda_2].
$$
Thus $K^0(A)$ is generated by $[\Lambda_1]$ and $[\Lambda_2]$ and
$$
a_1[\Lambda_1]+a_2[\Lambda_2]=0\Leftrightarrow\exists b_1\in\Z:
a_1=-b_1\frac{n}{d_1},\ a_2=b_1\frac{n}{d_2}.
$$

Proceeding by induction, assume that the lemma holds for $N-1$. By the
induction hypothesis $K^0(C)$ is generated by
$[\Lambda_2|_C],[\Lambda_3|_C],\dots,[\Lambda_N|_C]$. Note that
$$
\psi\bigl(a_1[\Lambda_1|_B],
(a_2[\Lambda_2|_C]+\dots+a_N[\Lambda_N|_C])\bigr)=
\sum_{i=1}^N a_i[\Lambda_i],
$$
such that $A$ is generated by $[\Lambda_1],\dots,[\Lambda_N]$. It also follows
that
$$
a_1[\Lambda_1]+a_2[\Lambda_2]+\dots+a_N[\Lambda_N]=0
$$
if and only if there exists $k\in\Z$ such that
$$
-k\frac{n}{d_1}[\Lambda_1|_B]=a_1[\Lambda_1|_B],\quad
k\frac{n}{d_N}[\Lambda_N|_C]=a_2[\Lambda_2|_C]+\dots+a_N[\Lambda_N|_C].
$$
By the induction hypothesis this happens if and only if there exist
$k,c_2,\dots,c_N\in\Z$ such that $\sum_{i=2}^N c_i=0$ and
$$
a_1=-k\frac n{d_1},\qquad a_i=c_i\frac{n}{d_i},\quad i=2,3,\dots,N-1,\qquad
a_N-k\frac{n}{d_N}=c_N\frac{n}{d_N}.
$$
The desired conclusion follows easily from these equations.
\end{proof}
\end{lemma}

\begin{prop}\label{K^0c}
Let $A=A(n,d_1,d_2,\dots,d_N)$ be a building block.
Then $K^0(A)$ is generated by $[\Lambda_1],[\Lambda_2],\dots,[\Lambda_N]$.
Furthermore, for $a_1,a_2,\dots,a_N\in\Z$ we have that
$$
a_1[\Lambda_1]+a_2[\Lambda_2]+\dots +a_N[\Lambda_N]=0
$$
if and only if there exist $b_1,b_2,\dots,b_N\in\Z$ such that
$\sum_{i=1}^N b_i=0$ and
$$
a_i=b_i\frac{n}{d_i},\quad i=1,2,\dots,N.
$$

\begin{proof}
Choose $t_1,t_2,\dots,t_N\in ]0,1[$ such that $e^{2\pi it_k}$,
$k=1,2,\dots,N$, are the exceptional points for $A$. Set
$$
B=\{f\in C[0,1]\otimes M_n:f(t_k)\in M_{d_k},\ k=1,2,\dots,N\}.
$$
Define a *-homomorphism $\iota:A\to B$ by $\iota(f)(t)=f(e^{2\pi it})$.
Let $\pi : A\to M_n$ be evaluation at $1\in\T$. Let
$\alpha : M_n\to M_n\oplus M_n$ denote the map $\alpha(x)=(x,x)$.
Let $\beta: B\to M_n\oplus M_n$ be the map $\beta(f)=(f(0),f(1))$.
We have a pull-back diagram
$$ 
  \begin{CD}
      A @> \pi >> M_n    \\   
      @V\iota VV      @VV \alpha V \\
      B   @>> \beta > M_n\oplus M_n
  \end{CD}
$$
and hence by \cite[Theorem 21.5.1]{BKOA} a six-term exact sequence of the form
$$
\begin{CD}
K^0(M_n\oplus M_n) @> (-\alpha^*,\beta^*) >> K^0(M_n)\oplus K^0(B)
@> \pi^*+\iota^* >> K^0(A)  \\
@AA A  @.  @VV V \\
K^1(A) @<<\pi^*+\iota^* < K^1(M_n)\oplus K^1(B)
@<< (-\alpha^*,\beta^*) < K^1(M_n).
\end{CD}
$$
$K^0(M_n\oplus M_n)\cong\Z\oplus\Z$ is generated by $[\pi_1]$ and $[\pi_2]$
where $\pi_1,\pi_2:M_n\oplus M_n\to M_n$ are the coordinate projections.
$K^0(M_n)\cong\Z$ is generated by the class of the identity map $id$ on $M_n$.
Note that
\begin{gather*}
\pi^*([id])=\frac{n}{d_1}[\Lambda_1^A], \\
\iota^*([\Lambda_i^B])=[\Lambda_i^A],\qquad i=1,2,\dots,N, \\
(-\alpha^*,\beta^*)(a[\pi_1]+b[\pi_2])=
(-(a+b)[id],(a+b)\frac{n}{d_1}[\Lambda_1^B]).
\end{gather*}
As $\pi^*+\iota^*$ maps onto $K^0(A)$ (because $K^1(M_n)=0$) and as
$im(\pi^*)\subseteq im(\iota^*)$,
we see that $\iota^*$ is surjective. Assume that $\iota^*(x)=0$. Then
$(0,x)\in im(-\alpha^*,\beta^*)$ and hence $x=0$ by the above.
Thus $\iota^*$ is an isomorphism and the conclusion follows from
Lemma \ref{K^0i}.
\end{proof}
\end{prop}

\begin{prop}\label{sf}
Let $A=A(n,d_1,d_2,\dots,d_N)$ and $B=A(m,e_1,e_2,\dots,e_M)$ be building
blocks and let $h:K^0(B)\to K^0(A)$ be a group homomorphism. For every
$j=1,2,\dots,M$, $i=1,2,\dots,N$, there is a uniquely determined integer
$h_{ji}$, with $0\le h_{ji}<\frac{n}{d_i}$ for $i\neq N$, such that
$$
\begin{pmatrix}
h([\Lambda_1^B]) \\
h([\Lambda_2^B]) \\
\vdots \\
h([\Lambda_M^B])
\end{pmatrix}
=
\begin{pmatrix}
h_{11} & h_{12} & \hdots & h_{1N} \\
h_{21} & h_{22} & \hdots & h_{2N} \\
\vdots & \vdots &        & \vdots \\
h_{M1} & h_{M2} & \hdots & h_{MN}
\end{pmatrix}
\begin{pmatrix}
{[\Lambda_1^A]} \\
{[\Lambda_2^A]} \\
\vdots \\
{[\Lambda_N^A]}
\end{pmatrix}.
$$
This will be called the standard form for $h$.

The integers determined by $h$ above satisfy the equations
\begin{gather*}
\frac{m}{e_j}h_{ji} \equiv \frac{m}{e_M}h_{Mi}
\text{\ \ mod\ \ } \frac{n}{d_i},
\quad j=1,2,\dots,M,\ i=1,2,\dots,N, \\
\frac{m}{e_j}\sum_{i=1}^N h_{ji}d_i = \frac{m}{e_M}\sum_{i=1}^N h_{Mi}d_i,
\quad j=1,2,\dots,M.
\end{gather*}

\begin{proof}
By Proposition \ref{K^0c}, or simply because homotopic *-homomorphisms
$A\to M_n$ define the same elements in $K^0(A)$, we have that
$$
\frac{n}{d_N}[\Lambda_N^A]=\frac{n}{d_i}[\Lambda_i^A],\quad i=1,2,\dots,N.
$$
From this the existence follows.

To check uniqueness, assume
$$
\begin{pmatrix}
h_{11} & h_{12} & \hdots & h_{1N} \\
h_{21} & h_{22} & \hdots & h_{2N} \\
\vdots & \vdots &        & \vdots \\
h_{M1} & h_{M2} & \hdots & h_{MN}
\end{pmatrix}
\begin{pmatrix}
{[\Lambda_1^A]} \\
{[\Lambda_2^A]} \\
\vdots \\
{[\Lambda_N^A]}
\end{pmatrix} =0
$$
where
$$
-\frac{n}{d_i}<h_{ji}<\frac{n}{d_i},\quad i=1,2,\dots,N-1,\ j=1,2,\dots,M.
$$
Fix some $j=1,2,\dots,M$. By Proposition \ref{K^0c} there exist integers
$b_{ji}$ such that $h_{ji}=b_{ji}\frac{n}{d_i}$, $i=1,2,\dots,N$. Therefore
$h_{ji}=0$ for $i=1,2,\dots,N$.

Finally, to prove the equations above, fix again some $j=1,2,\dots,M$.
Note that
$$
0=h(0)=h(-\frac{m}{e_j}[\Lambda_j^B]+\frac{m}{e_M}[\Lambda_M^B])
=\sum_{i=1}^N(-\frac{m}{e_j}h_{ji}+\frac{m}{e_M}h_{Mi})[\Lambda_i^A].
$$
Hence there exist integers $b_{ji}$, $i=1,2,\dots,N$, such that
$\sum_{i=1}^N b_{ji}=0$ and
$$
-\frac{m}{e_j}h_{ji}+\frac{m}{e_M}h_{Mi}=b_{ji}\frac{n}{d_i}.
$$
The desired conclusion follows easily from these equations.
\end{proof}
\end{prop}

From now on, let $A=A(n,d_1,d_2,\dots,d_N)$ and $B=A(m,e_1,e_2,\dots,e_M)$
be building blocks. Define a group homomorphism
$$
\Gamma : KK(A,B)\to Hom(K^0(B),K^0(A))\oplus K_1(B)
$$
by
$$
\Gamma(\kappa)=(\kappa^*,\kappa_*[v^A]).
$$
We want to show that $\Gamma$ is an isomorphism in certain cases.

\begin{prop}\label{homex}
Let $h:K^0(B)\to K^0(A)$ be a group homomorphism with standard form
$$
\begin{pmatrix}
h([\Lambda_1^B]) \\
h([\Lambda_2^B]) \\
\vdots \\
h([\Lambda_M^B])
\end{pmatrix}
=
\begin{pmatrix}
h_{11} & h_{12} & \hdots & h_{1N} \\
h_{21} & h_{22} & \hdots & h_{2N} \\
\vdots & \vdots &        & \vdots \\
h_{M1} & h_{M2} & \hdots & h_{MN}
\end{pmatrix}
\begin{pmatrix}
{[\Lambda_1^A]} \\
{[\Lambda_2^A]} \\
\vdots \\
{[\Lambda_N^A]}
\end{pmatrix}
$$
where $h_{jN}\ge\frac{n}{d_N}$ for $j=1,2,\dots,M$, and
$\sum_{i=1}^N h_{Mi}d_i=e_M$.
Let $\chi\in K_1(B)$. There is a unital *-homomorphism
$\varphi:A\to B$ such that $\Gamma([\varphi])=(h,\chi)$.

\begin{proof}
Let $1\le i\le N$. By Proposition \ref{sf} there is an integer $s_i$,
$0\le s_i<\frac{n}{d_i}$, and integers $l_{ji}$, $j=1,2,\dots,M$, such that
\begin{equation}\label{s_i}
\frac{m}{e_j}h_{ji}=l_{ji}\frac{n}{d_i}+s_i.
\end{equation}
Note that $l_{ji}\ge 0$ for $i=1,2,\dots,N-1$, and
$l_{jN}\ge 1$. By Proposition \ref{sf} we see that for $j=1,2,\dots,M$,
$$
m=\frac m{e_M} \sum_{i=1}^N h_{Mi}d_i=
\frac m{e_j} \sum_{i=1}^N h_{ji}d_i=\sum_{i=1}^N (l_{ji}n+s_id_i).
$$
By (\ref{s_i}) there exists a unitary $V_j\in M_m$ such that the matrix
$$
V_j\text{diag}\bigl(\Lambda_1^{s_1}(f),\dots,\Lambda_N^{s_N}(f),
\underbrace{f(x_1),\dots,f(x_1)}_{l_{j1}\text{\ times}},\dots,
\underbrace{f(x_N),\dots,f(x_N)}_{l_{jN}\text{\ times}}\bigr)V_j^*
$$
belongs to $M_{e_j}\subseteq M_m$ for all $f\in A$.

Set
$$
L=\frac 1n(m-\sum_{i=1}^N s_id_i)=\sum_{i=1}^N l_{ji},\quad j=1,2,\dots,M.
$$
Let $x_1,x_2,\dots,x_N$ denote the exceptional points of $A$ and let
$y_1,y_2,\dots,y_M$ be those of $B$.
Choose continuous functions $\lambda_1,\lambda_2,\dots,\lambda_{L-1}:\T\to\T$
such that
$$
\bigl(\lambda_1(y_j),\lambda_2(y_j),\dots,\lambda_{L-1}(y_j)\bigr)=
\bigl(\underbrace{x_1,\dots,x_1}_{l_{j1} \text{\ times}},
\dots,
\underbrace{x_{N-1},\dots,x_{N-1}}_{l_{j(N-1)} \text{\ times}},
\underbrace{x_N,\dots,x_N}_{l_{jN}-1\text{\ times}}\bigr)
$$
as ordered tuples. Choose a unitary $U\in C(\T)\otimes M_m$ such that
$U(y_j)=V_j$. Define a unital *-homomorphism $\psi:A\to B$ by
$$
\psi(f)(z)=U(z)\text{diag}\bigl(\Lambda_1^{s_1}(f),\dots,\Lambda_N^{s_N}(f),
f(\lambda_1(z)),\dots,f(\lambda_{L-1}(z)),f(x_N)\bigr)U(z)^*.
$$
By Theorem \ref{K_1c} we have that $\chi=\sum_{j=1}^M a_j[U_j^B]$ for some
$a_1,a_2,\dots,a_M\in\Z$. Let
$$
\psi_*[v^A]=\sum_{i=1}^M b_j[U_j^B]
$$
in $K_1(B)$. Define $\xi:\T\to\T$ by
\begin{equation}\label{xideteq}
\xi(z)=\prod_{j=1}^M Det(U_j^B(z))^{a_j-b_j},
\end{equation}
and define $\lambda_L:\T\to\T$ by $\lambda_L(z)=\xi(z)x_N$.
Note that $\lambda_L(y_j)=x_N$, $j=1,2,\dots,M$. Define $\varphi:A\to B$ by
$$
\varphi(f)(z)=U(z)\text{diag}\bigl(\Lambda_1^{s_1}(f),\dots,\Lambda_N^{s_N}(f),
f(\lambda_1(z)),\dots,f(\lambda_{L}(z))\bigr)U(z)^*.
$$
By Lemma \ref{K1det} and (\ref{xideteq}) we see that in $K_1(B)$,
\begin{align*}
\varphi_*[v^A]=&\ \psi_*[v^A]+
[z\mapsto U(z)\text{diag}(1,1,\dots,1,v^A(\lambda_L(z))v^A(x_N)^*)U(z)^*] \\
=&\ \psi_*[v^A]+\sum_{j=1}^M (a_j-b_j)[U_j^B]=\sum_{j=1}^M a_j[U_j^B].
\end{align*}
Since $\varphi(f)(y_j)=\psi(f)(y_j)$, $f\in A$, $j=1,2,\dots,M$, we conclude
that
$$
\varphi^*([\Lambda_j^B])=[\Lambda_j^B\circ\varphi]=[\Lambda_j^B\circ\psi]=
\sum_{i=1}^N (s_i+l_{ji}\frac{n}{d_i})\frac{e_j}{m}[\Lambda_i^A]=
\sum_{i=1}^N h_{ji}[\Lambda_i^A]=h([\Lambda_j^B]).
$$
\end{proof}
\end{prop}

\begin{lemma}\label{K^0surj}
Let $h:K^0(B)\to K^0(A)$ be a group homomorphism and assume that there exists
a homomorphism $h': K^0(B)\to K^0(A)$ with standard form
$$
\begin{pmatrix}
h'([\Lambda_1^B]) \\
h'([\Lambda_2^B]) \\
\vdots \\
h'([\Lambda_M^B])
\end{pmatrix}
=
\begin{pmatrix}
h'_{11} & h'_{12} & \hdots & h'_{1N} \\
h'_{21} & h'_{22} & \hdots & h'_{2N} \\
\vdots  & \vdots  &        & \vdots  \\
h'_{M1} & h'_{M2} & \hdots & h'_{MN}
\end{pmatrix}
\begin{pmatrix}
{[\Lambda_1^A]} \\
{[\Lambda_2^A]} \\
\vdots \\
{[\Lambda_N^A]}
\end{pmatrix}
$$
where $h_{jN}'\ge\frac{n}{d_N}$ for $j=1,2,\dots,M$, and
$\sum_{i=1}^N h_{Mi}'d_i=e_M$. Then there is a $\kappa\in KK(A,B)$ such that
$\kappa^*=h$ in $Hom(K^0(B),K^0(A))$.

\begin{proof}
By Proposition \ref{homex} there exists an element $\nu\in KK(A,B)$ such that
$\nu^*=h'$. Let $h\in Hom(K^0(B),K^0(A))$ have standard form
$$
\begin{pmatrix}
h([\Lambda_1^B]) \\
h([\Lambda_2^B]) \\
\vdots \\
h([\Lambda_M^B])
\end{pmatrix}
=
\begin{pmatrix}
h_{11} & h_{12} & \hdots & h_{1N} \\
h_{21} & h_{22} & \hdots & h_{2N} \\
\vdots & \vdots &        & \vdots \\
h_{M1} & h_{M2} & \hdots & h_{MN}
\end{pmatrix}
\begin{pmatrix}
{[\Lambda_1^A]} \\
{[\Lambda_2^A]} \\
\vdots \\
{[\Lambda_N^A]}
\end{pmatrix}.
$$
By adding an integer-multiple of $h'$ we may assume that $h_{jN}\ge 0$ for
$j=1,2,\dots, M$.
Define $l_{ji}$ and $s_i$, $i=1,2,\dots, N$,
as in the proof of Proposition \ref{homex}. Let
$$
c=\frac m{e_M} \sum_{i=1}^N h_{Mi}d_i=\frac m{e_j} \sum_{i=1}^N h_{ji}d_i
=\sum_{i=1}^N (l_{ji}n+s_id_i),\quad j=1,2,\dots,M.
$$
Choose a positive integer $d$ such that $c\le dm$.
Choose for each $j=1,2,\dots,M$, a unitary $V_j\in M_{dm}$ such that
the matrix
$$
V_j\text{diag}\bigl(\Lambda_1^{s_1}(f),\dots,\Lambda_N^{s_N}(f),
\underbrace{f(x_1),\dots,f(x_1)}_{l_{j1}\text{\ times}},\dots,
\underbrace{f(x_N),\dots,f(x_N)}_{l_{jN}\text{\ times}},
\underbrace{0,\dots,0}_{dm-c}
\bigl)V_j^*
$$
belongs to $M_{de_j}\subseteq M_{dm}$ for all $f\in A$.

As in the proof of Proposition \ref{homex} these matrices can be connected
to define a *-homomorphism
$\varphi:A\to M_d(B)$. We leave it with the reader to check that
$\varphi^*=h$ on $K^0(B)$. Set $\kappa=[\varphi]$.
\end{proof}
\end{lemma}

\begin{prop}\label{KKiso}
Assume that there exists a homomorphism $h': K^0(B)\to K^0(A)$ with
standard form
$$
\begin{pmatrix}
h'([\Lambda_1^B]) \\
h'([\Lambda_2^B]) \\
\vdots \\
h'([\Lambda_M^B])
\end{pmatrix}
=
\begin{pmatrix}
h'_{11} & h'_{12} & \hdots & h'_{1N} \\
h'_{21} & h'_{22} & \hdots & h'_{2N} \\
 \vdots &  \vdots &        & \vdots  \\
h'_{M1} & h'_{M2} & \hdots & h'_{MN}
\end{pmatrix}
\begin{pmatrix}
{[\Lambda_1^A]} \\
{[\Lambda_2^A]} \\
\vdots \\
{[\Lambda_N^A]}
\end{pmatrix}
$$
where $h_{jN}'\ge\frac{n}{d_N}$ for $j=1,2,\dots,M$, and
$\sum_{i=1}^N h_{Mi}'d_i=e_M$.
Then the map $\Gamma:KK(A,B)\to Hom(K^0(B),K^0(A))\oplus K_1(B)$ is an
isomorphism.

\begin{proof}
By Theorem \ref{K_1c} there exist finite abelian groups $G$ and $H$ such that
$K_1(A)\cong\Z\oplus G$, $K_1(B)\cong\Z\oplus H$. By the 
universal coefficient theorem, \cite[Theorem  1.17]{RSUCT},
\begin{align*}
KK(A,B)\ \cong\ &Ext(K_0(A),K_1(B))\oplus Ext(K_1(A),K_0(B))\oplus \\
&Hom(K_0(A),K_0(B))\oplus Hom(K_1(A),K_1(B)) \\
\ \cong\
&0\oplus G\oplus \Z\oplus Hom(G,H)\oplus K_1(B).
\end{align*}
By the universal coefficient theorem again, $K^0(A)\cong K_1(A)$ and
$K^0(B)\cong K_1(B)$. Hence
$$
Hom(K^0(B),K^0(A))\oplus K_1(B)\cong K_1(A)\oplus Hom(H,G)\oplus K_1(B).
$$
Note that $Hom(G,H)\cong Hom(H,G)$. Thus $Hom(K^0(B),K^0(A))\oplus K_1(B)$
and $KK(A,B)$ are isomorphic groups.
Since any surjective endomorphism of a finitely generated abelian group
is an isomorphism, it suffices to show that $\Gamma$ is surjective.

Let $(h,\chi)\in Hom(K^0(B),K^0(A))\oplus K_1(B)$.
By Lemma \ref{K^0surj} there exists an element $\kappa\in KK(A,B)$
such that $\Gamma(\kappa)=(h-h',\eta)$ for some $\eta\in K_1(B)$.
Next, by Proposition \ref{homex} there exists a $\nu\in KK(A,B)$ such that
$\Gamma(\nu)=(h',\chi-\eta)$. Thus $\Gamma(\kappa+\nu)=(h,\chi)$.
\end{proof}
\end{prop}

\begin{thm}\label{KKliftc}
Let $A=A(n,d_1,d_2,\dots,d_N)$ and $B=A(m,e_1,e_2,\dots,e_M)$ be building
blocks such that $s(B)\ge Nn$ and assume that there exists an element
$\kappa$ in $KK(A,B)_e$.
Then the map $\Gamma:KK(A,B)\to Hom(K^0(B),K^0(A))\oplus K_1(B)$ is an
isomorphism and there exists a unital *-homomorphism $\varphi:A\to B$ such that
$[\varphi]=\kappa$.

\begin{proof}
Let $\kappa^*:K^0(B)\to K^0(A)$ have standard form
$$
\begin{pmatrix}
\kappa^*([\Lambda_1^B]) \\
\kappa^*([\Lambda_2^B]) \\
\vdots \\
\kappa^*([\Lambda_M^B])
\end{pmatrix}
=
\begin{pmatrix}
h_{11} & h_{12} & \hdots & h_{1N} \\
h_{21} & h_{22} & \hdots & h_{2N} \\
\vdots & \vdots &        & \vdots \\
h_{M1} & h_{M2} & \hdots & h_{MN}
\end{pmatrix}
\begin{pmatrix}
{[\Lambda_1^A]} \\
{[\Lambda_2^A]} \\
\vdots \\
{[\Lambda_N^A]}
\end{pmatrix}.
$$
Let $\cdot$ denote the Kasparov product. By assumption we have that
$[1_A]\cdot\kappa=[1_B]$ in $KK(\C,B)\cong K_0(B)$. Thus
$$
[1_B]\cdot[\Lambda_j^B] = [1_A]\cdot\kappa\cdot[\Lambda_j^B] =
[1_A]\cdot (\sum_{i=1}^N h_{ji}[\Lambda_i^A])
$$
in $KK(\C,\C)\cong\Z$. Hence $e_j=\sum_{i=1}^N h_{ji} d_i$ for $j=1,2,\dots,M$.
This implies that $h_{jN}>\frac{n}{d_N}$ since
$$
Nn\le e_j=\sum_{i=1}^N h_{ji} d_i < \sum_{i=1}^{N-1} \frac{n}{d_i}d_i+h_{jN}d_N
=(N-1)n+h_{jN}d_N.
$$
Therefore $\Gamma$ is an isomorphism by Proposition \ref{KKiso}.
By Proposition \ref{homex} there is a unital *-homomorphism
$\varphi:A\to B$ such that $\Gamma([\varphi])=\Gamma(\kappa)$. Thus
$[\varphi]=\kappa$.
\end{proof}
\end{thm}

\section{The commutator subgroup of the unitary group}

In this section we analyze the unitary group modulo the closure of its
commutator subgroup for building blocks.

\begin{lemma}\label{cano}
Let $A$ be a unital inductive limit of a sequence of finite direct sums
of building blocks. Then the canonical maps $\pi_0(U(A))\to K_1(A)$ and
$\pi_1(U(A))\to K_0(A)$ are isomorphisms.

\begin{proof}
Following \cite{TNS} we let $k_n(\cdot)=\pi_{n+1}(U(\cdot))$ for every integer
$n\ge -1$. By \cite[Proposition 2.6]{TNS} it suffices to show that the
canonical maps $k_{-1}(A)\to k_{-1}(A\otimes\mathcal{K})\cong K_1(A)$ and
$k_0(A)\to k_0(A\otimes\mathcal{K})\cong K_0(A)$ are isomorphisms,
where $\mathcal{K}$ denotes the set of compact
operators on a separable infinite dimensional Hilbert-space.
As noted in \cite{TNS} it follows from \cite[Proposition 4.4]{HKAF} that
$k_n$ is a continuous functor. Since it is obviously additive, we may assume
that $A$ is a building block.

As in the proof of Theorem \ref{K_1c} we see that there exists finite
dimensional $C^*$-algebras $F_1$ and $F_2$ such that we have a short exact
sequence of the form
$$
  \begin{CD}
      0 @>  >> SF_1  @>  >> A @>  >> F_2 @>  >> 0. \\   
  \end{CD}
$$
Apply \cite[Proposition 2.5]{TNS} to this short exact sequence and
the one obtained by tensoring with $\mathcal K$ to obtain two long exact
sequences for $k_n$. It is well-known that the canonical maps
$k_i(F_2)\to k_i(F_2\otimes\mathcal K)$ and
$k_i(SF_1)\to k_i(SF_1\otimes\mathcal K)$ are isomorphisms for $i=-1,0$
(cf. \cite[Lemma 2.3]{TNS}), so the theorem follows from the five lemma in
algebra.
\end{proof}
\end{lemma}

Let $A$ be a unital $C^*$-algebra. Let
$\pi_A:U(A)/\overline{DU(A)}\to K_1(A)$ denote the group homomorphism
$\pi_A(q_A'(u))=[u]$.

\begin{prop}\label{exact}
Let $A$ be a unital inductive limit of a sequence of finite direct sums
of building blocks. There exists a group homomorphism
\begin{gather*}
\lambda_A:\text{Aff}\, T(A)/\overline{\rho_A(K_0(A))}\to
U(A)/\overline{DU(A)}, \\
\lambda_A(q_A(\widehat{a}))=q_A'(e^{2\pi ia}),\quad a\in A_{\text{sa}}.
\end{gather*}
This map is an isometry when $\text{Aff}\, T(A)/\overline{\rho_A(K_0(A))}$
is equipped with the metric $d_A$, and it gives rise to a short exact sequence
of abelian groups
$$
  \begin{CD}
      0 @>  >> \text{Aff}\, T(A)/\overline{\rho_A(K_0(A))}  @>\lambda_A  >>
U(A)/\overline{DU(A)} @>\pi_A  >> K_1(A) @>  >> 0. \\   
  \end{CD}
$$
This sequence is natural in $A$ and splits unnaturally.

\begin{proof}
Combine Lemma \ref{cano} with \cite[Lemma 6.4]{TATD}.
\end{proof}
\end{prop}

\begin{prop}\label{detunitgr}
Let $A=A(n,d_1,d_2,\dots,d_N)$ be a building block. Let $u\in A$ be a
unitary. Assume that
\begin{gather*}
Det(u(z))=1,\quad z\in\T, \\
Det(\Lambda_i(u))=1,\quad i=1,2,\dots,N.
\end{gather*}
Then $u\in\overline{DU(A)}$.

\begin{proof}
First note that $[u]=0$ in $K_1(A)$ by Lemma \ref{K1det}. Hence
$q_A'(u)=q_A'(e^{2\pi ia})$ by Proposition \ref{exact} for some self-adjoint
element $a\in A$. Since
$$
Det(u(z))=Det(e^{2\pi ia(z)})=e^{2\pi iTr(a(z))}
$$
it follows that $Tr(a(z))=k$ for some $k\in\Z$ and all $z\in\T$. Hence
$\widehat{a}=\frac kn \widehat{1}$ in $\text{Aff}\, T(A)$ by Lemma \ref{trace}.
By applying $\lambda_A$ we get that
$q_A'(u)=q_A'(e^{2\pi ia})=q_A'(\lambda 1)$, where
$\lambda=e^{2\pi i\frac kn}$. Since $Det(\Lambda_i(u))=1$ we see that
$\lambda^{d_i}=1$, $i=1,2,\dots,N$. Thus $\lambda^d=1$ where
$d=gcd(d_1,d_2,\dots,d_N)$. But then $\frac kn=\frac ld$ for some $l\in\Z$.
It follows by Corollary \ref{imK_0} that
$\frac kn \widehat{1}\in\rho_A(K_0(A))$ and
hence by Proposition \ref{exact} we get that $\lambda 1\in \overline{DU(A)}$.
\end{proof}
\end{prop}

\begin{lemma}\label{detrange}
Let $A=A(n,d_1,d_2,\dots,d_N)$ be a building block with exceptional points
$x_1,x_2,\dots,x_N$. Let $g:\T\to\T$ be a continuous function and let
$h_i\in\T$ be such that $h_i^{\frac{n}{d_i}}=g(x_i)$, $i=1,2,\dots,N$.
There exists a unitary $u\in A$ such that
\begin{gather*}
Det(u(z))=g(z),\quad z\in\T, \\
Det(\Lambda_i(u))=h_i,\quad i=1,2,\dots,N.
\end{gather*}

\begin{proof}
Choose a continuous function $f:\T\to\T$ such that $f(x_i)^{d_i}=h_i$.
Define a unitary $v\in A$ by $v=f\otimes 1$.
Since
$$
f(x_i)^n=h_i^{\frac n{d_i}}=g(x_i),
$$
we can define a unitary $w\in A$ by
$$
w(z)=\text{diag}(g(z)f(z)^{-n},1,1,\dots,1),\quad z\in\T.
$$
Set $u=wv$.
\end{proof}
\end{lemma}

Let $A=A(n,d_1,d_2,\dots,d_N)$ be a building block. By Lemma \ref{detrange}
there exist unitaries $w_1^A,w_2^A,\dots,w_N^A\in A$ such that
$Det(w_k^A(z))=1$, $z\in\T$, $k=1,2,\dots,N$, and such that
$$
Det(\Lambda_l(w_k^A))=
\begin{cases}
1 & l\neq k, \\
\exp(2\pi i\frac {d_l}n) & l=k.
\end{cases}
$$

Let $A=A(n,d_1,d_2,\dots,d_N)$ and $B=A(m,e_1,e_2,\dots,e_M)$ be building
blocks. Let $\varphi : A\to B$ be a unital *-homomorphism. As in
\cite[Chapter 1]{TATD} we define $s^{\varphi}(j,i)$ to be the multiplicity of
the representation $\Lambda_i^A$ in the representation
$\Lambda_j^B\circ\varphi$ for $i=1,2,\dots,N$, $j=1,2,\dots,M$.

The following theorem shows that there is a connection between $KK(A,B)$ and
the torsion subgroups of $U(A)/\overline{DU(A)}$ and $U(B)/\overline{DU(B)}$
when $A$ and $B$ are building blocks.

\begin{thm}\label{sr}
Let $A=A(n,d_1,d_2,\dots,d_N)$ and $B=A(m,e_1,e_2,\dots,e_M)$ be building
blocks and let $\varphi,\psi:A\to B$ be unital *-homomorphisms. The
following are equivalent.
\begin{itemize}
\item [(i)]
$\varphi^*=\psi^*$ in $Hom(K^0(B),K^0(A))$,
\item [(ii)]
$s^{\varphi}(j,i)\equiv s^{\psi}(j,i)\quad\text{mod}\ \frac n{d_i},\quad
i=1,2,\dots,N,\ j=1,2,\dots,M$,
\item [(iii)]
$\varphi^{\#}(x)=\psi^{\#}(x),\quad x\in Tor\bigl(U(A)/\overline{DU(A)}\bigr)$,
\item [(iv)]
$\varphi^{\#}(q_A'(w_k^A))=\psi^{\#}(q_A'(w_k^A)),\quad k=1,2,\dots,N$.
\end{itemize}

\begin{proof}
For each $i=1,2,\dots,N$, $j=1,2,\dots,M$, let $r_i^j$ and $s_i^j$ be the
integers with $0\le r_i^j,s_i^j<\frac n{d_i}$, and
$r_i^j\equiv s^{\varphi}(j,i)\ \text{mod}\ \frac n{d_i}$,
$s_i^j\equiv s^{\psi}(j,i)\ \text{mod}\ \frac n{d_i}$.
By Lemma \ref{repr} there exist
$a_1^j,\dots,a_{K_j}^j,b_1^j,\dots,b_{L_j}^j\in\T$ and unitaries
$u_j,v_j\in M_{e_j}$ such that
\begin{align}\label{equiveq1}
\Lambda_j^B\circ\varphi(f)&=u_j\text{diag}(\Lambda_1^{r_1^j}(f),
\Lambda_2^{r_2^j}(f),\dots,
\Lambda_N^{r_N^j}(f),f(a_1^j),f(a_2^j),\dots,f(a_{K_j}^j))u_j^*, \\
\label{equiveq2}
\Lambda_j^B\circ\psi(f)&=v_j\text{diag}(\Lambda_1^{s_1^j}(f),
\Lambda_2^{s_2^j}(f),\dots,
\Lambda_N^{s_N^j}(f),f(b_1^j),f(b_2^j),\dots,f(b_{L_j}^j))v_j^*.
\end{align}
Since
$$
e_j=\sum_{i=1}^N r_i^jd_i+K_jn=\sum_{i=1}^N s_i^jd_i+L_jn,
$$
we remark that if (ii) holds then $K_j=L_j$, $j=1,2,\dots,M$.

Note that
\begin{align*}
\varphi^*([\Lambda_j^B])&=
\sum_{i=1}^N r_i^j[\Lambda_i^A]+K_j\frac n{d_N}[\Lambda_N^A]=
\sum_{i=1}^{N-1} r_i^j[\Lambda_i^A]+(K_j\frac n{d_N}+r_N^j)[\Lambda_N^A], \\
\psi^*([\Lambda_j^B])&=
\sum_{i=1}^N s_i^j[\Lambda_i^A]+L_j\frac n{d_N}[\Lambda_N^A]=
\sum_{i=1}^{N-1} s_i^j[\Lambda_i^A]+(L_j\frac n{d_N}+s_N^j)[\Lambda_N^A].
\end{align*}
By Proposition \ref{sf} we see that $\varphi^*=\psi^*$ if and only if for every
$j=1,2,\dots,M$,
$$
K_j\frac n{d_N}+r_N^j=L_j\frac n{d_N}+s_N^j, \qquad\text{and}\qquad
r_i^j=s_i^j,\quad i=1,2,\dots,N-1.
$$
It follows that (i) holds if and only if $r_i^j=s_i^j$ and $K_j=L_j$
for every $i,j$. But this statement is equivalent to (ii) by the remark above.

Assume (ii) holds. To prove (iii), let $u\in A$ be a unitary such that
$q_A(u)$ has finite order in the group $U(A)/\overline{DU(A)}$. Then
$Det(u(\cdot))$ is constant. By (\ref{equiveq1}), (\ref{equiveq2}), and since
$K_j=L_j$, $j=1,2,\dots,M$, it follows that
$$
Det(\Lambda_j(\varphi(u)))=Det(\Lambda_j(\psi(u))),\quad j=1,2,\dots,M.
$$
In particular, $Det(\varphi(u)(\cdot))$ equals $Det(\psi(u)(\cdot))$ at the
exceptional points of $B$. On the other hand, $Det(\varphi(u)(\cdot))$
and $Det(\psi(u)(\cdot))$ are constant functions on $\T$ and are hence equal
everywhere. We may therefore use Proposition \ref{detunitgr} to conclude that
$\varphi^{\#}(q_A(u))=\psi^{\#}(q_A(u))$.
(iii) $\Rightarrow$ (iv) is trivial. Assume (iv). By (\ref{equiveq1}) and
(\ref{equiveq2}),
$$
\exp(2\pi i\frac{d_k}n r_k^j)=
Det(\Lambda_j^B\circ\varphi(w_k^A))=
Det(\Lambda_j^B\circ\psi(w_k^A))=
\exp(2\pi i\frac{d_k}n s_k^j).
$$
Hence $r_k^j=s_k^j$ for $k=1,2,\dots,N$, $j=1,2,\dots,M$, and
we have (ii).
\end{proof}
\end{thm}

\begin{prop}\label{KKunid}
Let $A$ and $B$ be finite direct sums of building blocks, let
$\varphi,\psi:A\to B$ be unital *-homomorphisms, and let $x$ be an element of
finite order in the group $U(A)/\overline{DU(A)}$. If $[\varphi]=[\psi]$ in
$KK(A,B)$ then $\varphi^{\#}(x)=\psi^{\#}(x)$.

\begin{proof}
We may assume that $B$ is a building block rather than a finite direct sum
of building blocks. Let $A=A_1\oplus A_2\oplus\dots\oplus A_R$ where each $A_i$
is a building block, and let $\iota_i:A_i\to A$ denote the inclusion.
Let $p_1,p_2,\dots,p_R$ be the minimal non-zero central
projections in $A$. Since $\varphi_*[p_i]=\psi_*[p_i]$ in $K_0(B)$, it follows
from Lemma \ref{cd} that there is a unitary $u\in B$ such that
$u\varphi(p_i)u^*=\psi(p_i)$, $i=1,2,\dots,R$. Hence we may assume that
$\varphi(p_i)=\psi(p_i)$, $i=1,2,\dots,R$. Set $q_i=\varphi(p_i)$.

Let $\varphi_i,\psi_i: A_i\to q_iBq_i$ be the induced maps and let
$\epsilon_i:q_iBq_i\to B$ be the inclusion, $i=1,2,\dots,R$. If $q_i\neq 0$
then $[\epsilon_i]\in KK(q_iBq_i,B)$ is a $KK$-equivalence by
\cite[Theorem 7.3]{RSUCT}. Thus
$$
[\varphi_i]=[\epsilon_i]^{-1}\cdot [\varphi]\cdot [\iota_i]=
[\epsilon_i]^{-1}\cdot [\psi]\cdot [\iota_i]=[\psi_i]
$$
in $KK(A_i,q_iBq_i)$.
Let $x=q_A'(u)$ where $u\in A$ is a unitary. Let $u=\sum_{i=1}^R \iota_i(u_i)$
where $u_i\in A_i$. By Theorem \ref{sr} and
Corollary \ref{cutdown} we see that
$\varphi_i(u_i)=\psi_i(u_i)$ mod $\overline{DU(q_iBq_i)}$
and thus $\epsilon_i\circ\varphi_i(u_i)+(1-q_i)=
\epsilon_i\circ\psi_i(u_i)+(1-q_i)$ mod $\overline{DU(B)}$. Hence
\begin{align*}
\varphi(u)&=\prod_{i=1}^R \varphi(\iota_i(u_i)+(1-p_i))=
\prod_{i=1}^R (\epsilon_i\circ\varphi_i(u_i)+(1-q_i)) \\
&=\prod_{i=1}^R (\epsilon_i\circ\psi_i(u_i)+(1-q_i))=
\prod_{i=1}^R \psi(\iota_i(u_i)+(1-p_i))=\psi(u)
\end{align*}
modulo $\overline{DU(B)}$.
\end{proof}
\end{prop}

\section{Homomorphisms between building blocks}

In this section we improve a result of Thomsen on
*-homomorphisms between building blocks that will be needed in the next
section.

Whenever $\theta_1,\theta_2,\dots,\theta_L$ are real numbers such that
$$
\theta_1\le\theta_2\le\dots\le\theta_L\le\theta_1+1,
$$
it will be convenient for us in the following to define $\theta_n$ for every
$n\in\Z$ by the formula $\theta_{pL+r}=\theta_r+p$, where
$p\in\Z$, $r=1,2,\dots,L$. Note that for every $n\in\Z$,
$$
\theta_n\le\theta_{n+1}\le\dots\le\theta_{n+L}\le\theta_n+1,
$$
and
$$
(e^{2\pi i\theta_1},e^{2\pi i\theta_1},\dots,e^{2\pi i\theta_L})=
(e^{2\pi i\theta_n},e^{2\pi i\theta_{n+1}},\dots,e^{2\pi i\theta_{n+L}})
$$
as unordered $L$-tuples.

\begin{lemma}\label{ordert}
Let $a_1,a_2,\dots,a_L\in\T$ and let $k$ be an integer. There exist real
numbers $\theta_1,\theta_2,\dots,\theta_L$ such that
$$
\theta_1\le\theta_2\le\dots\le\theta_L\le\theta_1+1,
$$
such that $\sum_{r=1}^L \theta_r\in [k,k+1[$ and such that
$$
(a_1,a_2,\dots,a_L)=(e^{2\pi i\theta_1},e^{2\pi i\theta_2},\dots,
e^{2\pi i\theta_L})
$$
as unordered $L$-tuples.

\begin{proof}
Choose $\omega_1,\omega_2,\dots,\omega_L\in [0,1[$ such that
$\omega_1\le\omega_2\le\dots\le \omega_L$ and such that
$$
(a_1,a_2,\dots,a_L)=(e^{2\pi i\omega_1},e^{2\pi i\omega_2},\dots,
e^{2\pi i\omega_L})
$$
as unordered $L$-tuples. Let $l$ be the integer such that
$\sum_{r=1}^L \omega_r\in [l,l+1[$. Set $\theta_r=\omega_{r+k-l}$.
\end{proof}
\end{lemma}

\begin{lemma}\label{glue}
Assume that
$$
\bigl(\exp(2\pi i \theta_1),\dots,\exp(2\pi i \theta_L)\bigr)=
\bigl(\exp(2\pi i \omega_1),\dots,\exp(2\pi i \omega_L)\bigr)
$$
as unordered $L$-tuples, where
$\theta_1,\theta_2,\dots,\theta_L$ and $\omega_1,\omega_2,\dots,\omega_L$
are real numbers such that
\begin{gather*}
\theta_1\le\theta_2\le\dots\le\theta_L\le \theta_1+1, \\
\omega_1\le\omega_2\le\dots\le\omega_L\le \omega_1+1.
\end{gather*}
Then $\theta_j=\omega_{r+j}$, $j=1,2,\dots,L$, where
$r=\sum_{j=1}^L (\theta_j-\omega_j)$.

\begin{proof}
Choose $m\in\Z$ such that $\theta_m<\theta_{m+1}$ and choose $n\in\Z$
such that
$$
\theta_{m+1}=\omega_{n+1}>\omega_{n}.
$$
Assume that
$$
\theta_{m+p}=\omega_{n+q}+k.
$$
for some integers $p,q$ with $1\le p\le L$, $1\le q\le L$, and an integer $k$.
Then
\begin{gather*}
-1<\theta_{m+1}-\theta_{m+L}\le \theta_{m+1}-\theta_{m+p}=
\omega_{n+1}-\omega_{n+q}-k\le -k, \\
0\le \theta_{m+p}-\theta_{m+1}=\omega_{n+q}+k-\omega_{n+1}<
\omega_{n+q}+k-\omega_n\le 1+k.
\end{gather*}
Hence $k=0$. By assumption it follows that for every $x\in\R$,
$$
\#\{j=1,2,\dots,L:\theta_{m+j}=x\}=\#\{j=1,2,\dots,L:\omega_{n+j}=x\}.
$$
Thus
$$
\bigl( \theta_{m+1},\theta_{m+2},\dots,\theta_{m+L}\bigr) =
\bigl( \omega_{n+1},\omega_{n+2},\dots,\omega_{n+L}\bigr)
$$
as unordered $L$-tuples. Therefore
$$
\theta_{m+j}=\omega_{n+j},\quad j=1,2,\dots,L.
$$
Hence $\theta_j=\omega_{n-m+j}$ for $j=1,2,\dots,L$. From this it follows
that $r=n-m$.
\end{proof}
\end{lemma}

\begin{prop}\label{kont}
Let $\lambda_1,\lambda_2,\dots,\lambda_L:[0,1]\to\T$ be continuous functions
and let $k$ be an integer. There exist continuous functions
$F_1,F_2,\dots,F_L:[0,1]\to\R$ such that
$\sum_{j=1}^L F_j(0)\in [k,k+1[$ and such that for each
$t\in [0,1]$,
$$
F_1(t)\le F_2(t)\le\dots\le F_L(t)\le F_1(t)+1,
$$
and
$$
\bigl(\lambda_1(t),\lambda_2(t),\dots,\lambda_L(t)\bigr)=
\bigl(\exp(2\pi iF_1(t)),\exp(2\pi iF_2(t)),\dots,\exp(2\pi iF_L(t))\bigr)
$$
as unordered $L$-tuples.

\begin{proof}
Choose a positive integer $n$ such that
$$
|s-t|\le \frac 1n\ \Rightarrow\ \rho(\lambda_j(s),\lambda_j(t))<\frac 1{2L},
\quad s,t\in [0,1],\ j=1,2,\dots,L.
$$
We will prove by induction in $m$ that there exist continuous
functions $F_1,\dots,F_L$ that satisfy the above for $t\in [0,\frac mn]$.
The case $m=0$ follows from Lemma \ref{ordert}.

Now assume that for some $m=0,1,\dots,n-1$ we have constructed continuous
functions $F_1,F_2,\dots,F_L:[0,\frac mn]\to \R$ such that
$\sum_{j=1}^L F_j(0)\in [k,k+1[$, and such that for each $t\in [0,\frac mn]$,
$F_1(t)\le F_2(t)\le\dots\le F_L(t)\le F_1(t)+1$, and
$$
\bigl(\lambda_1(t),\lambda_2(t),\dots,\lambda_L(t)\bigr)=
\bigl(\exp(2\pi iF_1(t)),\exp(2\pi iF_2(t)),\dots,\exp(2\pi iF_L(t))\bigr)
$$
as unordered $L$-tuples.
Choose $\alpha_m\in\R$ such that
$\rho(e^{2\pi i\alpha_m},\lambda_j(\frac mn))\ge \frac 1{2L}$ for
$j=1,2,\dots,L$. Choose continuous functions
$G_j:[\frac mn,\frac{m+1}n]\to ]\alpha_m,\alpha_m+1[$
such that for each $t\in [\frac mn,\frac {m+1}n]$,
$$
G_1(t)\le G_2(t)\le\dots\le G_L(t)
$$
and
$$
\bigl(\lambda_1(t),\lambda_2(t),\dots,\lambda_L(t)\bigr)=
\bigl(\exp(2\pi iG_1(t)),\exp(2\pi iG_2(t)),\dots,\exp(2\pi iG_L(t))\bigr)
$$
as unordered $L$-tuples. Set for $j=1,2,\dots,L$, $p\in\Z$,
$$
G_{pL+j}(t)=G_j(t)+p,\quad t\in [\frac mn,\frac{m+1}n].
$$
By Lemma \ref{glue} there exists an integer $r$ such that for $j=1,2,\dots,L$,
$$
F_j(\frac mn)=G_{r+j}(\frac mn).
$$
Define for $j=1,2,\dots,L$, a continuous function
$F_j':[0,\frac{m+1}n]\to\R$ by
$$
F_j'(t)=
\begin{cases}
F_j(t)
& t\in [0,\frac mn], \\
G_{r+j}(t)
& t\in [\frac mn,\frac{m+1}n].
\end{cases}
$$
$F_1',F_2',\dots,F_L'$
satisfy the conclusion of the lemma for $t\in [0,\frac{m+1}n]$.
\end{proof}
\end{prop}

\begin{prop}\label{stand}
Let $A=A(n,d_1,d_2,\dots,d_N)$ and $B=A(m,e_1,e_2,\dots,e_M)$ be building
blocks and let $\varphi:A\to B$ be a unital *-homomorphism.
There exist integers $r_1,r_2,\dots,r_N$ with $0\le r_i<\frac n{d_i}$, an
integer $L\ge 0$, and a unitary $w\in M_m$ such that if $\psi:A\to B$ is
a unital *-homomorphism with
$\varphi^{\#}(q_A'(\omega_k^A))=\psi^{\#}(q_A'(\omega_k^A))$, $k=1,2,\dots,N$,
and if $\gamma:\T\to\R$ is a continuous function such that
$$
Det(\psi(v^A)(z))=Det(\varphi(v^A)(z))\exp(2\pi i\gamma(z)),\quad z\in\T,
$$
then $\varphi$ and $\psi$ are approximately unitarily
equivalent to *-homomorphisms of the form
\begin{gather*}
\varphi'(f)(e^{2\pi it})=u(t)\text{diag}\bigl(
\Lambda_1^{r_1}(f),\dots,\Lambda_N^{r_N}(f),
f(e^{2\pi iF_1(t)}),\dots,f(e^{2\pi iF_L(t)})\bigr)u(t)^*, \\
\psi'(f)(e^{2\pi it})=v(t)\text{diag}\bigl(
\Lambda_1^{r_1}(f),\dots,\Lambda_N^{r_N}(f),
f(e^{2\pi iG_1(t)}),\dots,f(e^{2\pi iG_L(t)})\bigr)v(t)^*,
\end{gather*}
where $u,v\in C[0,1]\otimes M_m$ are unitaries with $u(0)=v(0)=1$,
$u(1)=v(1)=w$, and $F_1,F_2,\dots,F_L:[0,1]\to\R$ and
$G_1,G_2,\dots,G_L:[0,1]\to\R$
are continuous functions such that for every $t\in [0,1]$,
\begin{gather*}
F_1(t)\le F_2(t)\le\dots\le F_L(t)\le F_1(t)+1, \\
G_1(t)\le G_2(t)\le\dots\le G_L(t)\le G_1(t)+1,
\end{gather*}
and such that $\gamma(e^{2\pi it})=\sum_{r=1}^L (G_r(t)-F_r(t))$
for every $t\in [0,1]$.

\begin{proof}
By \cite[Chapter 1]{TATD} it follows that $\varphi$ is approximately
unitarily equivalent to a *-homomorphism $\varphi_1:A\to B$ of the form
\begin{equation}\label{standphi}
\varphi_1(f)(e^{2\pi it})=u_0(t)\text{diag}(\Lambda_1^{r_1}(f),\dots,
\Lambda_N^{r_N}(f),f(\lambda_1(t)),\dots,f(\lambda_L(t)))u_0(t)^*
\end{equation}
for $t\in [0,1]$, where $r_1,r_2,\dots,r_N$ are integers,
$0\le r_i<\frac n{d_i}$, $i=1,2,\dots,N$,
$\lambda_1,\lambda_2,\dots,\lambda_L:[0,1]\to\T$ are continuous functions,
and $u_0\in C[0,1]\otimes M_m$ is a unitary.
Let $l$ denote the winding number of $Det(\varphi(v^A)(\cdot))$.
Let $y$ be a unitary $(Ln)\times (Ln)$ matrix such that
$$
y\, \text{diag}(a_1,a_2,\dots,a_L)\, y^*=\text{diag}(a_L,a_1,a_2,\dots,a_{L-1})
$$
for all $a_1,a_2,\dots,a_L\in M_m$. Set
$$
w=\text{diag}(\underbrace{1,1,\dots,1}_{m-Ln\text{\ times}},y^l).
$$

Let now $\psi:A\to B$ be given. As above $\psi$ is approximately unitarily
equivalent to a *-homomorphism $\psi_1:A\to B$ of the form
$$
\psi_1(f)(e^{2\pi it})=v_0(t)\text{diag}\bigl(\Lambda_1^{s_1}(f),
\dots,\Lambda_N^{s_N}(f),f(\mu_1(t)),\dots,f(\mu_K(t))\bigr)v_0(t)^*.
$$
Note that
\begin{gather*}
s^{\varphi}(j,i)\frac m{e_j}=s^{\varphi_1}(j,i)\frac m{e_j}=
r_i+\#\{r=1,2,\dots,L:\lambda_r(y_j)=x_i\}\frac n{d_i}, \\
s^{\psi}(j,i)\frac m{e_j}=s^{\psi_1}(j,i)\frac m{e_j}=
s_i+\#\{r=1,2,\dots,K:\mu_r(y_j)=x_i\}\frac n{d_i}.
\end{gather*}
By Theorem \ref{sr} it follows that $r_i=s_i$, $i=1,2,\dots,N$. And since
$$
m=Kn+\sum_{i=1}^N s_id_i=Ln+\sum_{i=1}^N r_id_i
$$
we see that $K=L$.

By Proposition \ref{kont} choose continuous functions
$F_1,F_2,\dots,F_L:[0,1]\to\R$ such that for every $t\in [0,1]$,
$$
F_1(t)\le F_2(t)\le\dots\le F_L(t)\le F_1(t)+1,
$$
and such that
$$
(\lambda_1(t),\lambda_2(t),\dots,\lambda_L(t))=
(e^{2\pi iF_1(t)},e^{2\pi iF_2(t)},\dots,e^{2\pi iF_L(t)})
$$
as unordered $L$-tuples for each $t\in [0,1]$. Again by Proposition \ref{kont}
there exist continuous functions $G_1,G_2,\dots,G_L:[0,1]\to\R$ such that
for every $t\in [0,1]$,
$$
G_1(t)\le G_2(t)\le\dots\le G_L(t)\le G_1(t)+1,
$$
such that
$$
(\mu_1(t),\mu_2(t),\dots,\mu_L(t))=
(e^{2\pi iG_1(t)},e^{2\pi iG_2(t)},\dots,e^{2\pi iG_L(t)})
$$
as unordered $L$-tuples for each $t\in [0,1]$, and such that
\begin{equation}\label{standdet}
|\sum_{r=1}^L G_r(0)-\sum_{r=1}^L F_r(0)+\gamma(1)|<1.
\end{equation}
It follows from (\ref{standphi}) that
$$
\bigl(\exp(2\pi iF_1(0)),\dots,\exp(2\pi iF_L(0))\bigr)=
\bigl(\exp(2\pi iF_1(1)),\dots,\exp(2\pi iF_L(1))\bigr)
$$
as unordered $L$-tuples.
Since $l=\sum_{r=1}^L (F_r(1)-F_r(0))$ we see by Lemma \ref{glue} that
$F_r(1)=F_{r+l}(0)$ for each $r=1,2,\dots,L$. Similarly, as
$Det(\varphi(v^A)(\cdot))$ and $Det(\psi(v^A)(\cdot))$ have the same winding
number, $G_r(1)=G_{r+l}(0)$, $r=1,2,\dots,L$.

Let $t_1,t_2,\dots,t_M\in ]0,1[$ be numbers such that $e^{2\pi it_j}$,
$j=1,2,\dots,M$, are the exceptional points of $B$.
By (\ref{standphi}) there exist a unitary $u_j\in M_m$ such that
$$
u_j\text{diag}\bigl(\Lambda_1^{r_1}(f),\dots,\Lambda_N^{r_N}(f),
f(e^{2\pi iF_1(t_j)}),\dots,f(e^{2\pi iF_L(t_j)})\bigr)u_j^*\in M_{e_j}
\subseteq M_m
$$
for every $f\in A$. Choose a unitary $u\in C[0,1]\otimes M_m$ such that
$u(t_j)=u_j$, $j=1,2,\dots,M$, $u(0)=1$ and $u(1)=w$.
Note that for every $f\in A$,
\begin{align*}
&u(0)\text{diag}\bigr(\Lambda_1^{r_1}(f),\dots,\Lambda_N^{r_N}(f),
f(e^{2\pi iF_1(0)}),\dots,f(e^{2\pi iF_L(0)})\bigl)u(0)^* \\
=\, &u(1)\text{diag}\bigr(\Lambda_1^{r_1}(f),\dots,\Lambda_N^{r_N}(f),
f(e^{2\pi iF_1(1)}),\dots,f(e^{2\pi iF_L(1)})\bigl)u(1)^*.
\end{align*}
It follows that we may define a unital *-homomorphism $\varphi':A\to B$ by
$$
\varphi'(f)(e^{2\pi it})=u(t)\text{diag}\bigl(
\Lambda_1^{r_1}(f),\dots,\Lambda_N^{r_N}(f),
f(e^{2\pi iF_1(t)}),\dots,f(e^{2\pi iF_L(t)})\bigr)u(t)^*, \\
$$
for $f\in A$, $t\in [0,1]$. Then for every $f\in A$, $z\in\T$,
$$
Tr(\varphi(f)(z))=Tr(\varphi_1(f)(z))=Tr(\varphi'(f)(z)).
$$
Hence $\varphi$ and $\varphi'$ are approximately unitarily equivalent
by \cite[Theorem 1.4]{TATD}.

Similarly we see that there exists a unitary $v\in C[0,1]\otimes M_m$ such that
$v(0)=1$ and $v(1)=w$, and such that
$$
\psi'(f)(e^{2\pi it})=v(t)\text{diag}\bigl(
\Lambda_1^{r_1}(f),\dots,\Lambda_N^{r_N}(f),
f(e^{2\pi iG_1(t)}),\dots,f(e^{2\pi iG_L(t)})\bigr)v(t)^*
$$
defines a *-homomorphism that is approximately unitarily equivalent to $\psi$.
Finally note that by (\ref{standdet}) we have that
$\gamma(e^{2\pi it})=\sum_{r=1}^L (G_r(t)-F_r(t))$ for every $t\in [0,1]$.
\end{proof}
\end{prop}

\section{Uniqueness}
The purpose of this section is to prove a uniqueness theorem, i.e a theorem
saying that two unital *-homomorphisms between (finite direct sums of)
building blocks are close in a suitable sense if they approximately agree
on the invariant. Many of the arguments here are inspired by similar arguments
in \cite{EAIS}, \cite{EATS}, \cite{NT}, \cite{TATD}, and \cite{JSSP}.

We start out with some definitions. Let $k$ be a positive integer.
A $k$-arc is an arc-segment of the form
$$
I=\{ e^{2\pi it} : t\in [\frac mk,\frac nk]\}
$$
where $m$ and $n$ are integers, $m<n$. We set
$$
I\pm \epsilon=\{ e^{2\pi it} : t\in [\frac mk-\epsilon,\frac nk+\epsilon]\}.
$$
Define a metric on the set of unordered $L$-tuples
consisting of elements from $\T$ by
$$
R_L\bigl((a_1,a_2,\dots,a_L),(b_1,b_2,\dots,b_L)\bigr)=
\min_{\sigma\in \Sigma_L}\bigl(\max_{1\le i\le L}\rho(a_i,b_{\sigma(i)})\bigr),
$$
where $\Sigma_L$ denotes the group of permutations of the set
$\{1,2,\dots,L\}$. It follows from Lemma \ref{cirperm} below that it suffices
to take the minimum over a certain subset of $\Sigma_L$.

\begin{lemma}\label{mar}
Let $a_1,a_2,\dots,a_L,b_1,b_2,\dots,b_L\in\T$ and let $\epsilon>0$.
Assume that there is a positive integer $k$ such that
$$
\#\{r:a_r\in I\}\le \#\{r:b_r\in I\pm\epsilon\}
$$
for all $k$-arcs $I$.
Then
$$
R_L\bigl((a_1,a_2,\dots,a_L),(b_1,b_2,\dots,b_L)\bigr) \le \epsilon+\frac 1k.
$$

\begin{proof}
For $j=1,2,\dots,L$, set
$$
X_j=\{x\in\T:\rho(x,a_j)\le \epsilon+\frac 1k\}
$$
and
$$
C_j=\{r=1,2,\dots,L:b_r\in X_j\}.
$$
Let $S\subseteq \{1,2,\dots,L\}$ be arbitrary. We will show that
$\#S\le \#\cup_{j\in S} C_j$.

Let $Y_1,Y_2,\dots,Y_m$ be the connected components of
$\cup_{j\in S} X_j$. Choose for each $n=1,2,\dots,m$ a $k$-arc $I_n$ such that
$I_n\pm\epsilon\subseteq Y_n$ and $\{a_j:j\in S\}\cap Y_n\subseteq I_n$. Then
$$
\#\{r:a_r\in I_n\}\le\#\{r:b_r\in I_n\pm\epsilon\}\le \#\{r:b_r\in Y_n\}.
$$
If $r\in S$ then $a_r\in I_n$ for some $n$. Hence
$$
\#S\le \#\{r:a_r\in\bigcup_{n=1}^m I_n\}\le\#\{r:b_r\in \bigcup_{n=1}^m Y_n\}=
\#\{r:b_r\in\bigcup_{j\in S} X_j\}=\#\bigcup_{j\in S} C_j.
$$

By Hall's marriage lemma, see e.g
\cite[Theorem 2.2]{BOLLO}, the sets $C_j$, $j=1,2,\dots,L$, have distinct
representatives. In other words, there exists a permutation $\sigma$
of $\{1,2,\dots,L\}$ such that $\rho(a_j,b_{\sigma(j)})\le\epsilon+\frac 1k$
for $j=1,2,\dots,L$.
\end{proof}
\end{lemma}

\begin{lemma}\label{intperm}
Let $a_1\le a_2\le\dots\le a_L$ and $b_1\le b_2\le\dots\le b_L$ be real
numbers and let $\sigma$ be a permutation of $\{1,2,\dots,L\}$. Then
$$
\max_{1\le j\le L} |a_j-b_j| \le \max_{1\le j\le L} |a_j-b_{\sigma(j)}|
$$

\begin{proof}
Let $\epsilon=\max |a_j-b_{\sigma(j)}|$. If e.g
$b_j<a_j-\epsilon$ for some $j$ then $\sigma$ must map the set
$\{1,2,\dots,j\}$ into $\{1,2,\dots,j-1\}$. Contradiction.
\end{proof}
\end{lemma}

The corresponding statement for the circle is slightly more complicated. It
can be viewed as a generalization of Lemma \ref{glue}.

\begin{lemma}\label{cirperm}
Let $\theta_1,\theta_2,\dots,\theta_L$ and $\omega_1,\omega_2,\dots,\omega_L$
be real numbers such that
\begin{gather*}
\theta_1\le\theta_2\le\dots\le\theta_L\le \theta_1+1, \\
\omega_1\le\omega_2\le\dots\le\omega_L\le \omega_1+1.
\end{gather*}
There exists an integer $p$ such that
$$
\max_{1\le j\le L} |\theta_j-\omega_{j+p}| =
R_L\bigl((e^{2\pi i\theta_1},\dots,e^{2\pi i\theta_L}),
(e^{2\pi i\omega_1},\dots,e^{2\pi i\omega_L})\bigr).
$$

\begin{proof}
Let $\epsilon=R_L\bigl((e^{2\pi i\theta_1},\dots,e^{2\pi i\theta_L}),
(e^{2\pi i\omega_1},\dots,e^{2\pi i\omega_L})\bigr)$. Note that
$0\le\epsilon\le\frac 12$. There exist $y_1,y_2,\dots,y_L\in\R$ such that
\begin{equation}\label{cirpermeq}
(e^{2\pi i\omega_1},e^{2\pi i\omega_2},\dots,e^{2\pi i\omega_L})=
(e^{2\pi iy_1},e^{2\pi iy_2},\dots,e^{2\pi iy_L})
\end{equation}
as unordered tuples, and such that
\begin{equation}\label{cyceq}
|\theta_j-y_j|\le\epsilon,\quad j=1,2,\dots,L.
\end{equation}
By Lemma \ref{intperm} we may assume that $y_1\le y_2\le\dots\le y_L$ and
still have that (\ref{cirpermeq}) and (\ref{cyceq}) hold.

Choose an integer $n$, $0\le n\le L-1$, such that $y_1,y_2,\dots,y_n<y_L-1$
and $y_L-1\le y_{n+1},y_{n+2},\dots,y_L$. Then
$y_1+1,\dots,y_n+1\in [y_L-1,y_L]$ since $y_L\le y_1+2$ by (\ref{cyceq}).
Choose $z_1,z_2,\dots,z_L\in [y_L-1,y_L]$ such
that $z_1\le z_2\le\dots\le z_L$ and
$$
(z_1,z_2,\dots,z_L)=(y_1+1,\dots,y_n+1,y_{n+1},\dots,y_L)
$$
as unordered $L$-tuples. By (\ref{cyceq}) and Lemma \ref{intperm} we see that
$\max |z_j-\theta_{n+j}|\le\epsilon$. By (\ref{cirpermeq}) and Lemma
\ref{glue} we have that $z_j=\omega_{j+m}$ for
some integer $m$. Hence
$$
\max_j |\theta_j-\omega_{j+m-n}|=\max_j |\theta_{n+j}-\omega_{j+m}|
\le\epsilon.
$$
The reversed inequality is trivial.
\end{proof}
\end{lemma}

The following lemma is fundamental in the proof of Theorem \ref{uns}.

\begin{lemma}\label{det}
Let $\theta_1,\theta_2,\dots,\theta_L$ and $\omega_1,\omega_2,\dots,\omega_L$
be real numbers such that
\begin{gather*}
\theta_1\le\theta_2\le\dots\le\theta_L\le \theta_1+1, \\
\omega_1\le\omega_2\le\dots\le\omega_L\le \omega_1+1,
\end{gather*}
and $|\sum_{j=1}^L (\theta_j-\omega_j)|<\delta$ for some $\delta>0$.
Let $\epsilon>0$ satisfy that $L\epsilon\le\delta$ and
$$
R_L\bigl( (e^{2\pi i\theta_1},e^{2\pi i\theta_2},\dots,e^{2\pi i\theta_L}),
(e^{2\pi i\omega_1},e^{2\pi i\omega_2},\dots,e^{2\pi i\omega_L})\bigr)\le
\epsilon.
$$
Assume finally that for some positive integer $s$,
\begin{equation}\label{sarceq}
\# \{j:e^{2\pi i\omega_j}\in I\} \ge 2\delta,\quad j=1,2,\dots,L,
\end{equation}
for every $s$-arc $I$. Then
$$
|\theta_j-\omega_j|<\epsilon+\frac 2s,
\quad j=1,2,\dots,L.
$$

\begin{proof}
By Lemma \ref{cirperm} there exists an integer $p$ such that
$$
|\theta_j-\omega_{j+p}|\le\epsilon, \quad j=1,2,\dots,L.
$$
Note that
$$
|p|=|\sum_{j=1}^L (\omega_{j+p}-\omega_j)|
\le |\sum_{j=1}^L (\omega_{j+p}-\theta_j)|+
|\sum_{j=1}^L (\theta_j-\omega_j)|<L\epsilon+\delta\le 2\delta.
$$
Fix some $j=1,2,\dots,L$. Set
$$
J=
\begin{cases}
\ \{e^{2\pi it}:\omega_j<t<\omega_{j+p}\}\quad \text{if}\ p\ge 0 \\
\ \{e^{2\pi it}:\omega_{j+p}<t<\omega_{j}\}\quad \text{if}\ p<0.
\end{cases}
$$
Since $\#\{j:e^{2\pi i\omega_j}\in J\}<|p|$ we see by (\ref{sarceq}) that
$J$ cannot contain an $s$-arc. Thus $|\omega_j-\omega_{j+p}|<\frac 2s$.
It follows that
$|\theta_j-\omega_j|<\epsilon+\frac 2s$, $j=1,2,\dots,L$.
\end{proof}
\end{lemma}

Let $A=A(n,d_1,d_2,\dots,d_N)$ be a building block and $p$ a positive integer.
Let $I$ be a $p$-arc. Choose a continuous function $f_A^I:\T\to [0,\frac 1n]$
such that $\emptyset\neq \text{supp}\, f_A^I\subseteq I$ and such that $f_A^I$
equals $0$ at all the exceptional points of $A$. Choose a continuous
function $g_A^I:\T\to [0,1]$ such that $g_A^I$ equals $1$ on $I$,
such that $\text{supp}\, g_A^I\subseteq I\pm \frac 1{2p}$, and such that
$\text{supp}\, g_A^I\setminus I$ contains no exceptional points of $A$. Set
\begin{gather*}
H(A,p)=\{f_A^I\otimes 1: I\ p \text{-arc}\}, \\
\widetilde{H}(A,p)=\{g_A^I\otimes 1: I\ p \text{-arc}\}.
\end{gather*}

\begin{thm}\label{uns}
Let $A=A(n,d_1,d_2,\dots,d_N)$ be a building block. Let $\epsilon>0$ and let
$F\subseteq A$ be a finite set. There exists a positive integer $l_0$
such that if $l$, $p$ and $q$ are positive integers with
$l_0\le l\le p\le q$, if $B=A(m,e_1,e_2,\dots,e_M)$ is a building block, if
$\varphi,\psi:A\to B$ are unital *-homomorphisms, and if $\delta>0$, such that
\begin{itemize}
\item[(i)]
$\widehat{\psi}(\widehat{h})>\frac 8p,\quad h\in H(A,l)$,
\vspace{1 mm}
\item[(ii)]
$\widehat{\psi}(\widehat{h})>\frac 2q,\quad h\in H(A,p)$,
\vspace{1 mm}
\item[(iii)]
$\|\widehat{\varphi}(\widehat{h})-\widehat{\psi}(\widehat{h})\|<\delta,
\quad h\in \widetilde{H}(A,2q)$,
\vspace{1 mm}
\item[(iv)]
$\widehat{\psi}(\widehat{h})>\delta,\quad h\in H(A,4q)$,
\vspace{1 mm}
\item[(v)]
$\varphi^{\#}(q_A'(\omega_k^A))=\psi^{\#}(q_A'(\omega_k^A)),
\quad k=1,2,\dots,N$,
\vspace{1 mm}
\item[(vi)]
$D_B\bigl( \varphi^{\#}(q_A'(v^A)),\psi^{\#}(q_A'(v^A))\bigr)<\frac 1q$;
\end{itemize}
then there exists a unitary $W\in B$ such that
$$
\|\varphi(f)-W\psi(f)W^*\|<\epsilon,\quad f\in F.
$$

\begin{proof}
Choose $l_0$ such that for $x,y\in\T$,
$$
\rho(x,y)\le\frac {6}{l_0}\ \Rightarrow\ \|f(x)-f(y)\|<\frac{\epsilon}6,
\quad f\in F.
$$

Let integers $q\ge p\ge l\ge l_0$, a building block
$B=A(m,e_1,e_2,\dots,e_M)$, and unital *-homomorphisms $\varphi,\psi:A\to B$ be
given such that (i)-(vi) are satisfied. Choose $c>0$ such that for $x,y\in\T$,
\begin{gather*}
\rho(x,y)<c\ \Rightarrow\ \|\varphi(f)(x)-\varphi(f)(y)\|<
\frac{\epsilon}6,\quad f\in F, \\
\rho(x,y)<c\ \Rightarrow\ \|\psi(f)(x)-\psi(f)(y)\|<
\frac{\epsilon}6,\quad f\in F.
\end{gather*}
Let $x_1,x_2,\dots,x_N$ denote the exceptional points of $A$ and let
$y_1,y_2,\dots,y_M$ be those of $B$.
Let for each $j=1,2,\dots,M$, $t_j\in ]0,1[$ be the number such that
$e^{2\pi it_j}=y_j$. Let $\tau:\T\to\T$ be a continuous function such that
$\rho(\tau(z),z)<c$ for every $z\in\T$, and such that for each $j=1,2,\dots,M$,
$\tau$ is constantly equal to $y_j$ on some arc
$$
I_j=\{e^{2\pi it} : t\in [a_j,b_j]\},
$$
where $0<a_j<t_j<b_j<1$. Define a unital *-homomorphism $\chi:B\to B$ by
$\chi(f)=f\circ\tau$. Set $\varphi_1=\chi\circ\varphi$ and
$\psi_1=\chi\circ\psi$. Then
\begin{gather*}
\|\varphi(f)-\varphi_1(f)\|<\frac{\epsilon}6,\quad f\in F, \\
\|\psi(f)-\psi_1(f)\|<\frac{\epsilon}6,\quad f\in F.
\end{gather*}
$\varphi_1$ and $\psi_1$ satisfy (i)-(vi). Let
$$
\varphi_1(v^A)=c\psi_1(v^A)e^{2\pi ib},
$$
where $c\in\overline{DU(B)}$ and $b\in B$ is a self-adjoint element with
$\|b\|<\frac 1q$. Note that
\begin{gather}
Det(\varphi_1(v^A)(z))=Det(\psi_1(v^A)(z))\exp(2\pi iTr(b(z))),\quad z\in\T, \\
\label{dettreq2}
Det(\Lambda_j\circ\varphi_1(v^A))=Det(\Lambda_j\circ\psi_1(v^A))
\exp(2\pi iTr(\Lambda_j(b))),\quad j=1,2,\dots,M.
\end{gather}

Fix some $j=1,2,\dots,M$. Let $\iota_j:M_{e_j}\to M_m$ denote the (unital)
inclusion. By Theorem \ref{sr} and (v) we have that
$s^{\varphi_1}(j,i)\equiv s^{\psi_1}(j,i)\ \text{mod}\ \frac n{d_i}$,
$i=1,2,\dots,N$, $j=1,2,\dots,M$. Choose $s_i^j$, $0\le s_i^j<\frac n{d_i}$,
such that $s_i^j\equiv s^{\varphi_1}(j,i)\ \text{mod}\ \frac n{d_i}$. By Lemma
\ref{repr} we see that for each $z\in I_j$,
\begin{gather*}
\varphi_1(f)(z)=\iota_j(y_1^j\text{diag}\bigl(\Lambda_1^{s_1^j}(f),\dots,
\Lambda_N^{s_N^j}(f),f(e^{2\pi i\theta_1^j}),\dots,
f(e^{2\pi i\theta_{D_j}^j})\bigr)y_1^{j*}), \\
\psi_1(f)(z)=\iota_j(y_2^j\text{diag}\bigl(\Lambda_1^{s_1^j}(f),\dots,
\Lambda_N^{s_N^j}(f),f(e^{2\pi i\omega_1^j}),\dots,
f(e^{2\pi i\omega_{D_j}^j})\bigr)y_2^{j*}),
\end{gather*}
for some unitaries $y_1^j,y_2^j\in M_{e_j}$ and numbers
$\theta_1^j,\dots,\theta_{D_j}^j,\omega_1^j,\dots,\omega_{D_j}^j\in\R$.
By changing $y_1^j$ and $y_2^j$ we may by (\ref{dettreq2}) and Lemma
\ref{ordert} assume that
\begin{gather*}
\theta_1^j\le\theta_2^j\le\dots\le\theta_{D_j}^j\le\theta_1^j+1, \\
\omega_1^j\le\omega_2^j\le\dots\le\omega_{D_j}^j\le\omega_1^j+1,
\end{gather*}
and
\begin{equation}\label{Treq}
\sum_{r=1}^{D_j} (\theta_r^j-\omega_r^j)=Tr(\Lambda_j(b)).
\end{equation}
Let $I$ be a $2q$-arc. By (iii),
\begin{align*}
&\#\{r: e^{2\pi i\theta_r^j}\in I\}n+\sum_{\{i:x_i\in I\}} s_i^jd_i \\
\le\
&Tr(\Lambda_j\circ\varphi_1(g_I^A\otimes 1)) \\
<\
&e_j\delta+Tr(\Lambda_j\circ\psi_1(g_I^A\otimes 1)) \\
\le\
&e_j\delta + \#\{r: e^{2\pi i\omega_r^j} \in I\pm\frac 1{4q}\}n+
\sum_{\{i:x_i\in \text{supp}\, g_I^A\}} s_i^jd_i \\
\le\
&\#\{r:e^{2\pi i\omega_r^j}\in I\pm\frac 1{2q}\}n+
\sum_{\{i:x_i\in \text{supp}\, g_I^A\}} s_i^jd_i.
\end{align*}
The last inequality uses (iv) and that $\|f_A^K\|_{\infty}\le 1$ for some
$4q$-arc $K$. Hence
$$
\#\{r: e^{2\pi i\theta_r^j}\in I\}\le
\#\{r:e^{2\pi i\omega_r^j}\in I\pm\frac 1{2q}\}.
$$
Therefore by Lemma \ref{mar},
$$
R_{D_j}\bigl((e^{2\pi i\theta_1^j},e^{2\pi i\theta_2^j},\dots,
e^{2\pi i\theta_{D_j}^j}),(e^{2\pi i\omega_1^j},e^{2\pi i\omega_2^j},\dots,
e^{2\pi i\omega_{D_j}^j})\bigr)\le \frac 1{2q}+\frac 1{2q}\le\frac 1q.
$$
By (ii), if $J$ is a $p$-arc then
$$
\#\{r:e^{2\pi i\omega_r^j}\in J\}\ge 2\frac{e_j}q,
$$
since $\|f_A^J\|_{\infty}\le\frac 1n$. Clearly $\frac {D_j}q\le\frac{e_j}q$.
Furthermore,
$$
|\sum_{r=1}^{D_j} (\theta_r^j-\omega_r^j)|=|Tr(\Lambda_j(b))|\le e_j\|b\|<
\frac {e_j}q.
$$
By Lemma \ref{det} it follows that
$$
|\theta_r^j-\omega_r^j|\le\frac 1q+\frac 2p\le \frac 3p,\quad r=1,2,\dots,D_j.
$$

Let $g_r^j:[a_j,b_j]\to\R$ be the continuous function such that
$g_r^j(a_j)=g_r^j(b_j)=\theta_r^j$, $g_r^j(t_j)=\omega_r^j$, and
such that $g_r^j$ is linear when restricted to each of the two intervals
$[a_j,t_j]$ and $[t_j,b_j]$. Note that
\begin{equation}\label{eig}
|g_r^j(t)-\theta_r^j|\le\frac 3p,\quad r=1,2,\dots,D_j.
\end{equation}
Finally, define a *-homomorphism $\xi_j:A\to C(I_j)\otimes M_m$ by
$$
\xi_j(f)(e^{2\pi it})=\iota_j(y_1^j\text{diag}\bigl(\Lambda_1^{s_1^j}(f),
\dots,\Lambda_N^{s_N^j}(f),f(e^{2\pi ig_1^j(t)}),\dots,
f(e^{2\pi ig_{D_j}^j(t)})\bigr)y_1^{j*}),
$$
for $t\in [a_j,b_j]$, $f\in A$.

Define a unital *-homomorphism $\xi:A\to B$ by
$$
\xi(f)(z)=
\begin{cases}
\xi_j(f)(z), &\ z\in I_j,\ j=1,2,\dots,M, \\
\varphi_1(f)(z), &\ z\in\T\, \backslash\cup_{j=1}^M I_j.
\end{cases}
$$
Then by (\ref{eig})
$$
\|\varphi_1(f)-\xi(f)\|<\frac{\epsilon}6,\quad f\in F.
$$
Note that $f\mapsto\Lambda_j\circ\xi(f)$ and
$f\mapsto\Lambda_j\circ\psi_1(f)$ are equivalent representations of $A$
on $M_{e_j}$, $j=1,2,\dots,M$. In particular,
$\xi$ and $\psi_1$ have the same small remainders, and hence
$\xi^{\#}(q_A'(\omega_k^A))=\psi_1^{\#}(q_A'(\omega_k^A))$, $k=1,2,\dots,N$
by Theorem \ref{sr}. Let $\eta:\T\to\R$ be the continuous function
$$
\eta(e^{2\pi it})=
\begin{cases}
\frac m{e_j}\sum_{r=1}^{D_j} (g_r^j(t)-\theta_r^j) &\ t\in [a_j,b_j],\
j=1,2,\dots,M, \\
0 &\ \text{otherwise}.
\end{cases}
$$
For $z\in\T$,
\begin{align*}
Det(\xi(v^A)(z))&=Det(\varphi_1(v^A)(z))\exp(2\pi i\eta(z)) \\
&=Det(\psi_1(v^A)(z))\exp(2\pi iTr(b(z)))\exp(2\pi i\eta(z)) \\
&=Det(\psi_1(v^A)(z))\exp(2\pi i\gamma(z)),
\end{align*}
where $\gamma:\T\to\R$ is defined by $\gamma(z)=\eta(z)+Tr(b(z))$. Note that
by (\ref{Treq})
$$
\gamma(y_j)=\eta(y_j)+Tr(b(y_j))=
\frac m{e_j}\sum_{r=1}^{D_j}(\omega_r^j-\theta_r^j)+
\frac m{e_j}Tr(\Lambda_j(b))=0,\quad j=1,2,\dots,M,
$$
and
$$\|\gamma\|_{\infty}\le \|\eta\|_{\infty}+\|Tr(b(\cdot))\|_{\infty}
<\frac m{e_j}\frac 3pD_j+\frac mq\le 3\frac mp+\frac mq\le 4\frac mp.
$$

By Proposition \ref{stand}, $\varphi_1$, $\psi_1$, and $\xi$ are
approximately unitarily equivalent to $\varphi_1'$, $\psi_1'$, and $\xi'$,
respectively, where $\varphi_1',\psi_1',\xi':A\to B$ are *-homomorphisms
of the form
\begin{gather*}
\varphi_1'(f)(e^{2\pi it})=u(t)\text{diag}\bigl( \Lambda_1^{r_1}(f),\dots,
\Lambda_N^{r_N}(f),f(e^{2\pi iF_1(t)}),\dots,
f(e^{2\pi iF_L(t)})\bigr) u(t)^*, \\
\psi_1'(f)(e^{2\pi it})=v(t)\text{diag}\bigl( \Lambda_1^{r_1}(f),\dots,
\Lambda_N^{r_N}(f),f(e^{2\pi iG_1(t)}),\dots,
f(e^{2\pi iG_L(t)})\bigr) v(t)^*, \\
\xi'(f)(e^{2\pi it})=w(t)\text{diag}\bigl( \Lambda_1^{r_1}(f),\dots,
\Lambda_N^{r_N}(f),f(e^{2\pi iH_1(t)}),\dots,f(e^{2\pi iH_L(t)})\bigr) w(t)^*,
\end{gather*}
for integers $r_1,r_2,\dots,r_N$ with $0\le r_i<\frac n{d_i}$, $i=1,2,\dots,N$,
unitaries $u,v,w$ in $C[0,1]\otimes M_m$ with $u(0)=v(0)=w(0)$,
$u(1)=v(1)=w(1)$, and continuous functions $F_r,G_r,H_r: [0,1]\to\R$,
$r=1,2,\dots,L$, such that for $t\in [0,1]$,
\begin{gather*}
F_1(t)\le F_2(t)\le\dots\le F_L(t)\le F_1(t)+1, \\
G_1(t)\le G_2(t)\le\dots\le G_L(t)\le G_1(t)+1, \\
H_1(t)\le H_2(t)\le\dots\le H_L(t)\le H_1(t)+1,
\end{gather*}
and such that for each $t\in [0,1]$,
\begin{equation}\label{GH}
\gamma(e^{2\pi it})=\sum_{r=1}^L (H_r(t)-G_r(t)).
\end{equation}
Hence
\begin{equation}\label{GHineq}
|\sum_{r=1}^L (H_r(t)-G_r(t))|< 4\frac mp.
\end{equation}

It follows from (\ref{eig}) that for each $t\in [0,1]$,
$$
R_L\bigl((e^{2\pi iF_1(t)},\dots,e^{2\pi iF_L(t)}),
(e^{2\pi iH_1(t)},\dots,e^{2\pi iH_L(t)})\bigr)\le\frac 3p.
$$
Let $t\in [0,1]$ and let $I$ be a $2q$-arc. Then by (iii) and (iv)
\begin{align*}
&\#\{r:e^{2\pi iF_r(t)}\in I\}n+\sum_{\{i:x_i\in I\}} r_id_i \\
\le\
&Tr(\varphi_1'(g_I^A\otimes 1)(e^{2\pi it})) \\
<\
&m\delta+Tr(\psi_1'(g_I^A\otimes 1)(e^{2\pi it})) \\
\le\
&m\delta + \#\{r:e^{2\pi iG_r(t)}\in I\pm\frac 1{4q}\}n+
\sum_{\{i:x_i\in \text{supp}\, g_I^A\}} r_id_i \\
\le\
&\#\{r:e^{2\pi iG_r(t)}\in I\pm\frac 1{2q}\}n+
\sum_{\{i:x_i\in \text{supp}\, g_I^A\}} r_id_i.
\end{align*}
Hence
$$
\#\{r:e^{2\pi iF_r(t)}\in I\}\le \#\{r:e^{2\pi iG_r(t)}\in I\pm\frac 1{2q}\}.
$$
It follows from Lemma \ref{mar} that for each $t\in [0,1]$,
$$
R_L\bigr((e^{2\pi iF_1(t)},\dots,e^{2\pi iF_L(t)}),
(e^{2\pi iG_1(t)},\dots,e^{2\pi iG_L(t)})\bigr)\le
\frac 1{2q}+\frac 1{2q}=\frac 1q.
$$
We conclude that
$$
R_L\bigr((e^{2\pi iG_1(t)},\dots,e^{2\pi iG_L(t)}),
(e^{2\pi iH_1(t)},\dots,e^{2\pi iH_L(t)})\bigr)\le
\frac 1q+\frac 3p\le\frac 4p.
$$

Since $f\mapsto\psi_1'(f)(y_j)$ and $f\mapsto\xi'(f)(y_j)$ are
equivalent representations of $A$ on $M_m$ for $j=1,2,\dots,M$, it follows that
$$
(e^{2\pi iG_1(t_j)},\dots,e^{2\pi iG_L(t_j)})=
(e^{2\pi iH_1(t_j)},\dots,e^{2\pi iH_L(t_j)})
$$
as unordered $L$-tuples. Therefore, as $\gamma(y_j)=0$, $j=1,2,\dots,M$, we
see by Lemma \ref{glue} and (\ref{GH}) that
$$
G_r(t_j)=H_r(t_j),\quad r=1,2,\dots,L,\ j=1,2,\dots,M.
$$
As $v(0)=w(0)$, $v(1)=w(1)$, we may thus define a *-homomorphism
$\mu:A\to B$ by
$$
\mu(f)(e^{2\pi it})=v(t)\text{diag}\bigl( \Lambda_1^{r_1}(f),\dots,
\Lambda_N^{r_N}(f),f(e^{2\pi iH_1(t)}),\dots,f(e^{2\pi iH_L(t)})\bigr) v(t)^*,
$$
for $f\in A$, $t\in [0,1]$. Since
$$
Tr(\mu(f)(z))=Tr(\xi'(f)(z))=Tr(\xi(f)(z)),\quad z\in\T,\ f\in A,
$$
we get from \cite[Theorem 1.4]{TATD} that $\mu$ and $\xi$ are approximately
unitarily equivalent.

By (i) we have that for every $l$-arc $J$,
$$
\#\{r:e^{2\pi iG_r(t)}\in J\} > 8\frac mp.
$$
As $L\frac 4p\le 4\frac mp$, we conclude from (\ref{GHineq}) and Lemma
\ref{det} that
$$
|G_r(t)-H_r(t)|\le\frac 4p+\frac 2l\le\frac 6l.
$$
Hence
$$
\|\mu(f)-\psi_1'(f)\|<\frac{\epsilon}6,\quad f\in F.
$$

Choose unitaries $U,V\in B$ such that
\begin{gather*}
\|\xi(f)-U\mu(f)U^*\|< \frac{\epsilon}6,\quad f\in F, \\
\|\psi_1'(f)-V\psi_1(f)V^*\|< \frac{\epsilon}6,\quad f\in F.
\end{gather*}
Set $W=UV$. Then for $f\in F$,
\begin{align*}
&\|\varphi(f)-W\psi(f)W^*\| \\
\le\ &\|\varphi(f)-\varphi_1(f)\|+ \|\varphi_1(f)-\xi(f)\|+
\|\xi(f)-U\mu(f)U^*\|+ \\
&\|U\mu(f)U^*-U\psi_1'(f)U^*\|+
\|U\psi_1'(f)U^*-UV\psi_1(f)V^*U^*\|+ \\
&\|W\psi_1(f)W^*-W\psi(f)W^*\| \\
<\ &\frac{\epsilon}6+\frac{\epsilon}6+\frac{\epsilon}6+
\frac{\epsilon}6+\frac{\epsilon}6+\frac{\epsilon}6=\epsilon.
\end{align*}
\end{proof}
\end{thm}

\begin{lemma}\label{cutdu}
Let $B=A(m,e_1,e_2,\dots,e_M)$ be a building block and let $r\in B$ be a
non-zero projection of rank $s\in\Z$. Let $C=rBr$ and let $u,v\in C$ be
unitaries. Then
$$
D_C\bigl(q_C'(u),q_C'(v)\bigr)\le
2\pi\frac ms D_B\bigl(q_B'(u+(1-r)),q_B'(v+(1-r))\bigr).
$$

\begin{proof}
Let $\epsilon=D_B\bigl(q_B'(u+(1-r)),q_B'(v+(1-r))\bigr)$. We may assume that
$\epsilon<1$. Let $b\in B$ be a self-adjoint element such that
$uv^*+(1-r)=e^{2\pi ib}$ modulo $\overline{DU(B)}$ and $\|b\|\le\epsilon$.
Define $c\in C$ by $c(z)=\frac 1s Tr(b(z))r$. Since $\widehat b=\widehat c$
we have that $e^{2\pi ib}=e^{2\pi ic}$ modulo $\overline{DU(B)}$. Thus
$$
uv^*+(1-r)=e^{2\pi ic}=re^{2\pi ic}r+(1-r)\quad\text{modulo}\ \overline{DU(B)}.
$$
$C$ is a building block by Corollary \ref{cutdown}, and therefore it
follows from Proposition \ref{detunitgr} that
$uv^*=re^{2\pi ic}r$ modulo $\overline{DU(C)}$. Thus
$$
D_C(q_C'(u),q_C'(v))\le \|re^{2\pi ic}r-r\|\le
\|e^{2\pi ic}-1\|\le 2\pi \|c\|\le 2\pi\frac ms\epsilon.
$$
\end{proof}
\end{lemma}

Let $A=A_1\oplus A_2\oplus\dots\oplus A_R$, where
$A_i=A(n_i,d_1^i,d_2^i,\dots,d_{N_i}^i)$ is a building block. For each
$i=1,2,\dots,R$, we define unitaries in $A$ by
\begin{align*}
v^A_i&=(1,\dots,1,v^{A_i},1,\dots,1). \\
w_{i,k}^A&=(1,\dots,1,w_k^{A_i},1,\dots,1),\quad k=1,2,\dots,N_i.
\end{align*}
Set $U^A=\cup_{i=1}^R \{w_{i,k}^A:k=1,2,\dots,N_i\}$.
If $p$ is a positive integer, we set
\begin{gather*}
H(A,p)=\cup_{i=1}^R \iota_i(H(A_i,p)), \\
\widetilde{H}(A,p)=\cup_{i=1}^R \iota_i(\widetilde{H}(A_i,p)),
\end{gather*}
where $\iota_i:A_i\to A$ denotes the inclusion, $i=1,2,\dots,R$.

\begin{thm}\label{un}
Let $A=A_1\oplus A_2\oplus\dots\oplus A_R$ be a finite direct sum of building
blocks. Let $p_1,p_2,\dots,p_R$ be the minimal non-zero central projections in
$A$. Let $\epsilon>0$ and let $F\subseteq A$ be a finite set. There exists a
positive integer $l$ such that if $p$ and $q$ are positive integers
with $l\le p\le q$, if $B$ is a finite
direct sum of building blocks, if $\varphi,\psi:A\to B$ are unital
*-homomorphisms, if $\delta>0$, if
\begin{itemize}
\item[(i)]\
$\widehat{\psi}(\widehat{h})>\frac 8p,\quad h\in H(A,l)$,
\vspace{1 mm}
\item[(ii)]\
$\widehat{\psi}(\widehat{h})>\frac 2q,\quad h\in H(A,p)\cup
\{p_1,p_2,\dots,p_R\}$,
\vspace{1 mm}
\item[(iii)]\
$\|\widehat{\varphi}(\widehat{h})-\widehat{\psi}(\widehat{h})\|<\delta,
\quad h\in \widetilde{H}(A,2q)$,
\vspace{1 mm}
\item[(iv)]\
$\widehat{\psi}(\widehat{h})>\delta,\quad h\in H(A,4q)$,
\vspace{1 mm}
\item[(v)]\
$D_B\bigl(\, \varphi^{\#}(q_A'(v^A_i))\, ,\, \psi^{\#}(q_A'(v^A_i))\, \bigr)<
\frac 1{4q^2},\quad i=1,2,\dots,R$;
\vspace{1 mm}
\end{itemize}
and if at least one of the two statements
\begin{itemize}
\item[(vi)]\
$[\varphi]=[\psi]\quad\text{in}\ KK(A,B)$,
\vspace{1 mm}
\item[(vii)]\
$\varphi_*=\psi_*$\ \ on $K_0(A)$\quad and \quad
$\varphi^{\#}(x)=\psi^{\#}(x)$,\ \ $x\in U^A$,
\end{itemize}
are true; then there exists a unitary $W\in B$ such that
$$
\|\varphi(f)-W\psi(f)W^*\|<\epsilon,\quad f\in F.
$$

\begin{proof}
For each $i=1,2,\dots,R$, let $\iota_i:A_i\to A$ be the inclusion and let
$\pi_i:A\to A_i$ be the projection. Choose by Theorem \ref{uns} a positive
integer $l_0^i$ with respect to the finite set $\pi_i(F)\subseteq A_i$ and
$\epsilon>0$. Set $l=\max_i l_0^i$.

Let integers $q\ge p\ge l$, a finite direct sum of building blocks
$B$, and unital *-homomorphisms $\varphi,\psi:A\to B$ be
given such that the above holds. Since (vi) implies (vii) by
Proposition \ref{KKunid}, we may assume that (vii) holds. It is easy to reduce
to the case where $B=A(m,e_1,e_2,\dots,e_M)$ is a single building block.

Since $\varphi_*[p_i]=\psi_*[p_i]$ in $K_0(B)$
for $i=1,2,\dots,R$, there is by Lemma \ref{cd} a unitary $u\in B$
such that $u\varphi(p_i)u^*=\psi(p_i)$ for every $i=1,2,\dots,R$.
Hence we may assume that $\varphi(p_i)=\psi(p_i)$, $i=1,2,\dots,R$.
Set $q_i=\psi(p_i)$. It follows from (ii) that $q_i\neq 0$, $i=1,2,\dots,R$.
Let $t_i$ be the (normalized) trace of $q_i$.

Let $\varphi_i,\psi_i:A_i\to q_iBq_i$ be the induced maps. Note that $q_iBq_i$
is a building block by Corollary \ref{cutdown}. Fix some $i=1,2,\dots,R$.

Every tracial state on $q_iBq_i$ is of the form
$\frac 1{t_i}\omega|_{q_iBq_i}$ for some $\omega\in T(B)$. Therefore
$\varphi_i$ and $\psi_i$ satisfy (i)-(iv) of Theorem \ref{uns}, with
$\delta$ replaced by $\frac{\delta}{t_i}$. Note that $t_i>\frac 2q$ by (ii).
Since
$$
D_B\bigl(q_B'(\varphi_i(v^{A_i})+(1-q_i)),q_B'(\psi_i(v^{A_i})+(1-q_i))\bigr)<
\frac 1{4q^2}
$$
by (vi), we have that
$$
D_{q_iBq_i}\bigl(\varphi_i^{\#}(q_{A_i}'(v^{A_i})),
\psi_i^{\#}(q_{A_i}'(v^{A_i}))\bigr)\le 2\pi\frac 1{t_i}\frac 1{4q^2}
<2\pi\frac q2\frac 1{4q^2}<\frac 1q
$$
by Lemma \ref{cutdu}, which is (v) of Theorem \ref{uns} for $\varphi_i$ and
$\psi_i$. Similarly we get that
$\varphi_i^{\#}(w_k^{A_i})=\psi_i^{\#}(w_k^{A_i})$, $k=1,2,\dots,N_i$, which
is (vi) of Theorem \ref{uns}. Hence there exists a unitary $W_i\in q_iBq_i$
such that
$$
\|\varphi_i(f)-W_i\psi_i(f){W_i}^*\|<\epsilon,\quad f\in\pi_i(F).
$$

Set $W=\sum_{i=1}^R W_i$. Then $W\in B$ is a unitary and
$$
\|\varphi(f)-W\psi(f)W^*\|<\epsilon,\quad f\in F.
$$
\end{proof}
\end{thm}

\section{Existence}

The goal of this section is to prove an existence theorem that is the
counterpart of the uniqueness theorem of the previous section.

Let $A$ and $B$ be building blocks and let $\varphi:A\to B$ be a
*-homomorphism. We say that continuous functions
$\lambda_1,\lambda_2,\dots,\lambda_N:\T\to\T$ are
eigenvalue functions for $\varphi$ if
$\lambda_1(z),\lambda_2(z),\dots,\lambda_N(z)$ are eigenvalues
for the matrix $\varphi(\iota\otimes 1)(z)$ (counting multiplicities) for
every $z\in\T$, where $\iota:\T\to\C$ denotes the inclusion.

\begin{thm}\label{ex1}
Let $A=A(n,d_1,d_2,\dots,d_N)$ be a building block, let $\epsilon>0$,
and let $C$ be a positive integer.
There exists a positive integer $K$ such that if
\begin{enumerate}
\item[(i)]
$B=A(m,e_1,e_2,\dots,e_M)$ is a building block with $s(B)\ge K$,
\item[(ii)]
$\kappa\in KK(A,B)_e$,
\item[(iii)]
$\lambda_1,\lambda_2,\dots,\lambda_C:\T\to\T$ are continuous functions,
\item[(iv)]
$u\in B$ is a unitary such that $\kappa_*[v^A]=[u]$ in $K_1(B)$;
\end{enumerate}
then there exists a unital *-homomorphism $\varphi : A\to B$ such that
$\lambda_1,\lambda_2,\dots,\lambda_C$ are eigenvalue functions for
$\varphi$, and such that
\begin{gather*}
[\varphi]=\kappa\quad \text{in}\ KK(A,B), \\
\varphi^{\#}(q_A'(v^A))=q_B'(u),\quad \text{in}\ U(B)/\overline{DU(B)},\\
\| \widehat{\varphi}(f)-\frac 1C\sum_{k=1}^C f\circ\lambda_k\|<\epsilon\|f\|,
\quad f\in\text{Aff}\, T(A),
\end{gather*}
when we identify $\text{Aff}\, T(A)$ and $\text{Aff}\, T(B)$ with
$C_{\R}(\T)$ as order unit spaces.

\begin{proof}
We may assume that $\epsilon<4$ and, by repeating the functions
$\lambda_1,\lambda_2,\dots,\lambda_C$, that $C>\frac{8}{\epsilon}$.
Let $K$ be a positive integer such that
$$
K\ge\frac{4(N+C+2)n}{\epsilon}.
$$
Let $B$, $\lambda_1,\lambda_2,\dots,\lambda_C$, $\kappa$, and $u$ be as above.
By Proposition \ref{sf} there are integers $h_{ji}$, $i=1,2,\dots,N$,
$j=1,2,\dots,M$, with $0\le h_{ji}<\frac{n}{d_i}$ for $i\neq N$, such that
$$
\begin{pmatrix}
\kappa^*([\Lambda_1^B]) \\
\kappa^*([\Lambda_2^B]) \\
\vdots \\
\kappa^*([\Lambda_M^B])
\end{pmatrix}
=
\begin{pmatrix}
h_{11} & h_{12} & \hdots & h_{1N} \\
h_{21} & h_{22} & \hdots & h_{2N} \\
\vdots & \vdots &        & \vdots \\
h_{M1} & h_{M2} & \hdots & h_{MN}
\end{pmatrix}
\begin{pmatrix}
{[\Lambda_1^A]} \\
{[\Lambda_2^A]} \\
\vdots \\
{[\Lambda_N^A]}
\end{pmatrix}.
$$
As in the proof of Theorem \ref{KKliftc} we see that
\begin{equation}\label{hji1}
\sum_{i=1}^N h_{ji}d_i=e_j,
\end{equation}
since $\kappa_*:K_0(A)\to K_0(B)$ preserves the order unit, and
$h_{jN}>\frac n{d_N}$, because $s(B)\ge Nn$.
By Proposition \ref{sf} we have for $i=1,2,\dots,N$, $j=1,2,\dots,M$,
\begin{equation}\label{eqs_i}
\frac{m}{e_j}h_{ji}=l_{ji}\frac{n}{d_i}+s_i,
\end{equation}
where $l_{ji}$ and $s_i$ are integers such that $0\le s_i<\frac n{d_i}$.
Note that $l_{ji}\ge 0$. For $j=1,2,\dots,M$, choose integers
$h_{jN}^o$, $0\le h_{jN}^o<\frac{n}{d_N}$, and $r_j\ge 0$
such that
\begin{equation}\label{eqh_{jN}^o}
h_{jN}=r_j\frac{n}{d_N}+h_{jN}^o,
\end{equation}
and note that
\begin{equation}\label{eqs_N}
\frac{m}{e_j}h_{jN}^o=l_{jN}^o\frac{n}{d_N}+s_N
\end{equation}
for some integers $l_{jN}^o\ge 0$, $j=1,2,\dots,M$. Then
\begin{equation}\label{eql_{jN}}
l_{jN}-l_{jN}^o=\frac{m}{e_j}r_j.
\end{equation}
Let for each $j$
\begin{equation}\label{eqr_j}
r_j=k_j(C+2)+u_j
\end{equation}
for some integers $k_j\ge 0$ and $0\le u_j<C+2$ and set
$$
b=\min_{1\le j\le M} k_j\frac{m}{e_j}.
$$
Note that for $j=1,2,\dots,M$,
\begin{align*}
e_j&=\sum_{i=1}^N h_{ji}d_i < (N-1)n + r_jn + h_{jN}^o d_N < Nn+r_jn \\
   &= (N+C+2)n+(r_j-(C+2))n\le \frac{\epsilon}4 e_j + (r_j-(C+2))n.
\end{align*}
Hence
$$
(1-\frac{\epsilon}4)e_j< (r_j-(C+2))n.
$$
Therefore
\begin{equation}\label{eqC}
nk_j(C+2)\frac{m}{e_j}=n(r_j-u_j)\frac{m}{e_j}
> n(r_j-(C+2))\, \frac{m}{e_j}>(1-\frac{\epsilon}4)m.
\end{equation}
Since by (\ref{hji1})
$$
nk_j(C+2)\frac{m}{e_j}\le nr_j\frac{m}{e_j}\le h_{jN}d_N\frac{m}{e_j}\le m,
$$
we see that
\begin{equation}\label{exeq}
nb\frac{8}{\epsilon}\le nbC\le nb(C+2)\le m.
\end{equation}
By this and (\ref{eqC}),
$$
m(1-\frac{\epsilon}4)<nb(C+2)\le nbC+\frac{\epsilon}4m.
$$
Hence from (\ref{exeq}) we conclude that
$$
0\le 1-\frac{nbC}m< \frac{\epsilon}2.
$$

Let $x_1,x_2,\dots,x_N$ denote the exceptional points of $A$ and let
$y_1,y_2,\dots,y_M$ be those of $B$. Set for $j=1,2,\dots,M$,
$$
a_j=\bigl(\prod_{i=1}^{N-1} Det(\Lambda_i(v^A))^{h_{ji}}\bigr)
Det(\Lambda_N(v^A))^{h_{jN}^o},
$$
then by (\ref{eqs_i}) and (\ref{eqs_N})
\begin{equation}\label{eqaj}
{a_j}^{\frac m{e_j}}=\bigl(\prod_{i=1}^{N-1} {x_i}^{l_{ji}}\bigr)
{x_N}^{l_{jN}^o}\prod_{i=1}^N Det(\Lambda_i(v^A))^{s_i}.
\end{equation}
Set for $j=1,2,\dots,M$,
$$
c_j=Det(\Lambda_j(u))
$$
and note that
\begin{equation}\label{eqcj}
{c_j}^{\frac m{e_j}}=Det(u(y_j)).
\end{equation}
By (\ref{eqC}) we see that $k_j\neq 0$, $j=1,2,\dots,M$, and hence there
exists a continuous function $\lambda_{C+1}:\T\to\T$ such that
\begin{equation}\label{eqeta}
(\lambda_{C+1}(y_j))^{-k_j}=a_j c_j^{-1}
\prod_{k=1}^C (\lambda_k(y_j))^{k_j},\quad j=1,2,\dots,M.
\end{equation}

Let for $f\in A$ and $j=1,2,\dots,M$, $D_j(f)$ be the $m\times m$ matrix
\begin{align*}
\Bigl(&\Lambda_1^{s_1}(f),\dots,\Lambda_N^{s_N}(f),
\underbrace{f(x_1),\dots,f(x_1)}_{l_{j1}\text{\ times}},\dots,
\underbrace{f(x_{N-1}),\dots,f(x_{N-1})}_{l_{j(N-1)}\text{\ times}}, \\
&\underbrace{f(x_N),\dots,f(x_N)}_{l_{jN}^o\text{\ times}},
\underbrace{f(\lambda_1(y_j)),\dots,f(\lambda_1(y_j))}_
{k_j\frac{m}{e_j}-b\text{\ times}},\dots, \\
&\underbrace{f(\lambda_{C+1}(y_j)),\dots,f(\lambda_{C+1}(y_j))}_
{k_j\frac{m}{e_j}-b\text{\ times}},
\underbrace{f(1),\dots,f(1)}_
{(k_j+u_j)\frac{m}{e_j}-b\text{\ times}}, \\
&\underbrace{f(\lambda_1(y_j)),\dots,f(\lambda_1(y_j))}_
{b\text{\ times}},\dots,
\underbrace{f(\lambda_{C+1}(y_j)),\dots,
f(\lambda_{C+1}(y_j))}_{b\text{\ times}}, \\
&\underbrace{f(1),\dots,f(1)}_
{b\text{\ times}}\, \Bigr).
\end{align*}
Since $D_j(f)$ is a block-diagonal matrix with $\frac m{e_j}h_{ji}$ blocks
of the form $\Lambda_i(f)$, $i=1,2,\dots,N-1$, $\frac m{e_j}h_{jN}^o$ blocks
of the form $\Lambda_N(f)$, $k_j\frac m{e_j}$ blocks of the form
$f(\lambda_k(y_j))$, $k=1,2,\dots,C+1$, and $(k_j+u_j)\frac m{e_j}$ blocks
of the form $f(1)$, there exists a unitary $W_j\in M_m$ such that
$$
W_jD_j(f)W_j^*\in M_{e_j}\subseteq M_m
$$
for every $f\in A$. Set $L=\frac 1n(m-\sum_{i=1}^N s_id_i)-(C+2)b$.
For each $j=1,2,\dots,M$, we have by (\ref{hji1}), (\ref{eqs_i}),
(\ref{eql_{jN}}), (\ref{eqr_j}) that
$$
L=\sum_{i=1}^N l_{ji} -(C+2)b
=\sum_{i=1}^{N-1} l_{ji} +l_{jN}^o +\frac m{e_j}(k_j(C+2)+u_j)-(C+2)b.
$$
Choose for $k=1,2,\dots,L$ continuous functions $\mu_k : \T \to \T$
such that for each $j=1,2,\dots,M,$
\begin{align*}
&\bigr(\mu_1(y_j),\mu_2(y_j),\dots,\mu_L(y_j)\bigl)=
\bigr(\underbrace{x_1,\dots,x_1}_{l_{j1}\ \text{times}},\dots,
\underbrace{x_{N-1},\dots,x_{N-1}}_{l_{j(N-1)}\ \text{times}},
\underbrace{x_{N},\dots,x_{N}}_{l_{jN}^o\ \text{times}}, \\
&\underbrace{\lambda_1(y_j),\dots,
\lambda_1(y_j)}_{k_j\frac{m}{e_j}-b\ \text{times}},\dots,
\underbrace{\lambda_{C+1}(y_j),\dots,
\lambda_{C+1}(y_j)}_{k_j\frac{m}{e_j}-b\ \text{times}},
\underbrace{1,1,\dots,1}_{(k_j+u_j)\frac{m}{e_j}-b\ \text{times}}\bigl)
\end{align*}
as ordered tuples.

Choose a unitary $W\in C(\T)\otimes M_m$ such that $W(y_j)=W_j$
for $j=1,2,\dots,M$. Define a continuous function $g:\T\to \T$ such that
$$
g(z)\, \prod_{k=1}^L \mu_k(z)
\prod_{k=1}^{C+1} (\lambda_k(z))^b\prod_{i=1}^N Det(\Lambda_i(v^A))^{s_i}=
Det(u(z)),\quad z\in\T.
$$
Then by (\ref{eqaj}), (\ref{eqcj}), and (\ref{eqeta}) we have that $g(y_j)=1$
for $j=1,2,\dots,M$. Define a unital *-homomorphism $\varphi:A\to B$ by
\begin{align*}
\varphi(f)(z)=W(z)\, \text{diag}\, \Bigl(\, &\Lambda_1^{s_1}(f),
\dots,\Lambda_N^{s_N}(f),f(\mu_1(z)),\dots,f(\mu_L(z)), \\
&\underbrace{f(\lambda_1(z)),\dots,f(\lambda_1(z))}_{b\ \text{times}},\dots,
\underbrace{f(\lambda_{C+1}(z)),\dots,
f(\lambda_{C+1}(z))}_{b\ \text{times}},\\
&f(g(z)),\underbrace{f(1),\dots,f(1)}_{b-1\ \text{times}}\, \Bigr)W(z)^*.
\end{align*}
By the remarks following the definition of $D_j(f)$ we see that
for $j=1,2,\dots,M$,
\begin{align*}
\varphi^*[\Lambda_j^B]&=\sum_{i=1}^{N-1} h_{ji}[\Lambda_i^A]+
h_{jN}^o[\Lambda_N^A]+(k_j(C+1)+(k_j+u_j))\frac n{d_N}[\Lambda_N^A] \\
&=\sum_{i=1}^{N-1} h_{ji}[\Lambda_i^A]+(h_{jN}^o+r_j\frac n{d_N})[\Lambda_N^A]=
\sum_{i=1}^N h_{ji}[\Lambda_i^A]=\kappa^*[\Lambda_j^B],
\end{align*}
and hence $\varphi^*=\kappa^*$ in $Hom(K^0(B),K^0(A))$.
Furthermore, as $Det(v^A(\cdot))$ is the identity map on $\T$, we have that
for $z\in\T$,
$$
Det(\varphi(v^A)(z))=\prod_{i=1}^N Det(\Lambda_i(v^A))^{s_i}
\prod_{k=1}^L \mu_k(z)\bigl(\prod_{k=1}^{C+1}\lambda_k(z)^{b}\bigr) g(z)=
Det(u(z)),
$$
and by (\ref{eqeta}), for $j=1,2,\dots,M$,
\begin{align*}
Det(\Lambda_j\circ\varphi(v^A))&=
\bigl(\prod_{i=1}^{N-1} Det(\Lambda_i(v^A))^{h_{ji}}\bigr)
Det(\Lambda_N(v^A))^{h_{jN}^o}\prod_{k=1}^{C+1} \lambda_k(y_j)^{k_j} \\
&=\bigl(\prod_{i=1}^{N-1} Det(\Lambda_i(v^A))^{h_{ji}}\bigr)
Det(\Lambda_N(v^A))^{h_{jN}^o}
{a_j}^{-1}c_j=Det(\Lambda_j(u)).
\end{align*}
Hence $q_B'(\varphi(v^A))=q_B'(u)$ in $U(B)/\overline{DU(B)}$ by Proposition
\ref{detunitgr}. It follows from Theorem \ref{KKliftc} that $[\varphi]=\kappa$
in $KK(A,B)$.

Finally, for $\omega\in T(B)$ and $f\in\text{Aff}\, T(A)\cong C_{\R}(\T)$,
\begin{align*}
&|\widehat{\varphi}(f)(\omega)-\frac 1C\sum_{i=1}^C f\circ\lambda_k(\omega)|
=|\omega(\varphi(f\otimes 1))-
\frac 1C\sum_{k=1}^C \omega((f\circ\lambda_k)\otimes 1)| \\
\le\ & |\frac 1m(m-Cbn)|\, \|f\|+\|\frac 1mbn\sum_{k=1}^C f\circ\lambda_k-
\frac 1C\sum_{k=1}^C f\circ\lambda_k\| \\
\le\ &|\frac 1m(m-Cbn)|\, \|f\|+|\frac 1mbn-\frac 1C|C\, \|f\|
=2|1-\frac{Cbn}m|\, \|f\|<\epsilon \|f\|. \\
\end{align*}
Hence
$$
\|\widehat{\varphi}(f)-\frac 1C\sum_{k=1}^C f\circ\lambda_k\|<\epsilon \|f\|.
$$
\end{proof}
\end{thm}

The following result is due to Li \cite[Theorem 2.1]{LIPHD}. It generalizes
a theorem of Thomsen \cite[Theorem 2.1]{TAIT} and it is the key stone in the
proof of Theorem \ref{ex2} below.

\begin{thm}\label{li}
Let $X$ be a path-connected compact metric space, let
$F\subseteq C_{\R}(X)$ be a finite subset and let $\epsilon>0$.
There exists a positive integer $L$ such that for all integers $N\ge L$, for
all compact metric spaces $Y$, and for all positive linear order unit
preserving maps $\Theta:C_{\R}(X)\to C_{\R}(Y)$, there exist continuous
functions $\lambda_k:Y\to X$, $k=1,2,\dots,N$, such that
$$
\|\Theta(f)-\frac 1N \sum_{k=1}^{N} f\circ \lambda_k\|<\epsilon , \quad f\in F.
$$
\end{thm}

\begin{thm}\label{ex2}
Let $A=A(n,d_1,d_2,\dots,d_N)$ be a building block, let $\epsilon>0$,
let $F\subseteq\text{Aff}\, T(A)$ be a finite set, and let $C$ be a
non-negative integer. There exists a positive integer $K$ such that if
\begin{enumerate}
\item[(i)]
$B=A(m,e_1,e_2,\dots,e_M)$ is a building block with $s(B)\ge K$,
\item[(ii)]
$\kappa\in KK(A,B)_e$,
\item[(iii)]
$\lambda_1,\lambda_2,\dots,\lambda_C:\T\to\T$ are continuous functions,
\item[(iv)]
$\Theta : \text{Aff}\, T(A)\to \text{Aff}\, T(B)$ is a positive linear
order unit preserving map,
\item[(v)]
$u\in B$ is a unitary such that $\kappa_*[v^A]=[u]$ in $K_1(B)$;
\end{enumerate}
then there exists a unital *-homomorphism $\varphi : A\to B$ such that
$\lambda_1,\lambda_2,\dots,\lambda_C$ are eigenvalue functions for
$\varphi$ and such that
\begin{gather*}
[\varphi]=\kappa\quad\text{in}\ KK(A,B), \\
\varphi^{\#}(q_A'(v^A))=q_B'(u)\quad\text{in}\ U(B)/\overline{DU(B)} \\
\|\widehat{\varphi}(f)-\Theta(f)\|<\epsilon,\quad f\in F.
\end{gather*}

\begin{proof}
We may assume that $\|f\|\le 1$, $f\in F$.
Choose by Theorem \ref{li} an integer $L$ with respect to
$F\subseteq\text{Aff}\, T(A)\cong C_{\R}(\T)$ and $\frac{\epsilon}3$.
We may assume that $L>C$ and that $1-\frac{L-C}{C+L}<\frac{\epsilon}3$.
Then choose by Theorem \ref{ex1} an integer $K$ with respect to $C+L$ and
$\frac{\epsilon}3$.

Now let $B$, $\Theta$, $\lambda_1,\lambda_2,\dots,\lambda_C$, $\kappa$
and $u$ be given as above. Choose continuous functions
$\lambda_{C+1},\lambda_{C+2},\dots,\lambda_{C+L}:\T\to\T$ such that
in $\text{Aff}\, T(B)\cong C_{\R}(\T)$,
$$
\|\Theta(f)-\frac 1L\sum_{k=C+1}^{C+L} f\circ\lambda_k\|<\frac{\epsilon}3,
\quad f\in F.
$$
By Theorem \ref{ex1} there exists a unital *-homomorphism $\varphi:A\to B$
such that $\lambda_1,\lambda_2,\dots,\lambda_{C+L}$ are eigenvalue functions
for $\varphi$ and such that
\begin{gather*}
[\varphi]=\kappa\quad\text{in}\ KK(A,B), \\
\varphi^{\#}(q_A'(v^A))=q_B'(u)\quad\text{in}\ U(B)/\overline{DU(B)}, \\
\| \widehat{\varphi}(f)-\frac 1{C+L}\sum_{k=1}^{C+L} f\circ\lambda_k\|<
\frac{\epsilon}3\|f\|,\quad f\in\text{Aff}\, T(A).
\end{gather*}
Since for $f\in \text{Aff}\, T(A)$,
\begin{align*}
&\|\frac 1{C+L} \sum_{k=1}^{C+L} f\circ\lambda_k-
\frac 1L\sum_{k=C+1}^{C+L} f\circ\lambda_k\| \\
\le\ &\|\frac 1{C+L} \sum_{k=C+1}^{C+L} f\circ\lambda_k-
\frac 1L \sum_{k=C+1}^{C+L} f\circ\lambda_k\|+
\|\frac 1{C+L} \sum_{k=1}^C f\circ\lambda_k\| \\
\le\ &|\frac 1{C+L} -\frac 1L|\, L\, \|f\|+\frac 1{C+L} C\, \|f\|=
(1-\frac{L-C}{C+L})\, \|f\|<\frac{\epsilon}3\, \|f\|,
\end{align*}
we get that
$$
\|\widehat{\varphi}(f)-\Theta(f)\|<\epsilon,\quad f\in F.
$$
\end{proof}
\end{thm}

\begin{lemma}\label{unitgrpAp}
Let $A=A(n,d_1,d_2,\dots,d_N)$ be a building block, let $p\in A$ be a non-zero
projection, and let $u\in A$ be a unitary. Then there exists a unitary
$w\in pAp$ such that
$$
q_A'(u)=q_A'(w+(1-p))\quad\text{in}\ U(A)/\overline{DU(A)}.
$$

\begin{proof}
Note that $pAp$ is a building block by Corollary \ref{cutdown}. Hence
by Lemma \ref{detrange} there exists a unitary $w\in pAp$ such that
\begin{gather*}
Det(w(z))=Det(u(z)),\quad z\in\T, \\
Det(\Lambda_i(w))=Det(\Lambda_i(u)),\quad i=1,2,\dots,N.
\end{gather*}
Then $q_A'(u)=q_A'(w+(1-p))\ \text{in}\ U(A)/\overline{DU(A)}$ by Theorem
\ref{detunitgr}.
\end{proof}
\end{lemma}

\begin{thm}\label{ex}
Let $A=A_1\oplus A_2\oplus\dots\oplus A_R$ be a finite direct sum of building
blocks. Let $F\subseteq \text{Aff}\, T(A)$ be a finite set and let
$\epsilon>0$. There exists a positive integer $K$ such that if
\begin{itemize}
\item[(i)]
$B=B_1\oplus B_2\oplus\dots\oplus B_S$ is a finite direct sum of
building blocks and $\kappa$ is an element in $KK(A,B)_e$,
\item[(ii)]
for every minimal non-zero central projection $p$ in $A$ we have that
$$
s(B)\, \rho_B(\kappa_*[p])\ge K\quad\text{in}\ \text{Aff}\, T(B),
$$
\item[(iii)]
there exists a linear positive order unit preserving map
$\Theta : \text{Aff}\, T(A) \to \text{Aff}\, T(B)$
such that the diagram
$$ 
  \begin{CD}
      K_0(A) @> \rho_A >> \text{Aff}\, T(A)    \\   
      @V \kappa_* VV      @VV \Theta V \\
      K_0(B)  @>> \rho_B > \text{Aff}\, T(B)
  \end{CD}
$$
commutes,
\item[(iv)]
$u_1,u_2,\dots,u_N\in B$ are unitaries such that
$$
\kappa_*[v^A_i]=[u_i]\quad\text{in}\ K_1(A),\quad i=1,2,\dots,R;
$$
\end{itemize}
then there exists a unital *-homomorphism $\varphi:A\to B$ such that
$[\varphi]=\kappa$ in $KK(A,B)$, and such that
\begin{gather*}
\|\widehat\varphi (f)-\Theta(f)\|<\epsilon,\quad f\in F, \\
\varphi^{\#}(q_A'(v_i^A))=q_B'(u_i)\quad\text{in}\ U(B)/\overline{DU(B)},
\quad i=1,2,\dots,R.
\end{gather*}

\begin{proof}
Let $\pi_i^A : A\to A_i$ be the projection and $\iota_i^A : A_i\to A$ be the
inclusion, $i=1,2,\dots,R$. Let $p_1,p_2,\dots,p_R$ denote the minimal non-zero
central projections in $A$. Choose by Theorem \ref{ex2} an integer $K_i$ with
respect to $\widehat{\pi_i^A}(F)\subseteq \text{Aff}\, T(A_i)$,
$\epsilon>0$ and $C=0$. Set $K=\max_{1\le i\le R} K_i$.

Let $B$, $\kappa$, $\Theta$, and $u_1,u_2,\dots,u_N$ be as above. We may assume
that $S=1$. To see this, assume that the case $S=1$ has been settled. Let
$\pi_l^B:B\to B_l$ be the projection and let $\iota_l^B:B_l\to B$ be the
inclusion. As the diagram
$$ 
  \begin{CD}
      K_0(A) @> \rho_A >> \text{Aff}\, T(A)    \\   
      @V {\pi^B_l}_*\circ\kappa_* VV      @VV \widehat{\pi_l^B}\circ\Theta V \\
      K_0(B_l)  @>> \rho_{B_l} > \text{Aff}\, T(B_l)
  \end{CD}
$$
commutes for $l=1,2,\dots,S$, and since
$s(B_l)\, \rho_{B_l}({\pi^B_l}_*\circ\kappa_*[p_i])\ge K$
for $i=1,2,\dots,R,\ l=1,2,\dots,S$, we get unital *-homomorphisms
$\varphi_l : A\to B_l$ such that
\begin{gather*}
[\varphi_l]=[\pi_l^B]\cdot\kappa\qquad\text{in\ } KK(A,B_l), \\
\|\widehat{\varphi_l}(f)-\widehat{\pi_l^B}\circ\Theta(f)\|<\epsilon,
\quad f\in F, \\
\varphi_l^{\#}(q_A'(v_i^A))=q_{B_l}'(\pi_l^B(u_i))
\quad\text{in}\ U(B_l)/\overline{DU(B_l)},\quad i=1,2,\dots,R.
\end{gather*}
Define $\varphi : A\to B$ by
$\varphi(a)=(\varphi_1(a),\varphi_2(a),\dots,\varphi_S(a))$. Then
\begin{gather*}
[\varphi]=[\sum_{l=1}^S \iota_l^B\circ\varphi_l]=
\sum_{l=1}^S [\iota_l^B]\cdot [\pi_l^B]\cdot\kappa=
\kappa\qquad\text{in\ } KK(A,B), \\
\|\widehat\varphi(f)-\Theta(f)\|=
\max_l \|\widehat{\pi_l^B}\circ\widehat\varphi(f)-
\widehat{\pi_l^B}\circ\Theta(f)\|<\epsilon,\quad f\in F, \\
\varphi^{\#}(q_A'(v_i^A))=q_B'(u_i)\quad\text{in}\ U(B)/\overline{DU(B)},
\quad i=1,2,\dots,R.
\end{gather*}

So assume $B=A(m,e_1,e_2,\dots,e_M)$. Note that by assumption
$\kappa_*[p_i]>0$ in $K_0(B)$ for $i=1,2,\dots,R$.
Let $e=gcd(e_1,e_2,\dots,e_M)$.
Choose by Corollary \ref{K_0c} orthogonal non-zero projections
$q_i\in M_e\subseteq B$,
for $i=1,2,\dots,R$, with sum $1$ such that $\kappa_*[p_i]=[q_i]$.
Let $t_i>0$ be the normalized trace of $q_i$. Note that we have a well-defined
map $J_i:\text{Aff}\, T(A_i)\to \text{Aff}\, T(A)$ such that
$J_i(\widehat a)=\widehat{\iota_i^A(a)}$ for every self-adjoint element
$a\in A_i$. Define
$\Theta_i : \text{Aff}\, T(A_i) \to \text{Aff}\, T(q_iBq_i)$ by
$$
\Theta_i(f)(\frac{1}{t_i} \tau\circ\epsilon_i)=
\frac{1}{t_i} \Theta(J_i(f))(\tau),\quad \tau\in T(B),
$$
where $\epsilon_i : q_iBq_i\to B$ denotes the inclusion.

$\Theta_i$ is a linear positive map, and it preserves the order unit since
$$
\Theta_i(1)(\frac{1}{t_i} \tau\circ \epsilon_i)=
\frac{1}{t_i} \Theta(\widehat{p_i})(\tau)=
\frac{1}{t_i} \rho_B\circ \kappa_*[p_i](\tau)=1.
$$
By \cite[Theorem 7.3]{RSUCT} we get that $[\epsilon_i]\in KK(q_iBq_i,B)$ is
a $KK$-equivalence. Note that
$$
[\epsilon_i]^{-1}\cdot\kappa\cdot [\iota_i^A]\in KK(A_i,q_iBq_i)_e
$$
By Corollary \ref{cutdown} we have that
$q_iBq_i\cong A(t_im,t_ie_1,t_ie_2,\dots,t_ie_M)$. Choose by
Lemma \ref{unitgrpAp} a unitary $w_i\in q_iBq_i$ such that
$$
q_B'(w_i+(1-q_i))=q_B'(u_i)\quad \text{in}\ U(B)/\overline{DU(B)}.
$$
Since $t_ie_j\ge K$ for $j=1,2,\dots,M$, we get by Theorem \ref{ex2}
a unital *-homomorphism $\varphi_i : A_i\to q_iBq_i$ such that
\begin{gather*}
[\varphi_i]=[\epsilon_i]^{-1}\cdot\kappa\cdot [\iota_i^A]
\quad\text{in\ }KK(A_i,q_iBq_i), \\
\|\widehat{\varphi_i}(f)-\Theta_i(f)\|<\epsilon,
\quad f\in \widehat{\pi_i^A}(F), \\
\varphi_i(v^{A_i})=w_i
\quad \text{mod}\ \overline{DU(q_iBq_i)}.
\end{gather*}
Now define $\varphi : A\to B$ by
$$
\varphi(a)=\sum_{i=1}^R \epsilon_i\circ\varphi_i\circ\pi_i^A(a).
$$
$\varphi$ is a unital *-homomorphism and
$$
[\varphi]=\sum_{i=1}^R [\epsilon_i]\cdot [\varphi_i]\cdot [\pi_i^A] =
\sum_{i=1}^R  \kappa\cdot [\iota_i^A]\cdot [\pi_i^A]=\kappa
\quad \text{in\ } KK(A,B).
$$
For $f\in \text{Aff}\, T(A)$, $\tau\in T(B)$, we have that
$$
\Theta(f)(\tau)=\sum_{i=1}^R  \Theta(J_i(\widehat{\pi_i^A}(f)))(\tau)=
\sum_{i=1}^R   t_i\Theta_i(\widehat{\pi_i^A}(f))
(\frac{1}{t_i}\tau\circ\epsilon_i),
$$
and
$$
\widehat \varphi(f)(\tau)=f(\tau\circ\varphi)=
f(\sum_{i=1}^R t_i\frac{1}{t_i}\tau\circ\epsilon_i\circ\varphi_i\circ\pi_i^A)=
\sum_{i=1}^R t_i \widehat{\varphi_i}(\widehat{\pi_i^A}(f))
(\frac{1}{t_i}\tau\circ\epsilon_i).
$$
It follows that
$$
\|\widehat\varphi(f)-\Theta(f)\|<\epsilon,\quad f\in F.
$$
Finally, for $i=1,2,\dots,R$,
$$
\varphi(v^A_i)=\epsilon_i\circ\varphi_i(v^{A_i})+(1-q_i)=w_i+(1-q_i)=u_i
$$
modulo $\overline{DU(B)}$.
\end{proof}
\end{thm}

\section{Injective connecting maps}

The purpose of this section is to show that a simple unital infinite
dimensional inductive limit of a sequence of finite direct sums of building
blocks can be realized as an inductive limit of a sequence of finite direct
sums of building blocks with unital and injective connecting maps.

From now on, we will consider inductive limits in the category of order unit
spaces and linear positive order unit preserving maps, as introduced by
Thomsen \cite{TAIT}. It follows from \cite[Lemma 3.3]{TAIT}
that $\text{Aff}\, T(\cdot)$ is a continuous functor from the category of
separable unital $C^*$-algebras and unital *-homomorphisms to the category
of order unit spaces. We will also need Elliott's approximative intertwining
argument, see \cite[Theorem 2.1]{ERR0} or \cite{TAI}.

\begin{lemma}\label{ainj}
Let $A$ be a finite direct sum of building blocks, interval building blocks,
and matrix algebras. Let $\epsilon>0$ and let $F\subseteq A$ be a finite set.
There exists a finite set of positive non-zero elements $H\subseteq A$ such
that if $B$ is a building block or an interval building block, and
$\varphi:A\to B$ is a unital *-homomorphism with $\varphi(h)\neq 0$, $h\in H$,
then there exists a unital injective *-homomorphism $\psi:A\to B$ such
that $\|\varphi(f)-\psi(f)\|<\epsilon$, $f\in F$.

\begin{proof}
By Corollary \ref{cutdown} (and the corresponding result for interval building
blocks) we may assume that $A$ is a building block, an interval building block
or a matrix algebra
rather than a finite direct sum of such algebras. We will carry out the proof
in the case that $A=A(n,d_1,d_2,\dots,d_N)$ is a circle building block.
The proof in the case that $A$ is an interval building block is similar, and
the matrix algebra case is trivial.

Choose $\delta>0$ such that for $x,y\in\T$,
$$
\rho(x,y)<2\delta\ \Rightarrow\ \|f(x)-f(y)\|<\epsilon,\quad f\in F.
$$
Let $\T=\cup_{i=1}^K V_i$ where each $V_i$ is an open arc-segment of length
less than $\delta$. Choose for each $i=1,2,\dots,K$, a non-zero continuous
function $\chi_i:\T\to [0,1]$ with support in $V_i$ such that $\chi_i$ is
zero at every exceptional point of $A$. Set
$$
H=\{\chi_1\otimes 1,\chi_2\otimes 1,\dots,\chi_K\otimes 1\}.
$$
Let $\varphi:A\to B$ be given such that $\varphi(h)\neq 0$, $h\in H$.
By \cite[Chapter 1]{TATD} we may assume that
$$
\varphi(f)(e^{2\pi it})=u(t)\text{diag}\bigl(\Lambda_1^{s_1}(f),
\dots,\Lambda_N^{s_N}(f),f(\lambda_1(t)),\dots,
f(\lambda_L(t))\bigr)u(t)^*,
$$
if $B=A(m,e_1,e_2,\dots,e_M)$ is a circle building block and
$$
\varphi(f)(t)=u(t)\text{diag}\bigl(\Lambda_1^{s_1}(f),\dots,
\Lambda_N^{s_N}(f),f(\lambda_1(t)),\dots,
f(\lambda_L(t))\bigr)u(t)^*,
$$
if $B=I(m,e_1,e_2,\dots,e_M)$ is an interval building block.
Here $u\in C[0,1]\otimes M_m$ is a unitary,
$\lambda_1,\dots,\lambda_L:[0,1]\to\T$ are continuous functions,
and $s_1,s_2,\dots,s_N$ are non-negative integers.
Since $\varphi(h)\neq 0$, $h\in H$, it follows that the set
$\bigcup_{k=1}^L \lambda_k([0,1])$ intersects non-trivially with every $V_i$.

If $B$ is an interval building block, let $t_1,t_2,\dots,t_M\in [0,1]$ be the
exceptional points of $B$. If $B$ is a circle building block, let
$t_1,t_2,\dots,t_M\in [0,1]$ be numbers such that $e^{2\pi it_j}$,
$j=1,2,\dots,M$, are the exceptional points of $B$.

For each $k=1,2,\dots,L$, choose a continuous function
$\mu_k:[0,1]\to\T$ such that $\rho(\mu_k(t),\lambda_k(t))<2\delta$,
$t\in [0,1]$, such that $\mu_k(t)=\lambda_k(t)$ for
$t\in \{t_1,t_2,\dots,t_M,0,1\}$, and such that
$\bigcup_{i=1}^k \mu_k([0,1])=\bigcup_{i=1}^k V_i=\T$.
Define $\psi:A\to B$ by
$$
\psi(f)(e^{2\pi it})=u(t)\text{diag}\bigl(\Lambda_1^{s_1}(f),\dots,
\Lambda_N^{s_N}(f),f(\mu_1(t)),\dots,f(\mu_L(t))\bigr)u(t)^*
$$
if $B$ is a circle building block, and
$$
\psi(f)(t)=u(t)\text{diag}\bigl(\Lambda_1^{s_1}(f),\dots,
\Lambda_N^{s_N}(f),f(\mu_1(t)),\dots,f(\mu_L(t))\bigr)u(t)^*.
$$
if $B$ is an interval building block. Note that $\psi$ is injective and
$\|\varphi(f)-\psi(f)\|<\epsilon$, $f\in F$.
\end{proof}
\end{lemma}

\begin{lemma}\label{unital}
Let $A$ be a unital $C^*$-algebra that is the inductive limit of a sequence
$$
  \begin{CD}
      A_1  @> \alpha_1 >> A_2 @> \alpha_2 >>
A_3 @>\alpha_3  >> \dots \\   
  \end{CD}
$$
of finite direct sums of building blocks.
Then $A$ is the inductive limit of a similar sequence, with unital
connecting maps.

\begin{proof}
Note that we may assume that $\alpha_{n,\infty}(p)\neq 0$ for every
positive integer $n$ and every minimal non-zero central projection $p\in A_n$.
By Lemma \ref{ideal} it follows that $\alpha_{n,\infty}(q)\neq 0$ for
every projection $q\in A_n$.
Let $1_n\in A_n$ denote the unit.
Since $\{\alpha_{n,\infty}(1_n)\}_{n=1}^{\infty}$ is an approximate unit for
$A$ there exists a positive integer $N$ such that $\alpha_{k,\infty}(1_k)=1$
for all $k\ge N$. Hence $\alpha_k(1_k)=1_{k+1}$, $k\ge N$.
\end{proof}
\end{lemma}

\begin{lemma}\label{nice}
Let $X\subseteq \T$ be a closed set and let $G\subseteq X$ be a finite
subset. Let $\delta>0$ be given. There exist a closed subset $R\subseteq X$
with finitely many connected components such that $G\subseteq R$,
together with a continuous surjective map $g: X \to R$ such that
$g(z)=z$, $z\in G$, and $\rho(g(z),z)\le \delta$, $z\in X$.

\begin{proof}
Let $G=\{e^{2\pi it_j}:j=1,2,\dots,N\}$ where
$0\le t_1<t_2<\dots<t_N<1$. Set $t_{N+1}=t_1+1$. Set
$I_j=\{e^{2\pi it}:t\in [t_j,t_{j+1}]\}$. We may assume that
$t_{j+1}-t_j<\delta$ unless the interior of $I_j$ intersects non-trivially
with $X$. On each $I_j$ let either $g$
be the identity map (if $I_j\subseteq X$) or a continuous map onto
$\{e^{2\pi it_j},e^{2\pi it_{j+1}}\}$ that is constant on the set of boundary
points of $I_j$. Set $R=g(X)$.
\end{proof}
\end{lemma}

\begin{lemma}\label{nicebb}
Let $A$ be a quotient of a finite direct sum of building blocks. Let
$F\subseteq A$ be a finite set and let $\epsilon>0$. There exists a finite
direct sum of building blocks, interval building blocks and matrix algebras
$B$, and unital *-homomorphisms $\varphi:A\to B$ and $\psi:B\to A$ such
that $\psi$ is injective and $\|\psi\circ\varphi(f)-f\|<\epsilon$, $f\in F$.

\begin{proof}
We may assume that $A$ is a quotient of a building block rather than of a
finite direct sum of building blocks. Hence by Lemma \ref{ideal}
$$
A=\{f\in C(X)\otimes M_n:f(x_i)\in M_{d_i},\ i=1,2,\dots,N\}
$$
where $X\subseteq\T$ is a closed subset and $x_1,x_2,\dots,x_N\in X$.
Choose $\delta>0$ such that
$$
y,z\in X,\ \rho(y,z)\le\delta\ \Rightarrow\
\|f(y)-f(z)\|<\epsilon,\quad f\in F.
$$
Choose by Lemma \ref{nice} a closed subset $R\subseteq X$ with finitely many
connected components such that $x_1,x_2,\dots,x_N\in R$, and a continuous
surjective map $g:X\to R$ such that $g(x_i)=x_i$, $i=1,2,\dots,N$,
and such that $\rho(g(z),z)\le\delta$, $z\in X$. Let
$$
B=\{f\in C(R)\otimes M_n:f(x_i)\in M_{d_i},\ i=1,2,\dots,N\}.
$$
Define $\psi:B\to A$ by $\psi(f)=f\circ g$ and let $\varphi:A\to B$ be
restriction. Then $\|\psi\circ\varphi(f)-f\|<\epsilon$, $f\in F$.
\end{proof}
\end{lemma}

\begin{prop}\label{injcim}
Let $A$ be a unital simple inductive limit of a sequence of finite direct
sums of building blocks. Then $A$ is the inductive limit of a sequence
of finite direct sums of building blocks, interval building blocks and
matrix algebras, with unital and injective connecting maps.

\begin{proof}
By Lemma \ref{unital} we have that $A$ is the inductive limit of a sequence
$$
  \begin{CD}
      A_1  @> \alpha_1 >> A_2 @> \alpha_2 >>
A_3 @>\alpha_3  >> \dots \\   
  \end{CD}
$$
where each $\alpha_n$ is unital and injective and each $A_n$ is a quotient
of a finite direct sum of building blocks.
We will construct a strictly increasing sequence of positive integers
$\{n_k\}$, a sequence
$$
  \begin{CD}
      B_1  @> \beta_1 >> B_2 @> \beta_2 >>
B_3 @>\beta_3  >> \dots \\   
  \end{CD}
$$
of finite direct sums of building blocks, interval building
blocks and matrix algebras with unital connecting maps,
unital *-homomorphisms $\mu_k:A_{n_k}\to B_{k+1}$, and unital
injective *-homomorphisms $\psi_k:B_k\to A_{n_k}$ such that the diagram
$$
\xymatrix{
A_{n_1} \ar[r]^{\alpha_{n_1,n_2}}\ar[dr]^{\mu_1}  &
A_{n_2} \ar[r]^{\alpha_{n_2,n_3}}\ar[dr]^{\mu_2}  &
A_{n_3} \ar[r]^{\alpha_{n_3,n_4}}\ar[dr]^{\mu_3}  & \dots \\
B_1 \ar[u]_{\psi_1}\ar[r]_{\beta_1} &
B_2 \ar[u]_{\psi_2}\ar[r]_{\beta_2} &
B_3 \ar[u]_{\psi_3}\ar[r]_{\beta_3} & \dots}
$$
becomes an approximate intertwining. Furthermore $\beta_k$ should be injective
unless $B_{k+1}$ is finite dimensional. This is sufficient since the
proposition is trivial if $A$ is an AF-algebra.

It is easy to construct $B_1$, $n_1$ and $\psi_1$.
Assume that $B_k$, $n_k$ and $\psi_k$ have been constructed. Let $\epsilon>0$
and finite sets $F\subseteq A_{n_k}$ and $G\subseteq B_k$ be given. Choose
$H\subseteq B_k$ by Lemma \ref{ainj} with respect to $\epsilon>0$ and $G$.
Since $A$ is simple we may choose $n_{k+1}$ such that
$\widehat{\alpha_{n_k,n_{k+1}}}(\widehat{\psi_k(h)})>0$ for $h\in H$.
Choose $B_{k+1}$, $\varphi_{k+1}$ and $\psi_{k+1}$ by Lemma \ref{nicebb}
with respect to $\epsilon>0$ and $\alpha_{n_k,n_{k+1}}(F)$. Set
$\mu_k=\varphi_{k+1}\circ\alpha_{n_k,n_{k+1}}$. Then
$$
\|\psi_{k+1}\circ\mu_k(x)-\alpha_{n_k,n_{k+1}}(x)\|<\epsilon,\quad x\in F.
$$
Since $\widehat{\mu_k}\circ\widehat{\psi_k}(\widehat{h})>0$, $h\in H$, there
exists by Lemma \ref{ainj} a unital *-homomorphism $\beta_k:B_k\to B_{k+1}$
such that
$$
\|\mu_k\circ\psi_k(x)-\beta_k(x)\|<\epsilon,\quad x\in G,
$$
and such that $\beta_k$ is injective if $B_{k+1}$ is infinite dimensional.
\end{proof}
\end{prop}

\begin{lemma}\label{infty}
Let $A$ be a simple infinite dimensional inductive limit of a sequence
$$
  \begin{CD}
      A_1  @> \alpha_1 >> A_2 @> \alpha_2 >>
A_3 @>\alpha_3  >> \dots \\   
  \end{CD}
$$
of finite direct sums of building blocks, interval building blocks and
matrix algebras, with unital and injective connecting maps.
Then $s(A_m)\to\infty$.

\begin{proof}
The lemma is well-known if $A$ is an AF-algebra. We may therefore assume
that $A_k$ is infinite dimensional for some $k$.
Let $L$ be a positive integer. Let $b_1,b_2,\dots,b_L\in A_k$ be positive
non-zero mutually orthogonal elements. Since $A$ is simple and the connecting
maps are injective, there exists an integer $N\ge k$ such that
$$
\widehat{\alpha_{k,N}}(\widehat{b_j})>0,\quad j=1,2,\dots,L.
$$
Hence if $m\ge N$ and $\mu:A_m\to M_n$ is a unital *-homomorphism, we see
that the elements $\mu\circ\alpha_{1,m}(b_j)$, $j=1,2,\dots,L$, are
non-zero and mutually orthogonal. Thus $n\ge L$.
\end{proof}
\end{lemma}

\begin{prop}\label{injci}
Let $A$ be a simple unital infinite dimensional inductive limit of a sequence
of finite direct sums of building blocks. Then $A$ is the inductive limit of a
sequence of finite direct sums of building blocks and interval building blocks
with unital and injective connecting maps.

\begin{proof}
By Proposition \ref{injcim} we have that $A$ is the inductive limit
of a sequence
$$
  \begin{CD}
      A_1  @> \alpha_1 >> A_2 @> \alpha_2 >>
A_3 @>\alpha_3  >> \dots \\   
  \end{CD}
$$
where each $\alpha_k$ is unital and injective and each $A_k$ is of the form
$C_k\oplus F_k$ for a finite (possibly trivial) direct sum of building blocks
$C_k$ and a finite dimensional $C^*$-algebra $F_k$. Set
$B_k=C_k\oplus (C(\T)\otimes F_k)$ and let $\psi_k:A_k\to B_k$ be the canonical
*-homomorphism. It suffices to construct a strictly increasing sequence
of positive integers $\{n_k\}$, unital *-homomorphisms
$\mu_k:B_{n_k}\to A_{n_{k+1}}$, and unital injective *-homomorphisms
$\beta_k: B_{n_k}\to B_{n_{k+1}}$ such that the diagram
$$
\xymatrix{
A_{n_1} \ar[r]^{\alpha_{n_1,n_2}}\ar[d]^{\psi_{n_1}}  &
A_{n_2} \ar[r]^{\alpha_{n_2,n_3}}\ar[d]^{\psi_{n_2}}  &
A_{n_3} \ar[r]^{\alpha_{n_3,n_4}}\ar[d]^{\psi_{n_3}}  & \dots \\
B_{n_1} \ar[ur]_{\mu_1}\ar[r]_{\beta_1} &
B_{n_2} \ar[ur]_{\mu_2}\ar[r]_{\beta_2} &
B_{n_3} \ar[ur]_{\mu_3}\ar[r]_{\beta_3} & \dots}
$$
becomes an approximate intertwining. This is done by induction. Set $n_1=1$.

Assume that $n_k$ has been constructed. Let $\epsilon>0$ and a finite set
$G\subseteq B_{n_k}$ be given. It suffices to construct $n_{k+1}>n_k$, a
unital *-homomorphism $\mu_k:B_{n_k}\to A_{n_{k+1}}$ such that
$\mu_k\circ\psi_{n_k}$ and $\alpha_{n_k,n_{k+1}}$ are approximately unitarily
equivalent, and a unital injective *-homomorphism
$\beta_k:B_{n_k}\to B_{n_{k+1}}$ such that
\begin{equation}\label{injcieq}
\|\beta_k(x)-\psi_{n_{k+1}}\circ\mu_k(x)\|<\epsilon,\quad x\in G.
\end{equation}
Let $F_{n_k}=M_{m_1}\oplus M_{m_2}\oplus\dots\oplus M_{m_N}$ and let
$p_1,p_2,\dots,p_N$ be the minimal non-zero central projections in
$F_{n_k}\subseteq A_{n_k}$. Let $\pi_i:B_{n_k}\to C(\T)\otimes M_{m_i}$
be the projection, $i=1,2,\dots,N$. Choose by Lemma \ref{ainj} a finite set
$H_i\subseteq C(\T)\otimes M_{m_i}$ of positive non-zero elements with
respect to $\epsilon$ and $\pi_i(G)$. Let $h_i$ be the cardinality of $H_i$.

Since $A$ is simple there exists a $\delta>0$ such that
$$
\widehat{\alpha_{n_k,\infty}}(\widehat{p_i})>\delta,\quad i=1,2,\dots,N.
$$
By Lemma \ref{infty} there exists an integer $n_{k+1}>n_k$ such that
\begin{gather*}
\widehat{\alpha_{n_k,n_{k+1}}}(\widehat{p_i})>\delta,\quad i=1,2,\dots,N, \\
s(A_{n_{k+1}})>\delta^{-1}\max_i (h_im_i),\quad i=1,2,\dots,N.
\end{gather*}
Let $q_i=\alpha_{n_k,n_{k+1}}(p_i)$ and note that
$$
s(q_iA_{n_{k+1}}q_i)>\delta s(A_{n_{k+1}})> h_i m_i,\quad i=1,2,\dots,N.
$$
Hence there exists a unital *-homomorphism
$\lambda_i:C(\T)\otimes M_{m_i}\to q_iA_{n_{k+1}}q_i$ such that
$\widehat{\lambda_i}(\widehat{h})>0$, $h\in H_i$.
Let $\mu_k:B_{n_k}\to A_{n_{k+1}}$ be the *-homomorphism that agrees with
$\alpha_{n_k,n_{k+1}}$ on $C_{n_k}$ and with $\lambda_i$ on
$C(\T)\otimes M_{m_i}$.
The *-homomorphism $x\mapsto \lambda_i(1\otimes x)$ from $M_{m_i}$ to
$q_iA_{n_{k+1}}q_i$ is  by \cite[Chapter 1]{TATD} approximately unitarily
equivalent to the *-homomorphism induced by
$\alpha_{n_k,n_{k+1}}$. Hence $\mu_k\circ\psi_{n_k}$ and $\alpha_{n_k,n_{k+1}}$
are approximately unitarily equivalent.

Let $e_i=\psi_{n_{k+1}}(q_i)$, $i=1,2,\dots,N$ and let
$\xi_i:q_iA_{n_{k+1}}q_i\to e_iB_{n_{k+1}}e_i$ be the unital *-homomorphism
induced by $\psi_{n_{k+1}}$. Since
$\widehat{\xi_i}\circ\widehat{\lambda_i}(\widehat{h})>0$ there exists
by Lemma \ref{ainj} a unital injective *-homomorphism
$\varphi_i:C(\T)\otimes M_{m_i}\to e_iB_{n_{k+1}}e_i$ such that
$$
\|\varphi_i(x)-\xi_i\circ\lambda_i(x)\|<\epsilon,
\quad x\in\pi_i(G),\quad i=1,2,\dots,N.
$$
Let $\beta_k$ be the *-homomorphism that agrees with $\psi_{n_{k+1}}\circ\mu_k$
on $C_{n_k}$ and with $\varphi_i$ on $C(\T)\otimes M_{m_i}$. Note that
$\beta_k$ is unital and injective and that (\ref{injcieq}) holds.
\end{proof}
\end{prop}

It remains to replace interval building blocks with building blocks.
This turns out to be much more complicated than in \cite[Lemma 1.5]{NT} or
\cite[Lemma 4.7]{TATD}, since our building blocks may be unital projectionless.
We will use the following lemma, which resembles a uniqueness result for
interval building blocks. The proof is inspired by Elliott's proof of
the uniqueness lemma for interval algebras \cite{EAIS}.

\begin{lemma}\label{unint}
Let $A=I(n,d_1,d_2,\dots,d_N)$ be an interval building block. Let
$F\subseteq A$ be a finite set and let $\epsilon>0$ be given. There is a
finite set $H\subseteq A$ of positive elements of norm $1$ such that if
$B=I(m,e_1,e_2,\dots,e_M)$ is an interval building block with exceptional
points $y_1,y_2,\dots,y_M$, if $\varphi,\psi:A\to B$ are unital
*-homomorphisms and if $\delta>0$, such that
\begin{enumerate}
\item[(i)]
$\|\widehat{\varphi}(\widehat{h})-\widehat{\psi}(\widehat{h})\|<\delta,
\quad h\in H$,
\item[(ii)]
$\widehat{\varphi}(\widehat{h})>\delta,\quad h\in H$,
\item[(iii)]
$\widehat{\psi}(\widehat{h})>\delta,\quad h\in H$,
\item[(iv)]
$f\mapsto\varphi(f)(y_j)$ and $f\mapsto\psi(f)(y_j)$ are equivalent
representations of $A$ on $M_m$, $j=1,2,\dots,M$;
\end{enumerate}
then there is a unitary $W\in B$ such that
$$
\|\varphi(f)-W\psi(f)W^*\|<\epsilon,\qquad f\in F.
$$

\begin{proof}
We may assume that $\|f\|\le 1$ for $f\in F$. Let $x_1,x_2,\dots,x_N$ be the
exceptional points of $A$. Choose a positive integer $q$ such that
\begin{gather*}
\frac 2q < \min\{|x_i-x_j|:i\neq j\}, \\
|x-y|\le \frac 3q\ \Rightarrow\ \|f(x)-f(y)\|<\frac{\epsilon}3,\quad f\in F.
\end{gather*}
For $r=1,2,\dots,q$, define a continuous function $h_r:[0,1]\to [0,1]$ by
$$
h_r(t)=
\begin{cases}
0 & 0\le t\le \frac{r-1}{q}, \\
qt-(r-1)& \frac{r-1}q\le t \le \frac rq, \\
1 & \frac rq \le t \le 1.
\end{cases}
$$
Set
$$
H=\{h_1\otimes 1,h_2\otimes 1,\dots,h_q\otimes 1\}\cup
\{ (h_1-h_2)\otimes 1,
\dots, (h_{q-1}-h_q)\otimes 1\}.
$$

Let $\varphi,\psi : A\to B$ be unital *-homomorphisms that satisfy (i)-(iv).
By \cite[Chapter 1]{TATD} we see that $\varphi$ and $\psi$ are
approximately unitarily equivalent to *-homomorphisms of the form
\begin{gather*}
\varphi'(f)(t)=u(t)\text{diag}\bigl(
\Lambda_1^{r_1}(f),\dots,\Lambda_N^{r_N}(f),
f(\lambda_1(t)),\dots,f(\lambda_L(t))
\bigr) u(t)^* \\
\psi'(f)(t)=v(t)\text{diag}\bigl(
\Lambda_1^{s_1}(f),\dots,\Lambda_N^{s_N}(f),
f(\mu_1(t)),\dots,f(\mu_K(t))\bigr) v(t)^*
\end{gather*}
for continuous functions $\lambda_1\le\lambda_2\le\dots\le\lambda_L,
\mu_1\le\mu_2\le\dots\le\mu_K:[0,1]\to [0,1]$,
integers $r_i$ and $s_i$ with $0\le r_i<\frac n{d_i}$,
$0\le s_i<\frac n{d_i}$ for $i=1,2,\dots,N$,
and unitaries $u,v\in C[0,1]\otimes M_m$.
By (iv) we have that $f\mapsto\varphi'(f)(y_j)$ and 
$f\mapsto\psi'(f)(y_j)$ are equivalent representations of $A$ on $M_m$,
$j=1,2,\dots,M$. It follows that $r_i=s_i$, $i=1,2,\dots,N$, that $K=L$,
and that
$$
(\lambda_1(y_j),\lambda_2(y_j),\dots,\lambda_L(y_j))=
(\mu_1(y_j),\mu_2(y_j),\dots,\mu_L(y_j))
$$
as unordered $L$-tuples, $j=1,2,\dots,M$. Hence
\begin{equation}\label{mueqlambda}
\lambda_k(y_j)=\mu_k(y_j),\quad j=1,2,\dots,M,\ k=1,2,\dots,L.
\end{equation}

For every $t\in [0,1],\ r=2,3,\dots,q$, we have that
\begin{align*}
\#\{ k: \lambda_k(t)\ge \frac rq \} n\ + 
\sum_{ i:\, x_i\ge \frac rq} r_id_i
\ \le\ &Tr(\varphi'(h_r\otimes 1)(t))\\
< \ &m\delta + Tr(\psi'(h_r\otimes 1)(t))\\
<\ &Tr(\psi'(h_{r-1}\otimes 1)(t)) \\
\le\ &\# \{ k:\mu_k(t)\ge \frac{r-2}q\} n\ + 
\sum_{i:\, x_i\ge\frac {r-2}q} r_id_i.
\end{align*}
As $[\frac{r-2}q,\frac rq]$ at most contains one of the exceptional points
of $A$, we see that
$$
\#\{k:\lambda_k(t)\ge\frac rq\} n\ 
<\ \#\{ k:\mu_k(t)\ge \frac{r-2}q\}n + n.
$$
Thus
$$
\#\{k:\lambda_k(t)\ge\frac rq\}\ \le\ \#\{ k:\mu_k(t)\ge \frac{r-2}q\}.
$$
It follows that $\lambda_k(t)\le \mu_k(t)+\frac 3q$. Symmetry allows us to
conclude that for all $t\in [0,1]$,
$$
|\lambda_k(t)-\mu_k(t)|\le \frac 3q, \quad k=1,2,\dots,L.
$$
By (\ref{mueqlambda}) we can define a *-homomorphism $\beta:A\to B$ by
$$
\beta(f)(t)=v(t)\text{diag}\bigl(
\Lambda_1^{r_1}(f),\dots,\Lambda_N^{r_N}(f),
f(\lambda_1(t)),\dots,f(\lambda_L(t))\bigr)v(t)^*
$$
Note that
$$
\|\psi'(f)-\beta(f)\|<\frac{\epsilon}3,\quad f\in F.
$$
Since $Tr(\beta(f)(t))=Tr(\varphi(f)(t))$, $f\in A$, $t\in [0,1]$, it follows
that $\beta$ and $\varphi$ are approximately unitarily equivalent by
\cite[Corollary 1.5]{TATD}. Hence there exists a unitary $U\in B$ such that
$$
\|\varphi(f)-U\beta(f)U^*\|<\frac{\epsilon}3,\quad f\in F.
$$
Choose a unitary $V\in B$ such that
$$
\|\psi'(f)-V\psi(f)V^*\|<\frac{\epsilon}3,\quad f\in F.
$$
Set $W=UV$. Then for $f\in F$,
\begin{align*}
&\|\varphi(f)-W\psi(f)W^*\| \\
\le\ &\|\varphi(f)-U\beta(f)U^*\|+
\|U\beta(f)U^*-U\psi'(f)U^*\|+\|U\psi'(f)U^*-UV\psi(f)V^*U^*\| \\
<\ &\frac{\epsilon}3+\frac{\epsilon}3+\frac{\epsilon}3=\epsilon.
\end{align*}
\end{proof}
\end{lemma}

Define a continuous function $\kappa:\T\to [0,1]$ by
$$
\kappa(e^{2\pi it})=
\begin{cases}
2t & t\in [0,\frac 12], \\
2-2t & t\in [\frac 12,1].
\end{cases}
$$
Define continuous functions $\iota_1,\iota_2:[0,1]\to\T$ by
$\iota_1(t)=e^{\pi it}$, $\iota_2(t)=e^{-\pi it}$.
Note that $\kappa\circ\iota_1=\kappa\circ\iota_2=id_{[0,1]}$.

Let $A=I(n,d_1,d_2,\dots,d_N)$ be an interval building block with
exceptional points $t_1,t_2,\dots,t_N$. Define a circle building block by
$$
A^{\T}=\{f\in C(\T)\otimes M_n:
f(\iota_1(t_i)),f(\iota_2(t_i))\in M_{d_i},\ i=1,2,\dots,N\}.
$$
Define unital *-homomorphisms $\xi_A:A\to A^{\T}$ by
$\xi_A(f)=f\circ\kappa$, $f\in A$,
and $j_A^1,j_A^2:A^{\T}\to A$ by
$j_A^1(g)=g\circ\iota_1$, $j_A^2(g)=g\circ\iota_2$, $g\in A^{\T}$. Then
$j_A^1\circ\xi_A=j_A^2\circ\xi_A=id_A$.

Let $A$ be a finite direct sum of building blocks and interval building
blocks. It follows from the above that there exists a finite direct sum of
building blocks $A^{\T}$ together with unital *-homomorphisms
$\xi_A:A\to A^{\T}$ and $j_A^1,j_A^2:A^{\T}\to A$ such that
$j_A^1\circ\xi_A=j_A^2\circ\xi_A=id_A$ and
\begin{equation}\label{injeq}
j_A^1(f)=j_A^2(f)=0\ \Rightarrow\ f=0,\quad f\in A^{\T}.
\end{equation}

\begin{thm}\label{inj}
Let $A$ be a simple unital infinite dimensional inductive limit of a
sequence of finite direct sums of circle building blocks. Then $A$ is the
inductive limit of a sequence of finite direct sums of circle building blocks
with unital and injective connecting maps.

\begin{proof}
By Proposition \ref{injci} we see that $A$ is the inductive limit of a
sequence
$$ 
  \begin{CD}
      A_1  @> \alpha_1 >> A_2 @> \alpha_2 >>
A_3 @>\alpha_3  >> \dots \\   
  \end{CD}
$$
where each $A_n$ is a finite direct sum of circle and interval building
blocks and each $\alpha_n$ is a unital and injective *-homomorphism.

By passing to a subsequence, if necessary, we may assume that either, every
$A_n$ is a circle or an interval building block or,
every $A_n$ is a finite direct sum of at least two circle or interval
building blocks.

Let us first assume that the latter is the case.

Let $A_n=A^n_1\oplus A^n_2\oplus\dots\oplus A^n_{N_n}$ where each $A_i^n$ is
a circle or an interval building block. For each $n$ let
$\pi_i^n:A_n\to A^n_i$ denote the coordinate projections, $i=1,2,\dots,N_n$.
First we claim that we may assume that all the maps
$\pi_i^{n+1}\circ\alpha_n$ are injective.

By Elliott's approximate intertwining argument it suffices to show that
given a finite set $G\subseteq A_n$ and $\epsilon>0$ there exists an
integer $m>n$ and a unital *-homomorphism $\psi:A_n\to A_m$ such that
$\|\alpha_{n,m}(g)-\psi(g)\|<\epsilon$, $g\in G$, and such that
$\pi_i^m\circ\psi$ is injective, $i=1,2,\dots,N_m$.
Choose by Lemma \ref{ainj} a finite set $H\subseteq A_n$ of positive
non-zero elements with respect to
$G$ and $\epsilon$. As $A$ is simple and the connecting maps are injective,
we have that $\widehat{\alpha_{n,\infty}}(\widehat{h})>0$, $h\in H$. Thus there
exists an integer $m>n$ such that
$\widehat{\alpha_{n,m}}(\widehat{h})>0$, $h\in H$. Hence
$\pi_i^m\circ\alpha_{n,m}(h)\neq 0$, $i=1,2,\dots,N_m$, and the claim
follows by $N_m$ applications of Lemma \ref{ainj}.

Define a unital *-homomorphism $\psi_n:A_n^{\T}\to A_{n+1}$ by
$$
\psi_n(x)=(\pi_1^{n+1}\circ\alpha_n\circ j_{A_n}^1(x),
\pi_2^{n+1}\circ\alpha_n\circ j_{A_n}^2(x),\dots,
\pi_{N_{n+1}}^{n+1}\circ\alpha_n\circ j_{A_n}^2(x)).
$$
Since the maps $\pi_i^{n+1}\circ\alpha_n$ are injective, $i=1,2,\dots,N_{n+1}$,
and as $N_{n+1}\ge 2$, it follows from (\ref{injeq}) that $\psi_n$ is
injective. The theorem therefore follows in this case from the commutativity of
the diagram
$$
\xymatrix{
A_1 \ar[r]^{\alpha_1}\ar[d]^{\xi_{A_1}}  & 
A_2 \ar[r]^{\alpha_2}\ar[d]^{\xi_{A_2}}  &
A_3 \ar[r]^{\alpha_3}\ar[d]^{\xi_{A_3}}  & \dots \\
A_1^{\T} \ar[ur]_{\psi_1}\ar[r]_{\xi_{A_2}\circ\psi_1} &
A_2^{\T} \ar[ur]_{\psi_2}\ar[r]_{\xi_{A_3}\circ\psi_2} &
A_3^{\T} \ar[ur]_{\psi_3}\ar[r]_{\xi_{A_4}\circ\psi_3}  & \dots}
$$

It remains to prove the theorem in the first case. By passing to a subsequence
we may assume that each $A_n$ is an interval building block.
Let $\epsilon>0$, let $k$ be a positive integer, and let $F\subseteq A_k$
be finite. Again by Elliott's approximative intertwining argument, it suffices
to show that there exists an integer $l>k$ and a unital
and injective *-homomorphism $\psi:A_k^{\T}\to A_l$ such that
$$
\|\alpha_{k,l}(x)-\psi\circ\xi_{A_k}(x)\|<\epsilon,\quad x\in F.
$$

Choose by Lemma \ref{unint} a finite set $H\subseteq A_k$ of
positive elements of norm $1$ with respect to $F$ and $\epsilon$.
Since $A$ is simple and the connecting maps are injective there exists a
$\delta>0$ such that $\widehat{\alpha_{k,\infty}}(\widehat{h})>2\delta$,
$h\in H$. Let $A_k=I(n,d_1,d_2,\dots,d_N)$. By Lemma \ref{infty} there exists
an integer $l>k$ such that $s(A_l)>\frac{2n}{\delta}$ and such that
$$
\widehat{\alpha_{k,l}}(\widehat{h})>2\delta,\quad h\in H.
$$
Let $A_l=I(m,e_1,e_2,\dots,e_M)$. By \cite[Chapter 1]{TATD}
$\alpha_{k,l}\circ j_{A_k}^1:A_k^{\T}\to A_l$ is approximately unitarily
equivalent to a *-homomorphism $\beta: A_k^{\T}\to A_l$ of the form
$$
\beta(f)(t)=u(t)\text{diag}\bigl(\Lambda_1^{r_1}(f),\dots,\Lambda_N^{r_N}(f),
f(\mu_1(t)),\dots,f(\mu_L(t))\bigr)u(t)^*,\quad t\in [0,1],
$$
where $u\in C[0,1]\otimes M_m$ is a unitary and
$\mu_1,\mu_2,\dots,\mu_L:[0,1]\to\T$ are continuous functions.
Choose a continuous function $\mu_1':[0,1]\to\T$ such that $\mu_1'=\mu_1$
at the exceptional points of $A_l$ and such that $\mu_1'$ is surjective.
Define $\varphi:A_k^{\T}\to A_l$ by
$$
\varphi(f)(t)=u(t)\text{diag}\bigl(\Lambda_1^{r_1}(f),\dots,\Lambda_N^{r_N}(f),
f(\mu_1'(t)),f(\mu_2(t)),\dots,f(\mu_L(t))\bigr)u(t)^*.
$$
Note that $\varphi$ is injective, and that for $h\in H$,
$$
\|\widehat{\varphi}\circ\widehat{\xi_{A_k}}(\widehat{h})-
\widehat{\alpha_{k,l}}(\widehat{h})\|=
\|\widehat{\varphi}(\widehat{\xi_{A_k}}(\widehat{h}))-
\widehat{\alpha_{k,l}}\circ\widehat{j_{A_k}^1}(\widehat{\xi_{A_k}}
(\widehat{h}))\|\le
\|\widehat{\varphi}-\widehat{\beta}\|\le\frac{2n}m<\delta.
$$
Finally, as $\Lambda_j\circ\varphi=\Lambda_j\circ\beta$, $j=1,2,\dots,M$,
we see by Lemma \ref{unint} that there exists a unitary $W\in A_l$ such that
$$
\|W\varphi\circ\xi_{A_k}(f)W^*-\alpha_{k,l}(f)\|<\epsilon,
\quad f\in F.
$$
Set $\psi(x)=W\varphi(x)W^*$, $x\in A_k^{\T}$.
\end{proof}
\end{thm}

\section{Construction of a certain map}

In \cite{RKL} R\o rdam defined the bifunctor $KL$ to be a certain quotient of
$KK$. Some of our main results are more elegantly formulated in terms of $KL$
than $KK$, and we will therefore from now on use $KL$ instead of $KK$.
Recall from \cite{RKL} that the Kasparov product yields a product
$KL(B,C)\times KL(A,B)\to KL(A,C)$. Furthermore, if $K_*(A)$ is finitely
generated then $KL(A,\cdot)\cong KK(A,\cdot)$, and this functor is continuous
by \cite[Theorem 1.14]{RSUCT} and \cite[Theorem 7.13]{RSUCT}. Finally,
approximately unitarily equivalent *-homomorphisms
define the same element of $KL$ \cite[Proposition 5.4]{RKL}.
It should be noted that $KL$ is related to homomorphisms of $K$-theory with
coefficients, \cite{DLUMCT}.

Let $A$ and $B$ be unital $C^*$-algebras. Let $KL(A,B)_e$ be the set of
elements $\kappa\in KK(A,B)$ for which $\kappa_*:K_0(A)\to K_0(B)$
preserves the order unit. Let $KL(A,B)_T$ be those elements
$\kappa\in KL(A,B)_e$ for which there exists an affine continuous map
$\varphi_T:T(B)\to T(A)$ such that
$$
r_B(\omega)(\kappa_*(x))=r_A(\varphi_T(\omega))(x),
\quad x\in K_0(A),\ \omega\in T(B).
$$

\begin{lemma}\label{KKM}
Let $C$ be a finite direct sum of building blocks, let $\epsilon>0$, and let
$F\subseteq\text{Aff}\, T(C)$ be a finite set. Let $B$ be the inductive
limit of a sequence of finite direct sums of building blocks
$$
  \begin{CD}
      B_1  @> \beta_1 >> B_2 @> \beta_2 >>
B_3 @>\beta_3  >> \dots \\   
  \end{CD}
$$
with unital connecting maps. Let $J:\text{Aff}\, T(C)\to\text{Aff}\, T(B)$
be a linear positive order unit preserving map and let $\kappa\in KL(C,B)_e$.
There exists a positive integer $n$, a linear positive order unit preserving
map $M:\text{Aff}\, T(C)\to\text{Aff}\, T(B_n)$, and an element
$\omega\in KK(C,B_n)_e$ such that
\begin{gather*}
\|J(f)-\widehat{\beta_{n,\infty}}\circ M(f)\|<\epsilon,\quad f\in F, \\
\kappa=[\beta_{n,\infty}]\cdot \omega \quad\text{in\ }KL(C,B).
\end{gather*}

\begin{proof}
We may assume that $\|f\|\le 1$, $f\in F$.
Decompose $C=C_1\oplus C_2\oplus\dots\oplus C_N$ as a finite direct
sum of building blocks and let $\pi_i:C\to C_i$ denote the projection,
$i=1,2,\dots,N$.

For every $i=1,2,\dots,N$,
identify $\text{Aff}\, T(C_i)$ and $C_{\R}(\T)$. Choose open sets
$V_1,V_2,\dots,V_{k_i}\subseteq \T$ such that $\cup_{j=1}^{k_i} V_j=\T$
and such that
$$
x,y\in V_j\ \Rightarrow\  |f(x)-f(y)|<\frac{\epsilon}2,
\quad f\in\widehat{\pi_i}(F).
$$
Let $\{h_j:j=1,2,\dots,k_i\}$ be a continuous partition of
unity in $C_{\R}(\T)$ subordinate to the cover $\{V_j:j=1,2,\dots,k_i\}$
and let $x_j\in V_j$ be an arbitrary point, $j=1,2,\dots,k_i$. Define linear
positive order unit preserving maps
$T_i:\text{Aff}\, T(C_i)\to \R^{k_i}$ and
$S_i:\R^{k_i}\to \text{Aff}\, T(C_i)$ by
\begin{gather*}
T_i(f)=(f(x_1),f(x_2),\dots,f(x_{k_i})), \\
S_i(t_1,t_2,\dots,t_{k_i})=\sum_{j=1}^{k_i} t_jh_j.
\end{gather*}
Note that
$$
\|S_i\circ T_i(f)-f\|<\frac{\epsilon}2,\quad f\in\widehat{\pi_i}(F).
$$
Hence there exist linear positive order unit preserving maps
\begin{gather*}
T: \text{Aff}\, T(C)\to \R^k, \\
S: \R^k\to \text{Aff}\, T(C),
\end{gather*}
where $k=\sum_{i=1}^N k_i$, such that
$$
\|S\circ T(f)-f\|<\frac{\epsilon}2,\quad f\in F.
$$
Let $\{e_j:j=1,2,\dots,k\}$ be the standard basis in $\R^k$. As
$\{J\circ S(e_j):j=1,2,\dots,k\}$
are positive elements with sum $1$ in $\text{Aff}\, T(B)$, there exist
a positive integer $l$ and positive elements
$x_1,x_2,\dots,x_k\in \text{Aff}\, T(B_l)$ such that $\sum_{j=1}^k x_j=1$
and 
$$
\|\widehat{\beta_{l,\infty}}(x_j)-J\circ S(e_j)\|<\frac{\epsilon}{2k},
\quad j=1,2,\dots,k.
$$
Define linear positive order unit preserving maps
$V:\R^k\to \text{Aff}\, T(B_l)$ by
$$
V(\sum_{j=1}^k t_je_j)=\sum_{j=1}^k t_jx_j,
$$
and $W:\text{Aff}\, T(C)\to\text{Aff}\, T(B_l)$ by $W=V\circ T$.
Since
$$
\|\widehat{\beta_{l,\infty}}\circ V-J\circ S\|<\frac{\epsilon}2
$$
we see that
$$
\|\widehat{\beta_{l,\infty}}\circ W(f)-J(f)\|<\epsilon,\quad f\in F.
$$
By continuity of $KL(C,\cdot)$ there exist an integer $m$ and an
element $\nu\in KL(C,B_m)$ such that $[\beta_{m,\infty}]\cdot\nu=\kappa$. As
$$
{\beta_{m,\infty}}_*\circ\nu_*[1]=\kappa_*[1]=[1]=
{\beta_{m,\infty}}_*[1]\quad\text{in\ }K_0(B)
$$
we see that there exists an integer $n\ge m,l$ such that
$[\beta_{m,n}]\cdot\nu\in KL(C,B_n)_e$.
Set $\omega=[\beta_{m,n}]\cdot\nu$ and $M=\widehat{\beta_{l,n}}\circ W$.
\end{proof}
\end{lemma}

\begin{prop}\label{Mcom}
Let $A$ be a simple unital inductive limit of a sequence of finite direct sums
of building blocks. Let $B$ be the inductive limit of a sequence
$$
  \begin{CD}
      B_1  @> \beta_1 >> B_2 @> \beta_2 >>
B_3 @>\beta_3  >> \dots \\   
  \end{CD}
$$
of finite direct sums of building blocks with unital connecting maps.
Assume that there exist a $\kappa\in KL(A,B)_e$ and an affine continuous map
$\varphi_T : T(B)\to T(A)$ such that
$$
r_B(\omega)(\kappa_*(x))=r_A(\varphi_T(\omega))(x),
\quad x\in K_0(A),\ \omega\in T(B).
$$
Let $C$ be a finite direct sum of building blocks and let $\psi:C\to A$ be
a unital *-homomorphism. Let $\epsilon>0$ and let
$F\subseteq \text{Aff}\, T(C)$ be a finite subset. There exist a
positive integer $m$ and a linear positive order unit preserving map
$M:\text{Aff}\, T(C)\to \text{Aff}\, T(B_m)$ such that
$$
\|\widehat{\beta_{m,\infty}}\circ M(f)-
\varphi_{T_*}\circ\widehat{\psi}(f)\|<\epsilon,\quad f\in F,
$$
and an element $\omega\in KL(C,B_m)_e$ such that
\begin{gather*}
[\beta_{m,\infty}]\cdot\omega=\kappa\cdot [\psi]
\quad\text{in\ }KL(C,B), \\
M\circ \rho_{C}=\rho_{B_m}\circ \omega_*\quad\text{on\ }K_0(C).
\end{gather*}

\begin{proof}
We may assume that $\|f\|\le 1$, $f\in F$. Decompose
$C=C_1\oplus C_2\oplus\dots\oplus C_N$ as a finite direct sum of building
blocks. Let $r_1,r_2,\dots,r_N\in C$ be projections such that
$[r_1],[r_2],\dots,[r_N]$ generate $K_0(C)$.
By factoring $\psi$ through the $C^*$-algebra obtained from $C$ by erasing
those direct summands $C_i$ for which $\psi(r_i)=0$, we may assume that
$\psi(r_i)\neq 0$, $i=1,2,\dots,N$.
There exist positive integers $d_1,d_2,\dots,d_N$ such that
$$
\sum_{i=1}^N d_i[r_i]=[1]\quad\text{in\ }K_0(C).
$$
Since $A$ is simple there exists a $\delta_0>0$ such that
$$
\widehat{\psi}(\widehat{r_i})>\delta_0,\quad i=1,2,\dots,N.
$$
Choose $\delta>0$ such that $\delta<\delta_0$ and
$\delta(1+\sum_{i=1}^N d_i)<\epsilon$.

By Lemma \ref{KKM} there exist a positive integer $l$ and a linear positive
order unit preserving map $V:\text{Aff}\, T(C)\to \text{Aff}\, T(B_l)$ such
that
$$
\|\widehat{\beta_{l,\infty}}\circ V(f)-
{\varphi_T}_*\circ\widehat{\psi}(f)\|<\delta,
\quad f\in F\cup \{\widehat{r_1},\widehat{r_2},\dots,\widehat{r_N}\},
$$
and an element $\nu\in KK(C,B_l)_e$ such that
$$
[\beta_{l,\infty}]\cdot\nu=\kappa\cdot [\psi]\quad\text{in\ } KL(C,B).
$$

Since by assumption $\rho_B\circ\kappa_*={\varphi_T}_*\circ\rho_A$
on $K_0(A)$ we see that for $i=1,2,\dots,N$,
$$
\widehat{\beta_{l,\infty}}\circ\rho_{B_l}\circ\nu_*[r_i]=
\rho_B\circ{\beta_{l,\infty}}_*\circ\nu_*[r_i]=
{\varphi_T}_*\circ\rho_A\circ{\psi}_*[r_i]=
{\varphi_T}_*\circ\widehat{\psi}(\widehat{r_i})>\delta_0.
$$
Hence
$$
\|\widehat{\beta_{l,\infty}}\circ\rho_{B_l}\circ\nu_*[r_i]-
\widehat{\beta_{l,\infty}}\circ V(\widehat{r_i})\|<\delta,\quad i=1,2,\dots,N.
$$
Choose $m>l$ such that for $i=1,2,\dots,N$,
\begin{gather*}
\widehat{\beta_{l,m}}\circ\rho_{B_l}\circ\nu_*[r_i]>\delta_0, \\
\|\widehat{\beta_{l,m}}\circ\rho_{B_l}\circ\nu_*[r_i]-
\widehat{\beta_{l,m}}\circ V(\widehat{r_i})\|<\delta.
\end{gather*}
Define $W:\text{Aff}\, T(C)\to \text{Aff}\, T(B_m)$ by
$W=\widehat{\beta_{l,m}}\circ V$. Define $\omega\in KK(C,B_m)_e$ by
$\omega=[\beta_{l,m}]\cdot\nu$.

Decompose $B_m=B^m_1\oplus B^m_2\oplus\dots\oplus B^m_L$ as
a finite direct sum of building blocks and let $\pi_j:B_m\to B_j^m$ be the
projection, $j=1,2,\dots,L$.
Identify $\text{Aff}\, T(B_m)$ with $\oplus_{j=1}^L C_{\R}(\T)$.
Fix some $j=1,2,\dots,L$.
Set $W_j=\widehat{\pi_j}\circ W$. $W_j(\widehat{r_i})$ is a strictly positive
function in $C_{\R}(\T)$, since $\delta<\delta_0$.
Thus for each $i=1,2,\dots,N$, we can define
$M_j:\text{Aff}\, T(A_n)\cong \oplus_{i=1}^N \text{Aff}\, T(C_i)\to
C_{\R}(\T)$ by
$$
M_j(f_1,f_2,\dots,f_N)=\sum_{i=1}^N W_j(0,\dots,0,f_i,0,\dots,0)
\frac 1{W_j(\widehat{r_i})}\widehat{\pi_j}(\rho_{B_m}\circ\omega_*[r_i]).
$$
$M_j$ is positive and linear, and it preserves the order unit since
$$
M_j(1)=\sum_{i=1}^N W_j(d_i\widehat{r_i})\frac 1{W_j(\widehat{r_i})}
\widehat{\pi_j}(\rho_{B_m}\circ\omega_*[r_i])=
\sum_{i=1}^N \widehat{\pi_j}(\rho_{B_m}\circ\omega_*(d_i[r_i]))=1.
$$
Let now $g\in C_{\R}(\T)\cong \text{Aff}\, T(C_i)$, $\|g\|\le 1$,
for $i=1,2,\dots,N$. Since
$$
-d_i\widehat{r_i}\le(0,\dots,0,g,0,\dots,0)\le d_i\widehat{r_i}
$$
in $\text{Aff}\, T(C)$ we have that
\begin{align*}
&\|M_j(0,\dots,g,\dots,0)-W_j(0,\dots,g,\dots,0)\| \\
=\ &\|W_j(0,\dots,g,\dots,0)\frac 1{W_j(\widehat{r_i})}
\bigr(\widehat{\pi_j}(\rho_{B_m}\circ\omega_*[r_i])-
{W_j(\widehat{r_i})}\bigl)\| \\
\le\ &d_i\|\widehat{\pi_j}(\rho_{B_m}\circ\omega_*[r_i])-W_j(\widehat{r_i})\|<
\delta d_i.
\end{align*}
Hence if $f\in\text{Aff}\, T(C)$, $\|f\|\le 1$, then
$$
\|M_j(f)-W_j(f)\|<\sum_{i=1}^N \delta d_i.
$$

Define $M:\text{Aff}\, T(C)\to \text{Aff}\, T(B_m)$ by
$$
M(f)=(M_1(f),M_2(f),\dots,M_L(f)).
$$
Then
$$
\|M(f)-W(f)\|<\sum_{i=1}^N \delta d_i,\quad f\in \text{Aff}\, T(C),\ 
\|f\|\le 1,
$$
and hence
$$
\|\widehat{\beta_{m,\infty}}\circ M(f)-
{\varphi_T}_*\circ\widehat{\psi}(f)\|<
\delta+\sum_{i=1}^N \delta d_i<\epsilon,\quad f\in F.
$$
Finally, $M(\widehat{r_i})=\rho_{B_m}\circ\omega_*[r_i]$, $i=1,2,\dots,N$.
It follows that $M\circ\rho_C=\rho_{B_m}\circ \omega_*$ on $K_0(C)$.
\end{proof}
\end{prop}

\begin{lemma}\label{torK_0}
Let $A$ be a unital simple inductive limit of a sequence of finite direct sums
of building blocks with $K_0(A)$ non-cyclic. Then
$\text{Aff}\, T(A)/\overline{\rho_A(K_0(A))}$ is torsion free.

\begin{proof}
The image of the canonical map $K_0(A)\to\text{Aff}\, SK_0(A)$ is dense
by \cite[Proposition 3.1]{BTAF}, since $K_0(A)$ is a simple countable
dimension group. By definition $\rho_A$ is the composition of this map with
the linear bounded map $\text{Aff}\, SK_0(A)\to\text{Aff}\, T(A)$ induced by
$r_A$. It follows that $\rho_A(K_0(A))$ is dense in some subspace of
$\text{Aff}\, T(A)$.
\end{proof}
\end{lemma}

\begin{lemma}\label{eta}
Let $A$ be an inductive limit of a sequence of finite direct sums of
building blocks
$$
  \begin{CD}
      A_1  @> \alpha_1 >> A_2 @> \alpha_2 >>
A_3 @>\alpha_3  >> \dots \\   
  \end{CD}
$$
with unital connecting maps. Assume that $\rho_A$ is injective and
that $\rho_A(K_0(A))$ is a discrete subgroup of $\text{Aff}\, T(A)$. Let $n$ be
a positive integer and let $x,y$ be elements of the torsion subgroup of
$U(A_n)/\overline{DU(A_n)}$ such that
$\alpha_{n,\infty}^{\#}(x)=\alpha_{n,\infty}^{\#}(y)$. There exists an
integer $k\ge n$ such that $\alpha_{n,k}^{\#}(x)=\alpha_{n,k}^{\#}(y)$.

\begin{proof}
Since ${\alpha_{n,\infty}}_*(\pi_{A_n}(x))={\alpha_{n,\infty}}_*(\pi_{A_n}(y))$
in $K_1(A)$ there is an integer $l\ge n$ such that
${\alpha_{n,l}}_*(\pi_{A_n}(x))={\alpha_{n,l}}_*(\pi_{A_n}(y))$. By Proposition
\ref{exact} we see that
$$
\alpha_{n,l}^{\#}(x-y)=\lambda_{A_l}(q_{A_l}(\frac 1m\rho_{A_l}(z)))
$$
for some positive integer $m$ and an element $z\in K_0(A_l)$. Since
$\rho_A(K_0(A))$ is discrete and since
$\lambda_A(q_A(\frac 1m\rho_A({\alpha_{l,\infty}}_*(z))))=0$ we see that
$\frac 1m\rho_A({\alpha_{l,\infty}}_*(z))=\rho_A({\alpha_{j,\infty}}_*(w))$
for some positive integer $j$ and an element $w\in K_0(A_j)$. Since $\rho_A$
is injective we may choose an integer $k\ge l,j$ such that
${\alpha_{l,k}}_*(z)={\alpha_{j,k}}_*(mw)$ in $K_0(A)$. Note that
$\alpha_{n,k}^{\#}(x-y)=
\lambda_{A_k}(q_{A_k}(\frac 1m\rho_{A_k}({\alpha_{k,l}}_*(z))))=0$.
\end{proof}
\end{lemma}

\begin{prop}\label{homotopy}
Let $A$ be a unital $C^*$-algebra and let $B$ be a unital inductive limit
of a sequence of finite direct sums
of building blocks such that the torsion subgroup of
$\text{Aff}\, T(B)/\overline{\rho_B(K_0(B))}$ is totally disconnected.
Let $\varphi,\psi:A\to B$ be unital *-homomorphisms that
are homotopic and let $x\in U(A)/\overline{DU(A)}$ be an element of finite
order. Then $\varphi^{\#}(x)=\psi^{\#}(x)$.

\begin{proof}
Let $u\in A$ be a unitary such that $x=q_A'(u)$.
Let $(\varphi_t)_{t\in [0,1]}$ be a homotopy connecting $\varphi$ to
$\psi$. We may assume that $\|\varphi_t(u)-\varphi_0(u)\|<1$ for
$t\in [0,1]$. Thus
$$
\varphi_t(u)\varphi_0(u)^*=e^{2\pi ib_t}
$$
where $t\mapsto b_t$ is a continuous path of self-adjoint elements in $B$.
Since $\lambda_B(q_B(\widehat{b_t}))=q_B'(e^{2\pi ib_t})$ we see that
$q_B(\widehat{b_t})$ has finite order in
$\text{Aff}\, T(B)/\overline{\rho_B(K_0(B))}$.
Thus $t\mapsto q_B(\widehat{b_t})$ is a continuous path in a totally
disconnected subset of a metric space. It follows
that it is constant and hence $q_B(\widehat{b_t})=0$
for every $t\in [0,1]$. We conclude that
$\varphi_0^{\#}(q_A'(u))=\varphi_1^{\#}(q_A'(u))$.
\end{proof}
\end{prop}

We leave it as an open question whether the torsion subgroup of the group
$\text{Aff}\, T(B)/\overline{\rho_B(K_0(B))}$ always is
totally disconnected.

\begin{prop}\label{KLu}
Let $A$ be a finite direct sum of building blocks and let $B$ be a unital
inductive limit of a sequence of finite direct sums
of building blocks such that the torsion subgroup of
$\text{Aff}\, T(B)/\overline{\rho_B(K_0(B))}$ is totally disconnected.
Let $\varphi,\psi:A\to B$ be unital *-homomorphisms such that
$[\varphi]=[\psi]$ in $KL(A,B)$. Let $x$ be an element of the torsion subgroup
of $U(A)/\overline{DU(A)}$. Then $\varphi^{\#}(x)=\psi^{\#}(x)$.

\begin{proof}
By \cite [Corollary 15.1.3]{LLS} and Theorem \ref{stable} there exist a
positive integer $m$ and
*-homomorphisms $\lambda,\mu:A\to B_m$ such that $\varphi$ is homotopic to
$\beta_{m,\infty}\circ\lambda$ and $\psi$ is homotopic to
$\beta_{m,\infty}\circ\mu$. By increasing $m$ we may assume that $\lambda$
and $\mu$ are unital. There exists an integer $k\ge m$ such that
$[\beta_{m,k}]\cdot[\lambda]=[\beta_{m,k}]\cdot[\mu]$ in $KL(A,B_k)$.
Thus $\beta_{m,k}^{\#}\circ\lambda^{\#}(x)=\beta_{m,k}^{\#}\circ\mu^{\#}(x)$
by Proposition \ref{KKunid}. Hence
$\varphi^{\#}(x)=\beta_{m,\infty}^{\#}\circ\lambda^{\#}(x)=
\beta_{m,\infty}^{\#}\circ\mu^{\#}(x)=\psi^{\#}(x)$ by
Proposition \ref{homotopy}.
\end{proof}
\end{prop}

\begin{lemma}\label{K_0nct}
Let $A$ be a simple unital inductive limit of a sequence of finite direct sums
of building blocks such that $K_0(A)$ is non-cyclic, and let $B$ be a unital
inductive limit of a sequence of finite direct sums of building blocks.
If there exists an element $\kappa\in KL(A,B)_T$ then
$\text{Aff}\, T(B)/\overline{\rho_B(K_0(B))}$ is torsion free.

\begin{proof}
By Lemma \ref{unital} we may assume that $A$ is the inductive limit of a
sequence
$$
  \begin{CD}
      A_1  @> \alpha_1 >> A_2 @> \alpha_2 >>
A_3 @>\alpha_3  >> \dots \\   
  \end{CD}
$$
of finite direct sums of building blocks with unital connecting
maps. Similarly $B$ is the inductive limit of a sequence of finite direct sums
of building blocks
$$
  \begin{CD}
      B_1  @> \beta_1 >> B_2 @> \beta_2 >>
B_3 @>\beta_3  >> \dots \\   
  \end{CD}
$$
with unital connecting maps.
Let $\epsilon>0$. There exists a positive integer $n$ such that for every
$t\in\R$ we have that $d_{A_n}'(q_{A_n}(t\widehat{1}),0)<\epsilon$. To see
this choose a positive integer $k$ such that $\frac 1k<\epsilon$.
Since $\text{Aff}\, T(A)/\overline{\rho_A(K_0(A))}$ is torsion free by
Lemma \ref{torK_0}, we may choose $n$ such that
$d_{A_n}'(q_{A_n}(\frac jk\widehat{1}),0)<\frac{\epsilon}2$, $j=1,2,\dots,k-1$.
Let $t\in\R$. We may assume that $0<t<1$. Choose $j=0,1,2,\dots,k$ such that
$|t-\frac jk|\le \frac 1{2k}<\frac{\epsilon}2$. Then
$d_{A_n}'(q_{A_n}(t\widehat{1}),0)<\epsilon$.

By Proposition \ref{Mcom} we get a positive integer $l$ and a
contractive group homomorphism
$S:\text{Aff}\, T(A_n)/\overline{\rho_{A_n}(K_0(A_n))}\to
\text{Aff}\, T(B_l)/\overline{\rho_{B_l}(K_0(B_l))}$ such that
$S(q_{A_n}(r\widehat{1}))=q_{B_l}(r\widehat{1})$ for every $r\in\R$.
Let $x\in \text{Aff}\, T(B)/\overline{\rho_{B}(K_0(B))}$ be an element of
order $m$. There is an integer $k\ge l$
such that $d_B'(x,q_B(\frac 1m\rho_B({\beta_{k,\infty}}_*(y))))<\epsilon$
for some element $y\in K_0(B_k)$.
We claim that $d_{B_k}'(q_{B_k}(\frac 1m\rho_{B_k}(y)),0)<\epsilon$.
To this end we may assume that $B_k$ is a building block. Then
$\rho_{B_k}(y)=w\widehat{1}$ for some $w\in\Q$. Hence
$$
d_{B_k}'(q_{B_k}(\frac 1m\rho_{B_k}(y)),0)=
d_{B_k}'(\widetilde{\beta_{l,k}}\circ S(q_{A_n}(\frac wm\widehat{1})),0)\le
d_{A_n}'(q_{A_n}(\frac wm\widehat{1}),0)<\epsilon.
$$
Thus $d_B'(x,0)<2\epsilon$. Since $\epsilon>0$ was arbitrary we conclude that
$x=0$.
\end{proof}
\end{lemma}

\begin{lemma}\label{konttor}
Let $A$ be a simple inductive limit of a sequence
$$
  \begin{CD}
      A_1  @> \alpha_1 >> A_2 @> \alpha_2 >>
A_3 @>\alpha_3  >> \dots \\   
  \end{CD}
$$
of finite direct sums of building blocks with unital connecting maps. Let
$y$ be an element in $U(A)/\overline{DU(A)}$ of order $k<\infty$. Then there
exist a positive integer $m$ and an element $w\in U(A_m)/\overline{DU(A_m)}$
of order $k$ such that $\alpha_{m,\infty}^{\#}(w)=y$.

\begin{proof}
By continuity of $K_1$ there exist a positive integer $l$ and an element
$z$ in $U(A_l)/\overline{DU(A_l)}$ such that
${\alpha_{l,\infty}}_*(\pi_{A_l}(z))=\pi_A(y)$ in $K_1(A)$. Since the
short exact sequence of Proposition \ref{exact} splits we may assume that
$kz=0$. Note that $\pi_A(\alpha_{l,\infty}^{\#}(z))=\pi_A(y)$
and hence
$$
y=\alpha_{l,\infty}^{\#}(z)+\lambda_A(q_A(f))
\quad\text{in}\ U(A)/\overline{DU(A)}
$$
for some $f\in\text{Aff}\, T(A)$ with $kq_A(f)=0$ in the group
$\text{Aff}\, T(A)/\overline{\rho_A(K_0(A))}$. If $K_0(A)$ is non-cyclic then
we see that $q_A(f)=0$ by Lemma \ref{torK_0}. Thus we
may assume that $K_0(A)\cong\Z$ such that $\rho_A(K_0(A))$ is a discrete
subgroup of $\text{Aff}\, T(A)$. It follows that
$f=\frac 1k\rho_A(x)$ for some $x\in K_0(A)$. By continuity of $K_0$ we have
that $x={\alpha_{m,\infty}}_*(h)$ for some integer $m\ge l$ and some
$h\in K_0(A_m)$. Define $w\in U(A_m)/\overline{DU(A_m)}$ by
$$
w=\alpha_{l,m}^{\#}(z)+\lambda_{A_m}(q_{A_m}(\frac 1k\rho_{A_m}(h))).
$$
Then
$$
\alpha_{m,\infty}^{\#}(w)=\alpha_{l,\infty}^{\#}(z)+
\lambda_A(q_A(\frac 1k\rho_A({\alpha_{m,\infty}}_*(h))))
=\alpha_{l,\infty}^{\#}(z)+\lambda_A(q_A(f))=y.
$$
Since $y$ has order $k$ and $kw=0$ it follows that $w$ has order $k$ as well.
\end{proof}
\end{lemma}

\begin{thm}
Let $A$ be a unital simple inductive limit of a sequence
of finite direct sums of building blocks and let $B$ be an inductive limit
of a similar sequence
$$
  \begin{CD}
      B_1  @> \beta_1 >> B_2 @> \beta_2 >>
B_3 @>\beta_3  >> \dots \\   
  \end{CD}
$$
with unital connecting maps such that $s(B_k)\to\infty$ and such that the
torsion subgroup of $\text{Aff}\, T(B)/\overline{\rho_B(K_0(B)}$ is totally
disconnected. Let $\kappa\in KL(A,B)_T$.
Let $C$ be a finite direct sum of building blocks and let $\varphi:C\to A$ be
a unital *-homomorphism. Then there is a unital *-homomorphism $\psi:C\to B$
such that $[\psi]=\kappa\cdot[\varphi]$ in $KL(C,B)$.

Moreover, if $C_1$ is another finite direct sum of building blocks,
if $\varphi_1:C_1\to A$ and $\psi_1:C_1\to B$
are unital *-homomorphisms such that $[\psi_1]=\kappa\cdot[\varphi_1]$ in
$KL(A,B)$, and if $x\in U(C)/\overline{DU(C)}$ and
$x_1\in U(C_1)/\overline{DU(C_1)}$ are elements of finite order
such that $\varphi^{\#}(x)=\varphi_1^{\#}(x_1)$, then
$\psi^{\#}(x)=\psi_1^{\#}(x_1)$.

\begin{proof}
Let a finite direct sum of building blocks $C$ and a unital
*-homomorphism $\varphi:C\to A$ be given.
Let $\varphi_T:T(B)\to T(A)$ be a continuous affine map such that
$$
r_B(\omega)(\kappa_*(x))=r_A(\varphi_T(\omega))(x),
\quad x\in K_0(A),\ \omega\in T(B).
$$
Let $p_1,p_2,\dots,p_N$ be the
minimal non-zero central projections in $C$. As in the proof of Proposition
\ref{Mcom} we see that we may assume that $\varphi(p_i)\neq 0$,
$i=1,2,\dots,N$. Choose $\delta>0$ such that
$\widehat{\varphi}(\widehat{p_i})>2\delta$. 
Choose an integer $K$ by Theorem \ref{ex} with respect to
$F=\emptyset$ and $\epsilon=1$. By Proposition \ref{Mcom} there
exist a positive integer $m$ and a linear positive order unit preserving map
$M:\text{Aff}\, T(C)\to \text{Aff}\, T(B_m)$ such that
$$
\|\widehat{\beta_{m,\infty}}\circ M(\widehat{p_i})-
\varphi_{T_*}\circ\widehat{\varphi}(\widehat{p_i})\|<\delta,
\quad i=1,2,\dots,N,
$$
and an element $\omega\in KL(C,B_m)_e$ such that
\begin{gather*}
[\beta_{m,\infty}]\cdot\omega=\kappa\cdot [\varphi]
\quad\text{in\ }KL(C,B), \\
M\circ \rho_{C}=\rho_{B_m}\circ \omega_*\quad\text{on\ }K_0(C).
\end{gather*}
Hence $\widehat{\beta_{m,\infty}}\circ M(\widehat{p_i})>\delta$,
$i=1,2,\dots,N$. Choose $k\ge m$ such that $s(B_k)\ge K\delta^{-1}$ and
such that $\widehat{\beta_{m,k}}\circ M(\widehat{p_i})>\delta$,
$i=1,2,\dots,N$. Then
$\rho_{B_k}({\beta_{m,k}}_*\circ\omega_*[p_i])>\delta$
and hence $s(B_k)\rho_{B_k}({\beta_{m,k}}_*\circ\omega_*[p_i])\ge K$.
Furthermore $\widehat{\beta_{m,k}}\circ M\circ\rho_C=
\rho_{B_k}\circ{\beta_{m,k}}_*\circ\omega_*$.
It follows from Theorem \ref{ex} that there exists a unital *-homomorphism
$\mu:C\to B_k$ such that $[\mu]=[\beta_{m,k}]\cdot\omega$. Set
$\psi=\beta_{k,\infty}\circ\mu$. This proves the first part of the theorem.

To prove the second part of the theorem, let us first note that
$$
\pi_B(\psi^{\#}(x))=\psi_*(\pi_C(x))=\kappa_*\circ\varphi_*(\pi_C(x))=
{\psi_1}_*(\pi_{C_1}(x_1))=\pi_B(\psi_1^{\#}(x_1)).
$$
Hence if $K_0(A)$ is non-cyclic then $\psi^{\#}(x)=\psi_1^{\#}(x_1)$
by Lemma \ref{K_0nct}.

We may therefore assume that $K_0(A)$ is cyclic.
By Lemma \ref{unital} we see that $A$ is the inductive limit of a
sequence
$$
  \begin{CD}
      A_1  @> \alpha_1 >> A_2 @> \alpha_2 >>
A_3 @>\alpha_3  >> \dots \\   
  \end{CD}
$$
where each $A_n$ is a finite direct sum of building blocks and each $\alpha_n$
is unital. By \cite[Corollary 15.1.3]{LLS} there exist a positive integer $n$
and *-homomorphisms $\lambda:C\to A_n$ and $\lambda_1:C_1\to A_n$ such that
$\varphi$ is homotopic to $\alpha_{n,\infty}\circ\lambda$ and $\varphi_1$ is
homotopic to $\alpha_{n,\infty}\circ\lambda_1$. Note that $\lambda$ and
$\lambda_1$ are unital. Since
$$
\alpha_{n,\infty}^{\#}\circ\lambda^{\#}(x)=
\alpha_{n,\infty}^{\#}\circ\lambda_1^{\#}(x_1)
$$
by Proposition \ref{homotopy}, there exists by Lemma \ref{eta} a positive
integer $k$ such that
$$
\alpha_{n,k}^{\#}\circ\lambda^{\#}(x)=
\alpha_{n,k}^{\#}\circ\lambda_1^{\#}(x_1).
$$
By the first part of the theorem there is a unital *-homomorphism
$\gamma:A_k\to B$ such that
$[\gamma]=\kappa\cdot [\alpha_{k,\infty}]$. Note that
\begin{gather*}
[\gamma]\cdot[\alpha_{n,k}]\cdot [\lambda]=
\kappa\cdot [\alpha_{n,\infty}]\cdot[\lambda]=\kappa\cdot[\varphi]=[\psi]
\quad\text{in}\ KL(C,B)  \\
[\gamma]\cdot[\alpha_{n,k}]\cdot [\lambda_1]=
\kappa\cdot [\alpha_{n,\infty}]\cdot[\lambda_1]=
\kappa\cdot[\varphi_1]=[\psi_1]\quad\text{in}\ KL(C_1,B) .
\end{gather*}
Hence
$$
\psi^{\#}(x)=\gamma^{\#}\circ\alpha_{n,k}^{\#}\circ\lambda^{\#}(x)=
\gamma^{\#}\circ\alpha_{n,k}^{\#}\circ\lambda_1^{\#}(x_1)=\psi_1^{\#}(x_1)
$$
by Proposition \ref{KLu}.
\end{proof}
\end{thm}

Let $A$, $B$ and $\kappa$ be as above. Let $y$ be an element in
$U(A)/\overline{DU(A)}$ of finite order. By Lemma \ref{konttor} there is a
finite direct sum of building blocks $C$, an element of finite order $x$ in
$U(C)/\overline{DU(C)}$, and a unital *-homomorphism $\varphi:C\to A$
such that $\varphi^{\#}(x)=y$. By the first part of the theorem above there
exists a unital *-homomorphism $\psi:C\to B$ such that
$[\psi]=\kappa\cdot[\varphi]$. Set $s_{\kappa}(y)=\psi^{\#}(x)$.
By the second part $s_{\kappa}(y)$ is independent of the choice of $\varphi$,
$\psi$ and $x$. Thus we have a well-defined map
$$
s_{\kappa}:Tor\bigl(U(A)/\overline{DU(A)}\bigr)
\to Tor \bigl(U(B)/\overline{DU(B)}\bigr).
$$
It follows easily from Lemma \ref{konttor} that $s_{\kappa}$ is a group
homomorphism. Note that if $\mu:A\to B$ is a unital *-homomorphism then
$s_{[\mu]}(y)=\mu^{\#}(y)$ for every $y$ in the torsion subgroup of
$U(A)/\overline{DU(A)}$. Finally, we note that $s_{\kappa}$ exists
for trivial reasons if $K_0(A)$ is non-cyclic (since
$\text{Aff}\, T(B)/\overline{\rho_B(K_0(B))}$ is torsion free in this case,
see Lemma \ref{K_0nct}). It is possible (as in \cite{TATD}) to prove our
classification theorem in the case of non-cyclic $K_0$-group without using the
map $s_{\kappa}$, but we have chosen to construct it in general in order to
obtain a unified proof of the classification theorem in the cases $K_0$ cyclic
and $K_0$ non-cyclic.

\begin{lemma}\label{lifttophi}
Let $A$ be a unital simple infinite dimensional inductive limit of a sequence
of finite direct sums of building blocks and let $B$ be an inductive limit
of a similar sequence
$$
  \begin{CD}
      B_1  @> \beta_1 >> B_2 @> \beta_2 >>
B_3 @>\beta_3  >> \dots \\   
  \end{CD}
$$
with unital connecting maps such that $s(B_k)\to\infty$ and such that the
torsion subgroup of $\text{Aff}\, T(B)/\overline{\rho_B(K_0(B))}$ is totally
disconnected. Let $\kappa\in KL(A,B)_e$ and let
$\varphi_T:T(B)\to T(A)$ be a continuous affine map such that
$$
r_B(\omega)(\kappa_*(x))=r_A(\varphi_T(\omega))(x),
\quad x\in K_0(A),\ \omega\in T(B).
$$
There exists a group homomorphism
$\Phi:U(A)/\overline{DU(A)}\to U(B)/\overline{DU(B)}$ such that
$\Phi(y)=s_{\kappa}(y)$ for $y$ in the torsion subgroup of
$U(A)/\overline{DU(A)}$ and such that the diagram
$$
  \begin{CD}
      0 @> >>\text{Aff}\, T(A)/\overline{\rho_A(K_0(A))} @> \lambda_A >>
      U(A)/\overline{DU(A)}   @> \pi_A >> K_1(A) @> >> 0   \\   
      @. @V\widetilde{\varphi_T} VV         @V\Phi VV      @VV\kappa_* V \\
      0 @> >> \text{Aff}\, T(B)/\overline{\rho_B(K_0(B))} @>> \lambda_B >
      U(B)/\overline{DU(B)}   @>> \pi_B > K_1(B) @> >> 0
  \end{CD}
$$
commutes.

\begin{proof}
It will be convenient to set $G_1=\text{Aff}\, T(A)/\overline{\rho_A(K_0(A))}$,
$G_2=K_1(A)$, and $H_1=\text{Aff}\, T(B)/\overline{\rho_B(K_0(B))}$,
$H_2=K_1(B)$. Note that $U(A)/\overline{DU(A)}\cong G_1\oplus G_2$ and
$U(B)/\overline{DU(B)}\cong H_1\oplus H_2$ by Proposition \ref{exact}.
Hence $s_{\kappa}$ can be identified with a matrix of the form
$$
\begin{pmatrix}
f_{11} & f_{12} \\
f_{21} & f_{22}
\end{pmatrix},
$$
where $f_{ij}:Tor(G_j)\to Tor(H_i)$ is a group homomorphism, $i,j=1,2$.

Let $z\in Tor(G_1)$. If $K_0(A)$ is cyclic then $z=q_A(\frac 1m\rho_A(h))$ for
some positive integer $m$ and $h\in K_0(A)$. Choose
a finite direct sum of building blocks $C$ and a unital *-homomorphism
$\varphi:C\to A$ such that $\varphi_*(g)=h$ for some $g\in K_0(C)$.
Choose a unital *-homomorphism
$\psi:C\to B$ such that $[\psi]=\kappa\cdot[\varphi]$. Since
${\varphi_T}_*\circ\rho_A=\rho_B\circ\kappa_*$ we see that
\begin{align*}
s_{\kappa}(\lambda_A(z))&=
s_{\kappa}(\lambda_A(q_A(\frac 1m\rho_A(\varphi_*(g)))))=
s_{\kappa}(\varphi^{\#}(\lambda_C(q_C(\frac 1m\rho_C(g))))) \\
&=\psi^{\#}(\lambda_C(q_C(\frac 1m\rho_C(g))))=
\lambda_B(q_B(\frac 1m\rho_B(\psi_*(g)))) \\
&=\lambda_B(q_B(\frac 1m{\varphi_T}_*(\rho_A(\varphi_*(g)))))=
\lambda_B(\widetilde{\varphi_T}(z)).
\end{align*}
Hence $f_{11}(z)=\widetilde{\varphi_T}(z)$ and $f_{21}(z)=0$.
By Lemma \ref{K_0nct} this conclusion also holds if $K_0(A)$ is non-cyclic.
Let $w\in Tor(G_2)$. Choose an element $y\in U(A)/\overline{DU(A)}$
of finite order such that $\pi_A(y)=w$. Choose a finite direct sum of
building blocks $C$ and a unital *-homomorphism $\varphi:C\to A$ such that
$\varphi^{\#}(x)=y$. Choose a unital *-homomorphism $\psi:C\to B$ such that
$[\psi]=\kappa\cdot[\varphi]$ in $KL(C,B)$. Since
$$
\pi_B(s_{\kappa}(y))=\pi_B(\psi^{\#}(x))=\psi_*(\pi_C(x))=
\kappa_*\circ\pi_A(\varphi^{\#}(x))=\kappa_*\circ\pi_A(y)
$$
we see that $f_{22}(w)=\kappa_*(w)$. Finally, since $H_1$ is a divisible
group there exists by \cite[Theorem 21.1]{FIAG} a group homomorphism
$\lambda:G_2\to H_1$ such that
$\lambda(w)=f_{12}(w)$ for every $w\in G_2$ of finite order. Set
$$
\Phi=\begin{pmatrix}
\widetilde{\varphi_T} & \lambda \\
0 & \kappa_*
\end{pmatrix}.
$$
It is easy to see that the diagram commutes.
\end{proof}
\end{lemma}

\section{Main results}

Consider the category of abelian groups, equipped with a complete and
translation invariant metric, and contractive group homomorphisms. Inductive
limits can be constructed in this category in a way similar to the way that
they are constructed in the category of $C^*$-algebras. Indeed, let
$$
  \begin{CD}
      G_1  @> \mu_1 >> G_2 @> \mu_2 >>
G_3 @>\mu_3  >> \dots \\   
  \end{CD}
$$
be an inductive system. Let $\rho_k$ denote the metric on
$G_k$. Let $H$ be the inductive limit in the category of groups. Define a
pseudo-metric $d$ on $H$ by
$$
d(\mu_{n,\infty}(x),\mu_{m,\infty}(y))=
\lim_{k\to\infty} \rho_k(\mu_{n,k}(x),\mu_{m,k}(y)).
$$
Form the quotient of $H$ by the subgroup $\{x\in H:d(x,0)=0\}$ and complete
with respect to the induced metric to obtain the inductive limit.

It is an elementary exercise to prove that
$U(\cdot)/\overline{DU(\cdot)}$ is a continuous functor from the category of
unital $C^*$-algebras and unital *-homomorphisms, to the category of abelian
groups equipped with a complete translation invariant metric, and contractive
group homomorphisms.

\begin{prop}\label{exhom}
Let $A$ be a simple inductive limit of a sequence
$$
  \begin{CD}
      A_1  @> \alpha_1 >> A_2 @> \alpha_2 >>
A_3 @>\alpha_3  >> \dots \\   
  \end{CD}
$$
of finite direct sums of building blocks with unital and injective connecting
maps. Let $B$ be an inductive limit of a similar sequence
$$
  \begin{CD}
      B_1  @> \beta_1 >> B_2 @> \beta_2 >>
B_3 @>\beta_3  >> \dots \\   
  \end{CD}
$$
with unital connecting maps such that $s(B_k)\to\infty$ and such that
the torsion subgroup of $\text{Aff}\, T(B)/\overline{\rho_B(K_0(B))}$ is
totally disconnected. Let $\varphi_T : T(B)\to T(A)$ be an affine continuous
map, let $\kappa\in KL(A,B)_e$ be an element such that
$$
r_B(\omega)(\kappa_*(x))=r_A(\varphi_T(\omega))(x),
\quad x\in K_0(A),\ \omega\in T(B),
$$
and let $\Phi:U(A)/\overline{DU(A)}\to U(B)/\overline{DU(B)}$ be a
homomorphism such that the diagram
$$
  \begin{CD}
      \text{Aff}\, T(A)/\overline{\rho_A(K_0(A))} @> \lambda_A >>
      U(A)/\overline{DU(A)} @> \pi_A >> K_1(A)    \\   
      @V\widetilde{\varphi_T} VV  @V\Phi VV      @VV\kappa_* V \\
      \text{Aff}\, T(B)/\overline{\rho_B(K_0(B))}
      @>> \lambda_B > U(B)/\overline{DU(B)} @>> \pi_B > K_1(B)
  \end{CD}
$$
commutes. Assume finally that
$$
s_{\kappa}(y)=\Phi(y),\quad y\in Tor(U(A)/\overline{DU(A)}).
$$

Let $n$ be a positive integer and let $F_1\subseteq\text{Aff}\, T(A_n)$
and $F_2\subseteq U(A_n)/\overline{DU(A_n)}$ be finite sets. There exist
a positive integer $m$ and a unital *-homomorphism $\psi:A_n\to B_m$
such that
\begin{gather*}
[\beta_{m,\infty}]\cdot [\psi]=\kappa\cdot [\alpha_{n,\infty}]
\quad\text{in\ }KL(A_n,B), \\
\|\widehat{\beta_{m,\infty}}\circ\widehat{\psi}(f)-
{\varphi_T}_*\circ\widehat{\alpha_{n,\infty}}(f)\|<\epsilon,\quad f\in F_1, \\
D_B\bigl(\, \beta_{m,\infty}^{\#}\circ\psi^{\#}(x)\, , \,
\Phi\circ\alpha_{n,\infty}^{\#}(x)\, \bigr)<\epsilon,\quad x\in F_2.
\end{gather*}

\begin{proof}
Let $A_n=C_1\oplus\dots\oplus C_R$ where each $C_i$ is a
building block. By Proposition \ref{exact} and Proposition \ref{K_1c}
there are for each $x\in U(A_n)/\overline{DU(A_n)}$ an element
$a_x$ in $\text{Aff}\, T(A_n)/\overline{\rho_{A_n}(K_0(A_n))}$, integers
$k_x^1,k_x^2,\dots,k_x^R$, and an element $y_x$ in the torsion subgroup of
$U(A_n)/\overline{DU(A_n)}$ such that
$$
x=\lambda_{A_n}(a_x)+\sum_{i=1}^R k_x^i q_{A_n}'(v^{A_n}_i)+y_x
\quad\text{in}\ U(A_n)/\overline{DU(A_n)}.
$$
Choose $b_x\in\text{Aff}\, T(A_n)$ such that $q_{A_n}(b_x)=a_x$. Set
$F_1'=F_1\cup \{b_x:x\in F_2\}$. Choose $0<\delta<\frac 12$ such that
$\delta<\epsilon$ and such that
$$
|e^{2\pi i\delta}-1|+\delta \sum_{i=1}^R k_i^x<\epsilon,\quad x\in F_2.
$$
Let $p_1,p_2,\dots,p_R$ denote the minimal non-zero central projections in
$A_n$. Since $A$ is simple and the connecting maps are injective, there
exists a $\gamma>0$ such that
$\widehat{\alpha_{n,\infty}}(\widehat{p_i})>\gamma$,
$i=1,2,\dots,R$. By Proposition \ref{Mcom} there exist a positive integer $l$,
a linear positive order unit preserving map
$M:\text{Aff}\, T(A_n)\to \text{Aff}\, T(B_l)$, and an element
$\omega\in KL(A_n,B_l)_e$ such that
\begin{gather*}
[\beta_{l,\infty}]\cdot\omega=\kappa\cdot [\alpha_{n,\infty}]
\quad\text{in\ }KL(A_n,B), \\
\|\widehat{\beta_{l,\infty}}\circ M(f)-
{\varphi_T}_*\circ\widehat{\alpha_{n,\infty}}(f)\|<\frac{\delta}2,
\quad f\in F_1', \\
M\circ\rho_{A_n}=\rho_{B_l}\circ\omega_*\quad\text{on\ }K_0(A_n).
\end{gather*}
Choose an integer $K$ by Theorem \ref{ex} with respect to
$F_1'\subseteq\text{Aff}\, T(A_n)$ and $\frac{\delta}2$. Choose a positive
integer $k$ and unitaries $u_1,u_2,\dots,u_R\in B_k$ such that
$$
D_B\bigl(\, \beta_{k,\infty}^{\#}(q_{B_k}'(u_i))\, ,\, 
\Phi\circ\alpha_{n,\infty}^{\#}(q_{A_n}'(v^{A_n}_i))\, \bigr)<\delta,
\quad i=1,2,\dots,R.
$$
Note that
$\kappa_*\circ{\alpha_{n,\infty}}_*[v^{A_n}_i]={\beta_{k,\infty}}_*[u_i]$
in $K_1(B)$. Hence
$$
{\beta_{l,\infty}}_*\circ\omega_*[v^{A_n}_i]=
{\beta_{k,\infty}}_*[u_i],\quad i=1,2,\dots,R.
$$
Since $\rho_B\circ\kappa_*={\varphi_T}_*\circ\rho_A$ we see that for
$i=1,2,\dots,R$,
$$
\widehat{\beta_{l,\infty}}(\rho_{B_l}(\omega_*[p_i]))=
\rho_B({\beta_{l,\infty}}_*\circ\omega_*[p_i])=
{\varphi_T}_*\circ\rho_A\circ {\alpha_{n,\infty}}_*[p_i]=
{\varphi_T}_*\circ\widehat{\alpha_{n,\infty}}(\widehat{p_i})>\gamma.
$$
Hence there exists an integer $m\ge k,l$ such that $s(B_m)\ge K\gamma^{-1}$ and
such that
\begin{gather*}
\widehat{\beta_{l,m}}(\rho_{B_l}(\omega_*[p_i]))>\gamma,\quad i=1,2,\dots,R, \\
{\beta_{l,m}}_*\circ\omega_*[v^{A_n}_i]=
{\beta_{k,m}}_*[u_i]\quad\text{in\ }K_1(B_m),\quad i=1,2,\dots,R.
\end{gather*}
It follows that $s(B_m) \rho_{B_m}({\beta_{l,m}}_*\circ\omega_*[p_i])\ge K$
and that
$$
\widehat{\beta_{l,m}}\circ M\circ\rho_{A_n}=
\widehat{\beta_{l,m}}\circ\rho_{B_l}\circ\omega_*=
\rho_{B_m}\circ {\beta_{l,m}}_*\circ\omega_*\quad\text{on}\ K_0(A_n).
$$
Therefore by Theorem \ref{ex} there exists a unital *-homomorphism
$\psi:A_n\to B_m$ such that
\begin{gather*}
[\psi]=[\beta_{l,m}]\cdot\omega\quad\text{in\ }KL(A_n,B_m), \\
\psi^{\#}(q_{A_n}'(v^{A_n}_i))=q_{B_m}'(\beta_{k,m}(u_i))
\quad\text{in\ }U(B_m)/\overline{DU(B_m)},\quad i=1,2,\dots,R, \\
\|\widehat{\psi}(f)-\widehat{\beta_{l,m}}\circ M(f)\|<\frac{\delta}2,
\quad f\in F_1'.
\end{gather*}
It follows that
\begin{gather}\label{skeq}
[\beta_{m,\infty}]\cdot [\psi]=\kappa\cdot [\alpha_{n,\infty}]
\quad\text{in\ }KL(A_n,B), \\
\|\widehat{\beta_{m,\infty}}\circ\widehat{\psi}(f)-{\varphi_T}_*\circ
\widehat{\alpha_{n,\infty}}(f)\|<\delta,\quad f\in F_1', \\
D_B\bigl(\, \beta_{m,\infty}^{\#}\circ\psi^{\#}(q_{A_n}'(v^{A_n}_i))\, , \,
\Phi\circ\alpha_{n,\infty}^{\#}(q_{A_n}'(v^{A_n}_i))\, \bigr)<\delta,
\quad i=1,2,\dots,R.
\end{gather}
Note that for $x\in F_2$,
\begin{align*}
d_B'\bigl(\, \widetilde{\beta_{m,\infty}}\circ\widetilde{\psi}(a_x)\, ,\,
\widetilde{\varphi_T}\circ\widetilde{\alpha_{n,\infty}}(a_x)\, \bigr)=\ &
d_B'\bigl(\, q_B(\widehat{\beta_{m,\infty}}\circ\widehat{\psi}(b_x))\, ,\,
q_B({\varphi_T}_*\circ\widehat{\alpha_{n,\infty}}(b_x))\, \bigr) \\
\le\ &\|\widehat{\beta_{m,\infty}}\circ\widehat{\psi}(b_x)-
{\varphi_T}_*\circ\widehat{\alpha_{n,\infty}}(b_x)\|<\delta<\frac 12.
\end{align*}
Hence
$$
d_B\bigl(\, \widetilde{\beta_{m,\infty}}\circ\widetilde{\psi}(a_x)\, ,\,
\widetilde{\varphi_T}\circ\widetilde{\alpha_{n,\infty}}(a_x)\, \bigr)<
|e^{2\pi i\delta}-1|,\quad x\in F_2.
$$
By Proposition \ref{exact}, $\lambda_B$ is an isometry when
$\text{Aff}\, T(B)/\overline{\rho_B(K_0(B))}$
is equipped with the metric $d_B$. It follows that
$$
D_B\bigl(\, \lambda_B\circ\widetilde{\beta_{m,\infty}}
\circ\widetilde{\psi}(a_x)\, ,\, \lambda_B\circ\widetilde{\varphi_T}\circ
\widetilde{\alpha_{n,\infty}}(a_x)\, \bigr)<
|e^{2\pi i\delta}-1|,\quad x\in F_2.
$$
Thus
$$
D_B\bigl(\, \beta_{m,\infty}^{\#}\circ\psi^{\#}\circ\lambda_{A_n}(a_x)\, ,\, 
\Phi\circ\alpha_{n,\infty}^{\#}\circ\lambda_{A_n}(a_x)\, \bigr)<
|e^{2\pi i\delta}-1|,\quad x\in F_2.
$$
Since $s_{\kappa}$ and $\Phi$ agree on the torsion subgroup of
$U(A)/\overline{DU(A)}$, we see by (\ref{skeq}) and the definition of
$s_{\kappa}$ that
$$
\beta_{m,\infty}^{\#}\circ\psi^{\#}(y_x)=
\Phi\circ\alpha_{n,\infty}^{\#}(y_x).
$$
Hence for $x\in F_2$,
\begin{align*}
&D_B\bigl(\, \beta_{m,\infty}^{\#}\circ\psi^{\#}(x)\, , \,
\Phi\circ\alpha_{n,\infty}^{\#}(x)\, \bigr) \\
\le\ & 
D_B\bigl(\, \beta_{m,\infty}^{\#}\circ\psi^{\#}(\lambda_{A_n}(a_x))\, ,\,
\Phi\circ\alpha_{n,\infty}^{\#}(\lambda_{A_n}(a_x))\, \bigr)+ \\
&\sum_{i=1}^R k_x^i\, 
D_B\bigl(\, \beta_{m,\infty}^{\#}\circ\psi^{\#}(q_{A_n}'(v^{A_n}_i))\, ,\, 
\Phi\circ\alpha_{n,\infty}^{\#}(q_{A_n}'(v^{A_n}_i))\, \bigr) \\
<\ &|e^{2\pi i\delta}-1|+\sum_{i=1}^R k_x^i\delta<\epsilon.
\end{align*}
\end{proof}
\end{prop}

\begin{thm}\label{e}
Let $A$ be a unital simple inductive limit of a sequence
$$ 
  \begin{CD}
      A_1  @> \alpha_1 >> A_2 @> \alpha_2 >>
A_3 @>\alpha_3  >> \dots \\   
  \end{CD}
$$
of finite direct sums of building blocks.
Let $B$ be an inductive limit of a similar sequence
$$
  \begin{CD}
      B_1  @> \beta_1 >> B_2 @> \beta_2 >>
B_3 @>\beta_3  >> \dots \\   
  \end{CD}
$$
with unital connecting maps such that $s(B_k)\to\infty$ and such that
the torsion subgroup of $\text{Aff}\, T(B)/\overline{\rho_B(K_0(B))}$ is
totally disconnected. Let $\varphi_T : T(B)\to T(A)$ be an affine continuous
map, let $\kappa\in KL(A,B)_e$ be an element such that
$$
r_B(\omega)(\kappa_*(x))=r_A(\varphi_T(\omega))(x),
\quad x\in K_0(A),\ \omega\in T(B),
$$
and let $\Phi:U(A)/\overline{DU(A)}\to U(B)/\overline{DU(B)}$ be a
homomorphism such that the diagram
$$
  \begin{CD}
      \text{Aff}\, T(A)/\overline{\rho_A(K_0(A))} @> \lambda_A >>
      U(A)/\overline{DU(A)} @> \pi_A >> K_1(A)    \\   
      @V\widetilde{\varphi_T} VV  @V\Phi VV      @VV\kappa_* V \\
      \text{Aff}\, T(B)/\overline{\rho_B(K_0(B))}
      @>> \lambda_B > U(B)/\overline{DU(B)} @>> \pi_B > K_1(B)
  \end{CD}
$$
commutes. Assume finally that
$$
s_{\kappa}(y)=\Phi(y),\quad y\in Tor(U(A)/\overline{DU(A)}).
$$
There exists a unital *-homomorphism $\psi:A\to B$ such that
$\psi^*=\varphi_T$ on $T(B)$, such that $\psi^{\#}=\Phi$ on
$U(A)/\overline{DU(A)}$, and such that $[\psi]=\kappa$ in $KL(A,B)$.

\begin{proof}
We may assume that $A$ is infinite dimensional. Hence by Theorem \ref{inj} we
may assume that each $\alpha_n$ is unital and injective. Let
$A_n=A^n_1\oplus A^n_2\oplus\dots\oplus A^n_{R_n}$ where each $A^n_i$ is a
building block and let $P_n$ be the set of minimal non-zero central
projections in $A_n$. For each positive integer $n$,
choose a finite set $G_n\subseteq A_n$ such that $G_n$ generates $A_n$ as
a $C^*$-algebra and such that $\alpha_n(G_n)\subseteq G_{n+1}$.
Choose by uniqueness, Theorem \ref{un}, a positive integer $l_n$ with respect
to $G_n\subseteq A_n$ and $2^{-n}$. Since $A$ is simple and the connecting
maps are injective there exists a positive integer $p_n$ such that
$$
\widehat{\alpha_{n,\infty}}(\widehat h)>\frac 8{p_n},\quad h\in H(A_n,l_n).
$$
Next, there exists a positive integer $q_n$ such that
$$
\widehat{\alpha_{n,\infty}}(\widehat h)>\frac 2{q_n},
\quad h\in H(A_n,p_n)\cup P_n.
$$
Finally choose $\delta_n>0$ such that $\delta_n<\frac 1{4q^2}$ and such that
$$
\widehat{\alpha_{n,\infty}}(\widehat h)>\delta_n,\quad h\in H(A_n,4q_n).
$$
Choose for each $n$ finite sets $F_n\subseteq \text{Aff}\, T(A_n)$ such that
$\widetilde{H}(A_n,2q_n)\subseteq F_n$, such that
$\widehat{\alpha_n}(F_n)\subseteq F_{n+1}$, and such that
$\cup_{n=1}^{\infty}\widehat{\alpha_{n,\infty}}(F_n)$ is dense in
$\text{Aff}\, T(A)$.

Next, choose finite sets $V_n\subseteq U(A_n)/\overline{DU(A_n)}$ such that
$q_{A_n}'(v^{A_n}_i)\in V_n$ for $i=1,2,\dots,R_n$,
such that $\alpha_n^{\#}(V_n)\subseteq V_{n+1}$, and such
that $\cup_{n=1}^{\infty}\alpha_{n,\infty}^{\#}(V_n)$ is dense in
$U(A)/\overline{DU(A)}$.

We will construct by induction strictly increasing sequences $\{n_k\}$ and
$\{m_k\}$ and unital *-homomorphisms $\psi_k:A_{n_k}\to B_{m_k}$ such that
\begin{itemize}
\item[(i)]
$\|\beta_{m_{k-1},m_k}\circ\psi_{k-1}(x)-\psi_k\circ\alpha_{n_{k-1},n_k}(x)\|<
2^{-n_{k-1}},\quad x\in G_{n_{k-1}},\ k\ge 2$,
\item[(ii)]
$\|\widehat{\beta_{m_k,\infty}}\circ\widehat{\psi_k}(f)-
{\varphi_T}_*\circ\widehat{\alpha_{n_k,\infty}}(f)\|<
\min \{2^{-n_k},\frac{\delta_{n_k}}2\},\quad f\in F_{n_k}$,
\item[(iii)]
$D_B\bigl(\, \beta_{m_k,\infty}^{\#}\circ\psi_k^{\#}(x)\, ,\,
\Phi\circ\alpha_{n_k,\infty}^{\#}(x)\bigr)<
\min\{2^{-n_k},\frac{\delta_{n_k}}2\},\quad x\in V_{n_k}$,
\item[(iv)]
$[\beta_{m_k,\infty}]\cdot [\psi_k]=\kappa\cdot [\alpha_{n_k,\infty}]
\quad \text{in}\ KL(A_{n_k},B)$.
\end{itemize}
The integers $n_k$, $m_k$, and the *-homomorphism $\psi_k$ are constructed
in step $k$. The case $k=1$ follows immediately from Proposition \ref{exhom}.

Assume that $n_k$, $m_k$, and $\psi_k$ have been constructed such that (i)-(iv)
hold. Choose $n_{k+1}>n_k$ such that
\begin{gather*}
\widehat{\alpha_{n_k,n_{k+1}}}(\widehat h)>\frac 8{p_{n_k}},
\quad h\in H(A_{n_k},l_{n_k}), \\
\widehat{\alpha_{n_k,n_{k+1}}}(\widehat h)>\frac 2{q_{n_k}},
\quad h\in H(A_{n_k},p_{n_k})\cup P_n, \\
\widehat{\alpha_{n_k,n_{k+1}}}(\widehat h)>\delta_{n_k},
\quad h\in H(A_{n_k},4q_{n_k}).
\end{gather*}

Choose by Proposition \ref{exhom} a positive integer $l$
and a unital *-homomorphism $\lambda:A_{n_{k+1}}\to B_l$ such that
\begin{gather*}
\|\widehat{\beta_{l,\infty}}\circ\widehat{\lambda}(f)-
{\varphi_T}_*\circ\widehat{\alpha_{n_{k+1},\infty}}(f)\|<
\min\{ 2^{-n_{k+1}},\frac{\delta_{n_k}}2,\frac{\delta_{n_{k+1}}}2\},
\quad f\in F_{n_{k+1}}, \\
D_B\bigl(\, \beta_{l,\infty}^{\#}\circ\lambda^{\#}(x)\, ,\,
\Phi\circ\alpha_{n_{k+1},\infty}^{\#}(x)\, \bigr)<
\min\{ 2^{-n_{k+1}},\frac{\delta_{n_k}}2,\frac{\delta_{n_{k+1}}}2\},
\quad x\in V_{n_{k+1}}, \\
[\beta_{l,\infty}]\cdot [\lambda]=\kappa\cdot [\alpha_{n_{k+1},\infty}]
\quad\text{in}\ KL(A_{n_{k+1}},B).
\end{gather*}
It follows that
\begin{gather*}
\|\widehat{\beta_{l,\infty}}\circ\widehat{\lambda}
\circ\widehat{\alpha_{n_k,n_{k+1}}}(f)-
\widehat{\beta_{m_k,\infty}}\circ\widehat{\psi_k}(f)\|<\delta_{n_k},
\quad f\in F_{n_k}, \\
D_B\bigl(\, \beta_{l,\infty}^{\#}\circ\lambda^{\#}\circ
\alpha_{n_k,n_{k+1}}^{\#}(x)\, ,\,
\beta_{m_k,\infty}^{\#}\circ\psi_k^{\#}(x)\, \bigr)<
\delta_{n_k}<\frac 1{4q_{n_k}}, \quad x\in V_{n_k}, \\
[\beta_{m_k,\infty}]\cdot [\psi_k]=
[\beta_{l,\infty}]\cdot [\lambda]\cdot [\alpha_{n_k,n_{k+1}}]
\quad \text{in}\ KL(A_k,B).
\end{gather*}
Hence there exists an integer $m_{k+1}\ge l$ such that
\begin{gather*}
\|\widehat{\beta_{l,m_{k+1}}}\circ\widehat{\lambda}
\circ\widehat{\alpha_{n_k,n_{k+1}}}(f)-
\widehat{\beta_{m_k,m_{k+1}}}\circ\widehat{\psi_k}(f)\|<\delta_{n_k},
\quad f\in F_{n_k}, \\
D_B\bigl(\, \beta_{l,m_{k+1}}^{\#}\circ\lambda^{\#}\circ
\alpha_{n_k,n_{k+1}}^{\#}(x)\, ,\,
\beta_{m_k,m_{k+1}}^{\#}\circ\psi_k^{\#}(x)\, \bigr)<\frac 1{4q_{n_k}},
\quad x\in V_{n_k}, \\
[\beta_{l,m_{k+1}}]\cdot [\lambda]\cdot [\alpha_{n_k,n_{k+1}}]=
[\beta_{m_k,m_{k+1}}]\cdot [\psi_k]
\quad \text{in}\ KL(A_k,B_{m_{k+1}}).
\end{gather*}
By uniqueness, Theorem \ref{un}, there exists a unitary $W\in B_{m_{k+1}}$ such
that
$$
\|\beta_{m_k,m_{k+1}}\circ\psi_k(x)-
W\beta_{l,m_{k+1}}\circ\lambda\circ\alpha_{n_k,n_{k+1}}(x)W^*\|<
2^{-n_k},\quad x\in G_{n_k}.
$$
Set $\psi_{k+1}(x)=W\beta_{l,m_{k+1}}\circ\lambda(x)W^*$,
$x\in A_{n_{k+1}}$. It is easily seen that (i)-(iv) are satisfied with $k+1$
in place of $k$. This completes the induction step.

By Elliott's approximate intertwining argument, see e.g \cite[Lemma 1]{TAI},
there exists a *-homomorphism $\psi:A\to B$ such that
$$
\psi(\alpha_{n,\infty}(x))=
\lim_{k\to\infty} \beta_{m_k,\infty}\circ\psi_k\circ\alpha_{n,n_k}(x),
\quad x\in A_n.
$$
Clearly, $\psi$ is unital.
Let $f\in F_n$, $\omega\in T(B)$. The sequence
$\omega\circ\beta_{m_k,\infty}\circ\psi_k\circ\alpha_{n,n_k}$
converges to $\omega\circ\psi\circ\alpha_{n,\infty}$ in $T(A_n)$ as
$k\to\infty$. Hence it follows that
$$
\widehat{\beta_{m_k,\infty}}\circ\widehat{\psi_k}
\circ\widehat{\alpha_{n,n_k}}(f)(\omega)\ \to\
\widehat{\psi}\circ\widehat{\alpha_{n,\infty}}(f)(\omega)\quad \text{as\ }
k\to\infty.
$$
On the other hand, from (ii) it follows that
$$
\widehat{\beta_{m_k,\infty}}\circ\widehat{\psi_k}
\circ\widehat{\alpha_{n,n_k}}(f)(\omega)\ \to\
{\varphi_T}_*\circ\widehat{\alpha_{n,\infty}}(f)(\omega)\quad \text{as\ }
k\to\infty.
$$
Hence $\widehat{\psi}={\varphi_T}_*$ on $\text{Aff}\, T(A)$ and thus
$\psi^*=\varphi_T$ on $T(B)$. If $y,z\in U(A)/\overline{DU(A)}$ then
$D_B(\Phi(y),\Phi(y))\le D_A(y,z)$. This is clear in the case that
$\pi_A(y)\neq\pi_A(z)$ since then $D_A(y,z)=2$, and otherwise it follows since
$\lambda_A$ and $\lambda_B$ are isometries and ${\varphi_T}_*$ is contractive
(with respect to $d_A$ and $d_B$). Thus $\Phi$ is continuous and by
arguments similar to those applied above we see that $\psi^{\#}=\Phi$.

Let finally $n$ be a positive integer. Since $A_n$ is semiprojective there
exists by \cite[Theorem 15.1.1]{LLS} a positive integer $l\ge n$
such that $\psi\circ \alpha_{n,\infty}$ is homotopic to
$\beta_{m_l,\infty}\circ\psi_l\circ\alpha_{n,n_l}$. Hence
$$
[\psi]\cdot [\alpha_{n,\infty}]=
[\beta_{m_l,\infty}]\cdot[\psi_l]\cdot[\alpha_{n,n_l}]=
\kappa\cdot [\alpha_{n_l,\infty}]\cdot [\alpha_{n,n_l}]=
\kappa\cdot[\alpha_{n,\infty}]
$$
in $KL(A_n,B)$. It follows from \cite[Lemma 5.8]{RKL} that $[\psi]=\kappa$
in $KL(A,B)$.
\end{proof}
\end{thm}

The following corollary generalizes a theorem of Thomsen
\cite[Theorem A]{TATD}.

\begin{cor}\label{exnc}
Let $A$ be a unital simple inductive limit of a sequence
$$ 
  \begin{CD}
      A_1  @> \alpha_1 >> A_2 @> \alpha_2 >>
A_3 @>\alpha_3  >> \dots \\   
  \end{CD}
$$
of finite direct sums of building blocks such that $K_0(A)$ is non-cyclic.
Let $B$ be an inductive limit of a similar sequence
$$
  \begin{CD}
      B_1  @> \beta_1 >> B_2 @> \beta_2 >>
B_3 @>\beta_3  >> \dots \\   
  \end{CD}
$$
with unital connecting maps such that $s(B_k)\to\infty$.
Let $\varphi_T : T(B)\to T(A)$ be an affine continuous
map, let $\kappa\in KL(A,B)_e$ be an element such that
$$
r_B(\omega)(\kappa_*(x))=r_A(\varphi_T(\omega))(x),
\quad x\in K_0(A),\ \omega\in T(B),
$$
and let $\Phi:U(A)/\overline{DU(A)}\to U(B)/\overline{DU(B)}$ be a
homomorphism such that the diagram
$$
  \begin{CD}
      \text{Aff}\, T(A)/\overline{\rho_A(K_0(A))} @> \lambda_A >>
      U(A)/\overline{DU(A)} @> \pi_A >> K_1(A)    \\   
      @V\widetilde{\varphi_T} VV  @V\Phi VV      @VV\kappa_* V \\
      \text{Aff}\, T(B)/\overline{\rho_B(K_0(B))}
      @>> \lambda_B > U(B)/\overline{DU(B)} @>> \pi_B > K_1(B)
  \end{CD}
$$
commutes. There exists a unital *-homomorphism $\psi:A\to B$ such that
$\psi^*=\varphi_T$ on $T(B)$, such that $\psi^{\#}=\Phi$ on
$U(A)/\overline{DU(A)}$, and such that $[\psi]=\kappa$ in $KL(A,B)$.

\begin{proof}
By Lemma \ref{K_0nct} we have that
$\text{Aff}\, T(B)/\overline{\rho_B(K_0(B))}$ is torsion free such that
$s_{\kappa}$ is defined. It follows by Proposition \ref{exact} that
$s_{\kappa}(y)=\Phi(y)$ for $y$ in the torsion subgroup of
$U(A)/\overline{DU(A)}$. Apply Theorem \ref{e}.
\end{proof}
\end{cor}

\begin{cor}
Let $A$ be a unital inductive limit of a sequence
$$ 
  \begin{CD}
      A_1  @> \alpha_1 >> A_2 @> \alpha_2 >>
A_3 @>\alpha_3  >> \dots \\   
  \end{CD}
$$
of finite direct sums of building blocks.
Let $B$ be an inductive limit of a similar sequence
$$
  \begin{CD}
      B_1  @> \beta_1 >> B_2 @> \beta_2 >>
B_3 @>\beta_3  >> \dots \\   
  \end{CD}
$$
with unital connecting maps such that $s(B_k)\to\infty$ and such that
the torsion subgroup of $\text{Aff}\, T(B)/\overline{\rho_B(K_0(B))}$ is
totally disconnected. Let $\varphi_T : T(B)\to T(A)$ be an affine continuous
map, let $\varphi_0:K_0(A)\to K_0(B)$ be an order unit preserving group
homomorphism such that
$$
r_B(\omega)(\varphi_0(x))=r_A(\varphi_T(\omega))(x),
\quad x\in K_0(A),\ \omega\in T(B),
$$
and let $\varphi_1:K_1(A)\to K_1(B)$ be a group homomorphism.
There exists a unital *-homomorphism $\psi:A\to B$ such that
$\psi^*=\varphi_T$ on $T(B)$, such that $\psi_*=\varphi_0$ on
$K_0(A)$, and such that $\psi_*=\varphi_1$ on $K_1(A)$.

\begin{proof}
Choose an element $\kappa\in KL(A,B)$ such that $\kappa_*=\varphi_0$
on $K_0(A)$ and such that $\kappa_*=\varphi_1$ on $K_1(A)$. By
Lemma \ref{lifttophi} there exists a group homomorphism
$\Phi:U(A)/\overline{DU(A)}\to U(B)/\overline{DU(B)}$ such that
$s_{\kappa}$ and $\Phi$ agree on the torsion subgroup of
$U(A)/\overline{DU(A)}$ and such that the diagram
$$
  \begin{CD}
      \text{Aff}\, T(A)/\overline{\rho_A(K_0(A))} @> \lambda_A >>
      U(A)/\overline{DU(A)} @> \pi_A >> K_1(A)    \\   
      @V\widetilde{\varphi_T} VV  @V\Phi VV      @VV\kappa_* V \\
      \text{Aff}\, T(B)/\overline{\rho_B(K_0(B))}
      @>> \lambda_B > U(B)/\overline{DU(B)} @>> \pi_B > K_1(B)
  \end{CD}
$$
commutes. The conclusion follows from Theorem \ref{e}.
\end{proof}
\end{cor}

\begin{thm}\label{aue}
Let $A$ and $B$ be unital inductive limits of sequences of finite direct
sums of building blocks, with $A$ simple.
Let $\varphi,\psi:A\to B$ be unital *-homomorphisms such that
$\varphi^*=\psi^*$ on $T(B)$, $\varphi^{\#}=\psi^{\#}$ on
$U(A)/\overline{DU(A)}$, and $[\varphi]=[\psi]$ in $KL(A,B)$. Then $\varphi$
and $\psi$ are approximately unitarily equivalent.

\begin{proof}
We may assume that $A$ is infinite dimensional, and hence by Theorem \ref{inj}
we see that $A$ is the inductive limit of a sequence
$$
  \begin{CD}
      A_1  @> \alpha_1 >> A_2 @> \alpha_2 >>
A_3 @>\alpha_3  >> \dots \\   
  \end{CD}
$$
of finite direct sums of building blocks with unital and injective connecting
maps. By Lemma \ref{unital} we have that $B$ is the inductive limit of a
sequence
$$
  \begin{CD}
      B_1  @> \beta_1 >> B_2 @> \beta_2 >>
B_3 @>\beta_3  >> \dots \\   
  \end{CD}
$$
of finite direct sums of building blocks with unital connecting maps.
Let $A_n=A_1^n\oplus A_2^n\oplus\dots\oplus A_{R_n}^n$ where each $A_i^n$ is
a building block. Let $P_n$ be the set of minimal non-zero central projections
in $A_n$.

Let $F\subseteq A$ be a finite set and let $\epsilon>0$. It suffices to see
that there exists a unitary $U\in B$ such that
$$
\|\varphi(x)-U\psi(x)U^*\|<\epsilon,\quad x\in F.
$$
We may assume that $F\subseteq\alpha_{n,\infty}(G)$ for a positive integer
$n$ and a finite set $G\subseteq A_n$.

Choose by uniqueness, Theorem \ref{un}, a positive integer $l$ with respect to
$G$ and $\frac{\epsilon}3$. Since $A$ is simple and the connecting maps are
injective there exists an integer $p\ge l$ such that
$$
\widehat{\alpha_{n,\infty}}(\widehat{h})>\frac 9p,\quad h\in H(A_n,l).
$$
Next choose $q\ge p$ such that
$$
\widehat{\alpha_{n,\infty}}(\widehat{h})>\frac 3q,
\quad h\in H(A_n,p)\cup P_n.
$$
Finally, choose $\delta>0$ such that $3\delta<\epsilon$,
$2\delta<\frac 1{4q^2}$, and such that
$$
\widehat{\alpha_{n,\infty}}(\widehat{h})>3\delta,\quad h\in H(A_n,4q).
$$
Since $A_n$ by Theorem \ref{stable} is semiprojective there exist by
\cite[Corollary 15.1.3]{LLS} a positive integer $r$ and *-homomorphisms
$\varphi_1,\psi_1:A_n\to B_r$ such that $\beta_{r,\infty}\circ\varphi_1$ is
homotopic to $\varphi\circ\alpha_{n,\infty}$ and
$\beta_{r,\infty}\circ\psi_1$ is homotopic to $\psi\circ\alpha_{n,\infty}$,
and such that if
$$
x\in G\cup H(A_n,l)\cup H(A_n,p)\cup H(A_n,4q)\cup\widetilde{H}(A_n,2q)\cup
P_n\cup\{v^{A_n}_1,v^{A_n}_2,\dots,v^{A_n}_{R_n}\}
$$
then
\begin{gather*}
\|\beta_{r,\infty}\circ\varphi_1(x)-\varphi\circ\alpha_{n,\infty}(x)\|<
\delta, \\
\|\beta_{r,\infty}\circ\psi_1(x)-\psi\circ\alpha_{n,\infty}(x)\|<\delta.
\end{gather*}
By increasing $r$ we may assume that $\varphi_1$
and $\psi_1$ are unital. Note that
\begin{gather*}
D_B\bigl(\, \beta_{r,\infty}^{\#}\circ\varphi_1^{\#}(q'_{A_n}(v_i^{A_n}))\, ,\,
\beta_{r,\infty}^{\#}\circ\psi_1^{\#}(q'_{A_n}(v_i^{A_n}))\bigr)<
2\delta<\frac 1{4q^2},\quad i=1,2,\dots,R_n, \\
\|\widehat{\beta_{r,\infty}}\circ\widehat{\varphi_1}(\widehat{h})-
\widehat{\beta_{r,\infty}}\circ\widehat{\psi_1}(\widehat{h})\|<2\delta,
\quad h\in \widetilde{H}(A_n,2q),
\end{gather*}
and
\begin{gather}
\widehat{\beta_{r,\infty}}\circ\widehat{\psi_1}(\widehat{h})>2\delta,
\quad h\in H(A_n,4q), \\
\widehat{\beta_{r,\infty}}\circ\widehat{\psi_1}(\widehat{h})>\frac 8p,
\quad h\in H(A_n,l), \\
\widehat{\beta_{r,\infty}}\circ\widehat{\psi_1}(\widehat{h})>\frac 2q,
\quad h\in H(A_n,p)\cup P_n, \\
\label{eqKKuni}
[\beta_{r,\infty}]\cdot[\psi_1]=[\beta_{r,\infty}]\cdot[\varphi_1]
\quad\text{in}\ KL(A_n,B).
\end{gather}
Choose an integer $m\ge r$ such that
\begin{gather}
\|\widehat{\beta_{r,m}}\circ\widehat{\varphi_1}(\widehat{h})-
\widehat{\beta_{r,m}}\circ\widehat{\psi_1}(\widehat{h})\|<
2\delta,\quad h\in\widetilde{H}(A_n,2q), \\
\widehat{\beta_{r,m}}\circ\widehat{\psi_1}(\widehat{h})>2\delta,
\quad h\in H(A_n,4q), \\
\widehat{\beta_{r,m}}\circ\widehat{\psi_1}(\widehat{h})>\frac 8p,
\quad h\in H(A_n,l), \\
\widehat{\beta_{r,m}}\circ\widehat{\psi_1}(\widehat{h})>\frac 2q,
\quad h\in H(A_n,p)\cup P_n, \\
D_{B_m}\bigl(\, \beta_{r,m}^{\#}\circ\varphi_1^{\#}(q'_{A_n}(v_i^{A_n}))\, ,\,
\beta_{r,m}^{\#}\circ\psi_1^{\#}(q'_{A_n}(v_i^{A_n}))\bigr)<
\frac 1{4q^2},\quad i=1,\dots,R_n, \\
\label{eqKKuni2}
[\beta_{r,m}]\cdot[\psi_1]=[\beta_{r,m}]\cdot[\varphi_1]
\quad\text{in}\ KL(A_n,B_m).
\end{gather}
By Theorem \ref{un} there exists a unitary $W\in B_m$ such that
\begin{equation}\label{eqKKuni3}
\|\beta_{r,m}\circ\varphi_1(x)-W\beta_{r,m}\circ\psi_1(x) W^*\|<
\frac{\epsilon}3,\quad x\in G.
\end{equation}
If we put $U=\beta_{m,\infty}(W)$ we have that
\begin{align*}
&\|\varphi\circ\alpha_{n,\infty}(x)-
U\psi\circ\alpha_{n,\infty}(x)U^*\| \\
\le\ &\|\varphi\circ\alpha_{n,\infty}(x)-\beta_{r,\infty}\circ\varphi_1(x)\|+
\|\beta_{r,\infty}\circ\varphi_1(x)-U\beta_{r,\infty}\circ\psi_1(x)U^*\|\ + \\
&\|\beta_{r,\infty}\circ\psi_1(x)-\psi\circ\alpha_{n,\infty}(x)\| \\
<\ &\delta+\frac{\epsilon}3+\delta<\epsilon,\quad x\in G.
\end{align*}
\end{proof}
\end{thm}

In view of Theorem \ref{uns} one might think that equality in $KL$ in the
above theorem could be replaced by equality in $K_0$. This is however
impossible in general, see \cite[p. 375-376]{DLUMCT} or
\cite[Theorem 8.4]{TATD}. But in some cases, e.g when $K_0(B)$ is cyclic, the
$KL$-condition can be relaxed:

\begin{thm}
Assume furthermore that $\rho_B$ is injective and $\rho_B(K_0(B))$ is a
discrete subgroup of $\text{Aff}\, T(B)$.
If $\varphi_*=\psi_*$ on $K_0(A)$, $\varphi^*=\psi^*$ on $T(B)$ and
$\varphi^{\#}=\psi^{\#}$ on $U(A)/\overline{DU(A)}$, then $\varphi$ and $\psi$
are approximately unitarily equivalent.

\begin{proof}
As above, but with the following changes. Instead of (\ref{eqKKuni}) we get
by Proposition \ref{homotopy} that
\begin{gather*}
\beta_{r,\infty}^{\#}\circ\psi_1^{\#}(x)=
\beta_{r,\infty}^{\#}\circ\varphi_1^{\#}(x),\quad x\in U^{A_n}, \\
\beta_{r,\infty}\circ{\psi_1}_*=\beta_{r,\infty}\circ{\varphi_1}_*
\quad\text{on\ }K_0(A_n).
\end{gather*}
By Lemma \ref{eta} we may now replace (\ref{eqKKuni2}) by
\begin{gather*}
\beta_{r,m}^{\#}\circ\psi_1^{\#}(x)=
\beta_{r,m}^{\#}\circ\varphi_1^{\#}(x),\quad x\in U^{A_n}, \\
\beta_{r,m}\circ{\psi_1}_*=\beta_{r,m}\circ{\varphi_1}_*
\quad\text{on\ }K_0(A_n).
\end{gather*}
Finally, (\ref{eqKKuni3}) follows again by Theorem \ref{un}.
\end{proof}
\end{thm}

\begin{thm}\label{i}
Let $A$ and $B$ be simple unital infinite dimensional inductive limits of
sequences of finite direct sum of building blocks. Let
$\varphi_0:K_0(A)\to K_0(B)$ be an isomorphism of groups with order units,
let $\varphi_1:K_1(A)\to K_1(B)$ be an isomorphism of groups, and let
$\varphi_T:T(B)\to T(A)$ be an affine homeomorphism such that
$$
r_B(\omega)(\varphi_0(x))=r_A(\varphi_T(\omega))(x),
\quad x\in K_0(A),\ \omega\in T(B).
$$
There exists a *-isomorphism $\varphi:A\to B$ such that
$\varphi_*=\varphi_0$ on $K_0(A)$, such that $\varphi_*=\varphi_1$ on $K_1(A)$,
and such that $\varphi_T=\varphi^*$ on $T(B)$. 

\begin{proof}
By Theorem \ref{inj} we may assume that $A$ is the inductive limit of a
sequence
$$ 
  \begin{CD}
      A_1  @> \alpha_1 >> A_2 @> \alpha_2 >>
A_3 @>\alpha_3  >> \dots \\   
  \end{CD}
$$
of finite direct sums of building blocks with unital and injective connecting
maps. Similarly we may assume that $B$ is the inductive limit of a sequence
$$
  \begin{CD}
      B_1  @> \beta_1 >> B_2 @> \beta_2 >>
B_3 @>\beta_3  >> \dots \\   
  \end{CD}
$$
of finite direct sums of building blocks with unital and injective connecting
maps. By Lemma \ref{infty} we have that
$s(A_n)\to\infty$ and $s(B_n)\to\infty$ as $n\to\infty$.

By \cite[Theorem 7.3]{RSUCT} there exists an invertible element
$\kappa\in KL(A,B)$ such that $\kappa_*=\varphi_0$ on $K_0(A)$ and
$\kappa_*=\varphi_1$ on $K_1(A)$.
By Lemma \ref{lifttophi} there exists a group isomorphism
$\Phi:U(A)/\overline{DU(A)}\to U(B)/\overline{DU(B)}$ such that the diagram
$$
  \begin{CD}
      0 @> >>\text{Aff}\, T(A)/\overline{\rho_A(K_0(A))} @> \lambda_A >>
      U(A)/\overline{DU(A)}   @> \pi_A >> K_1(A) @> >> 0   \\   
      @. @V\widetilde{\varphi_T} VV         @V\Phi VV      @VV\kappa_* V \\
      0 @> >> \text{Aff}\, T(B)/\overline{\rho_B(K_0(B))} @>> \lambda_B >
      U(B)/\overline{DU(B)}   @>> \pi_B > K_1(B) @> >> 0
  \end{CD}
$$
commutes and such that $s_{\kappa}(y)=\Phi(y)$ for $y$ in the torsion subgroup
of $U(A)/\overline{DU(A)}$.

By Theorem \ref{e} there exists a unital *-homomorphism $\lambda:A\to B$ such
that $\lambda^*=\varphi_T$ on $T(B)$, such that $\lambda^{\#}=\Phi$ on
$U(A)/\overline{DU(A)}$, and such that $[\lambda]=\kappa$ in $KL(A,B)$.
Note that $\kappa^{-1}\in KL(B,A)_T$.
It is easy to see that $s_{\kappa}$ is a bijection with inverse
$s_{\kappa^{-1}}$. Hence $s_{\kappa^{-1}}=\Phi^{-1}$ on
$Tor(U(B)/\overline{DU(B)})$.
Thus there exists a unital *-homomorphism $\psi:A\to B$ such that
$\psi^*={\varphi_T}^{-1}$ on $T(A)$, such that $\psi^{\#}=\Phi^{-1}$ on
$U(B)/\overline{DU(B)}$, and such that $[\psi]=\kappa^{-1}$ in $KL(B,A)$.

By Theorem \ref{aue} the *-homomorphisms $\psi\circ\lambda$ and
$id_A$ are approximately unitarily equivalent. Similarly
$\lambda\circ\psi$ and $id_B$ are approximately unitarily equivalent.
Thus there are sequences of unitaries $\{u_n\}$ and $\{v_n\}$, in $A$ and $B$
respectively, such that, upon setting $\lambda_n(x)=v_n\lambda(x)v_n^*$ and
$\psi_n(x)=u_n\psi(x)u_n^*$, the diagram
$$
\xymatrix{
A \ar[r]^{id_A}\ar[d]^{\lambda_1}  & A \ar[r]^{id_A}\ar[d]^{\lambda_2}  
& A \ar[r]^{id_A}\ar[d]^{\lambda_3}  & \dots \\
B \ar[ur]_{\psi_1}\ar[r]_{id_B} & B \ar[ur]_{\psi_2}\ar[r]_{id_B}
& B \ar[ur]_{\psi_3}\ar[r]_{id_B} & \dots}
$$
becomes an approximate intertwining. Hence by e.g \cite[Theorem 3]{TAI}
there is a *-isomorphism $\varphi:A\to B$ such that
$$
\varphi(x)=\lim_{n\to\infty} v_n\lambda(x)v_n^*,\quad x\in A.
$$
It follows that $\varphi^*=\lambda^*=\varphi_T$ on $T(B)$, that
$\varphi_*=\lambda_*=\varphi_0$ on $K_0(A)$, and that
$\varphi_*=\lambda_*=\varphi_1$ on $K_1(A)$.
\end{proof}
\end{thm}

\section{Range of the invariant}

The purpose of this section is to determine the range of the Elliott
invariant, i.e to answer the question which quadruples
$(K_0(A),K_1(A),T(A),r_A)$ occur as the Elliott invariant for simple unital
infinite dimensional $C^*$-algebras that are inductive limits of sequences
of finite direct sums of building blocks. Villadsen \cite{VREI} has answered
this question in the case where $A$ is an inductive limit of a sequence of
finite direct sums of circle algebras. Using this result Thomsen has
been able to determine the range of the Elliott invariant for those
$C^*$-algebras that are inductive limits of finite direct sums of building
blocks of the form $A(n,d,d,\dots,d)$, see below.

We start out by examining the restrictions on $(K_0(A),K_1(A),T(A),r_A)$.
Let $A$ be a simple unital infinite dimensional inductive limit of a sequence
$$
  \begin{CD}
      A_1  @> \alpha_1 >> A_2 @> \alpha_2 >>
A_3 @>\alpha_3  >> \dots \\   
  \end{CD}
$$
of finite direct sums of building blocks. We may by Theorem \ref{inj} assume
that each $\alpha_n$ is unital and injective.
By Corollary \ref{K_0c} each $K_0(A_k)$ is isomorphic (as an ordered group with
order unit) to the $K_0$-group of a finite dimensional $C^*$-algebra. Thus
$K_0(A)$ must be a countable dimension group. This group has to be simple as
$A$ is simple. 

If $K_0(A)\cong\Z$ then by passing to a subsequence, if necessary,
we may assume that $A$ is the inductive limit of a sequence of building blocks,
rather than finite direct sums of such algebras. By Lemma \ref{K_1ex} it
follows that $K_1(A)$ an inductive limit of groups of the form $\Z\oplus H$,
where $H$ is any finite abelian group.

If $K_0(A)$ is not cyclic our only immediate conclusion is that $K_1(A)$
is a countable abelian group.

$T(A)$ must be a metrizable Choquet simplex.
If $B$ is a building block then obviously $r_B:T(B)\to SK_0(B)$ maps extreme
points to extreme points. By \cite[Corollary 1.6]{VREI} and
\cite[Corollary 1.7]{VREI} the same must be the case for $r_A$. Finally,
$r_A$ is surjective by either \cite[Theorem 3.3]{BR} and \cite{HQT}, or
\cite[Corollary 9.18]{HTQT} (or more
elementary, because each $r_{A_k}:T(A_k)\to SK_0(A_k)$ is surjective).
It follows from Theorem \ref{rannc} and Corollary \ref{ran} that these are the
only restrictions.

As mentioned above, Thomsen has calculated the range of the invariant for
a subclass of the class we are considering.
By \cite[Theorem 9.2]{TATD} we have the following:

\begin{thm}\label{rannc}
Let $G$ be a countable simple dimension group with order unit, $H$
a countable abelian group, $\Delta$ a compact metrizable Choquet simplex, and
$\lambda:\Delta\to SG$ an affine continuous extreme point preserving
surjection. There exists a simple unital infinite dimensional inductive
limit of a sequence of finite direct sums of building blocks $A$ together with
an isomorphism $\varphi_0:K_0(A)\to G$ of ordered groups with order unit, an
isomorphism $\varphi_1:K_1(A)\to H$, and an affine homeomorphism
$\varphi_T:\Delta\to T(A)$ such that
$$
r_A(\varphi_T(\omega))(x)=\lambda(\omega)(\varphi_0(x)),
\quad\omega\in\Delta,\ x\in K_0(A)
$$
if and only if $G$ is non-cyclic.

$A$ can be realized as an inductive limit of a sequence of finite direct sums
of circle algebras and interval building blocks of the form $I(n,d,d)$.
\end{thm}

A different proof of this theorem could be based on Theorem \ref{ex2}
and \cite[Theorem 4.2]{VREI}. Combining the above theorem with Theorem \ref{i}
we get the following:

\begin{thm}\label{isointbbcir}
Let $A$ be a simple unital inductive limit of a sequence of finite direct sums
of building blocks such that $K_0(A)$ is non-cyclic. Then $A$ is the inductive
limit of a sequence of finite direct sums of circle algebras and interval
building blocks of the form $I(n,d,d)$.
\end{thm}

We are left with the case of cyclic $K_0$-group. Note that the equation
$$
r_A(\varphi_T(\omega))(x)=\lambda(\omega)(\varphi_0(x)),
\quad\omega\in\Delta,\ x\in K_0(A)
$$
is trivial when $A$ is a unital $C^*$-algebra with $K_0(A)\cong\Z$.

\begin{lemma}\label{projless}
Let $A$ be a simple unital inductive limit of a sequence of finite direct
sums of building blocks. Then $(K_0(A),K_0(A)^+,[1])\cong (\Z,\Z^+,1)$
if and only if $A$ is unital projectionless.
\begin{proof}
This follows easily from Theorem \ref{inj} and Lemma \ref{proj}.
\end{proof}
\end{lemma}

\begin{thm}\label{ranc}
Let $\Delta$ be a metrizable Choquet simplex, and let $H$ be the inductive
limit of a sequence
$$
  \begin{CD}
      \Z\oplus H_1  @>h_1 >> \Z\oplus H_2 @>h_2 >> \Z\oplus H_3 @>h_3 >>
 \dots \\   
  \end{CD}
$$
where each $H_k$ is a finite abelian group. There exists an infinite
dimensional simple unital projectionless $C^*$-algebra $A$ that is an
inductive limit of a sequence of building blocks, with $K_1(A)\cong H$ and
such that $T(A)$ is affinely homeomorphic to $\Delta$.

\begin{proof}
By \cite[Lemma 3.8]{TAIT} $\text{Aff}\, \Delta$ is isomorphic to an inductive
limit in the category of order unit spaces of a sequence
$$
  \begin{CD}
      C_{\R}[0,1]  @>  >> C_{\R}[0,1]  @>  >> C_{\R}[0,1]  @>  >> \dots \\   
  \end{CD}
$$
It is easy to see that this implies that $\text{Aff}\, \Delta$ is isomorphic to
an inductive limit of a sequence of the form
$$
  \begin{CD}
      C_{\R}(\T)  @> \Theta_1 >> C_{\R}(\T)  @> \Theta_2 >>
C_{\R}(\T)  @>\Theta_3  >> \dots \\   
  \end{CD}
$$
Choose a dense sequence $\{x_k\}_{k=1}^{\infty}$ in $C_{\R}(\T)$ and
a dense sequence $\{z_k\}_{k=1}^{\infty}$ in $\T$.

For every positive integer $k$ we will construct a unital projectionless
building block $A_k$ such that $K_1(A_k)\cong \Z\oplus H_k$, and a
unital and injective *-homomorphism $\alpha_k:A_k\to A_{k+1}$ such that
the (constant) functions $z\mapsto z_1,z\mapsto z_2,\dots,z\mapsto z_k$
are eigenvalue functions for
$\alpha_k$, such that ${\alpha_k}_*=h_k$ on $K_1(A_k)$ (under the
identification $K_1(A_k)\cong \Z\oplus H_k$) and such that
$$
\|\widehat{\alpha_k}(f)-\Theta_k(f)\|<2^{-k},\quad f\in F_k,
$$
under the identification $\text{Aff}\, T(A_k)\cong C_{\R}(\T)$, where
$$
F_k=\{x_1,x_2,\dots,x_k\}
\bigcup_{j=1}^{k-1} \Theta_{j,k}(\{x_1,x_2,\dots,x_k\})
\bigcup_{j=1}^{k-1} \widehat{\alpha_{j,k}}(\{x_1,x_2,\dots,x_k\}).
$$

First choose by Lemma \ref{K_1ex} a unital projectionless building block $A_1$
such that $K_1(A_1)\cong \Z\oplus H_1$.

Assume that $A_k$ has been constructed. We will construct $A_{k+1}$ and
$\alpha_k$. Choose $K$ by Theorem \ref{ex2} with respect to
$F_k\subseteq \text{Aff}\, T(A_k)$, $\epsilon=2^{-k}$ and the integer $k+1$.
By Lemma \ref{K_1ex} there exists a unital projectionless building block
$A_{k+1}$ such that $s(A_{k+1})\ge K$ and $K_1(A_{k+1})\cong \Z\oplus H_{k+1}$.
By Theorem \ref{ex2} there exists a unital *-homomorphism
$\alpha_k:A_k\to A_{k+1}$ such that the identity function on $\T$ and each of
the functions $z\mapsto z_1,z\mapsto z_2,\dots,z\mapsto z_k$
are among the eigenvalue functions for $\alpha_k$ and such that
\begin{gather*}
\|\widehat{\alpha_k}(f)-\Theta_k(f)\|<2^{-k},
\quad f\in F_k, \\
{\alpha_k}_*=h_k\quad\text{on}\ K_1(A_k).
\end{gather*}
This completes the construction.

Set $A=\varinjlim(A_k,\alpha_k)$. $A$ is
infinite dimensional since the connecting maps are injective, and it is unital
projectionless since the connecting maps are unital.
By \cite[Lemma 3.4]{TAIT} $\text{Aff}\, T(A)\cong
\varinjlim (C_{\R}[0,1],\Theta_k)\cong \text{Aff}\, \Delta$, and
hence $T(A)$ and $\Delta$ are affinely homeomorphic. Clearly $K_1(A)\cong H$.

Let $I\subseteq A$ be a closed two-sided ideal in $A$, $I\neq\{0\}$. By
(the proof of) \cite[Lemma 3.1]{BAF},
$$
I=\overline{\bigcup_{n=1}^{\infty}
\alpha_{n,\infty}({\alpha_{n,\infty}}^{-1}(I))}.
$$
Choose a positive integer $n$ such that
${\alpha_{n,\infty}}^{-1}(I)\neq \{0\}$. Choose
$f\in {\alpha_{n,\infty}}^{-1}(I)$ such that $f\neq 0$. Choose $k>n$
such that $f(z_k)\neq 0$. Then $\alpha_{n,l}(f)(z)\neq 0$ for every $z\in\T$
and $l>k$. Hence by Lemma \ref{ideal} we see that
${\alpha_{l,\infty}}^{-1}(I)=A_l$ for every $l>k$. It follows that $I=A$.
Thus $A$ is simple.
\end{proof}
\end{thm}

In the above theorem, let $H=0$ and $\Delta$ be a one-point set. Then
we obtain by Lemma \ref{projless} and Theorem \ref{i} the $C^*$-algebra
$\mathcal Z$ constructed by Jiang and Su \cite{JSSP}.

\begin{cor}\label{ran}
Let $d$ be a positive integer, let $\Delta$ be a metrizable Choquet simplex
and let $H$ be a countable abelian group.
There exists an infinite dimensional simple unital inductive limit of a
sequence of finite direct sums of building blocks $A$ such that
$(K_0(A),K_0(A)^+,[1])\cong (\Z,\Z^+,d)$, $T(A)\cong \Delta$ and
$K_1(A)\cong H$ if and only if $H$ is the inductive limit of a sequence
$$
  \begin{CD}
      \Z\oplus H_1  @> >> \Z\oplus H_2 @> >> \Z\oplus H_3 @> >> \dots \\   
  \end{CD}
$$
where each $H_k$ is a finite abelian group.

The $C^*$-algebra $A$ is isomorphic to $M_d(B)$ where $B$ is a simple unital
projectionless $C^*$-algebra that is an inductive limit of a
sequence of building blocks.

\begin{proof}
Combine Theorem \ref{ranc}, Lemma \ref{projless} and Theorem \ref{i}.
\end{proof}
\end{cor}

Theorem \ref{rannc} and Corollary \ref{ran} together determine the range of the
Elliott invariant for the class of $C^*$-algebras for which our classification
theorem applies. Let us conclude this paper by comparing our classification
theorem with the classification theorems of Thomsen \cite{TATD} and
Jiang and Su \cite{JSSP}.

It follows from \cite[Theorem 9.2]{TATD} that a $C^*$-algebra in our class
is contained in Thomsen's class if and only if $K_0$ is non-cyclic. By
calculating the range of the invariant for the $C^*$-algebras contained in
Jiang's and Su's class, one can show that
a $C^*$-algebra in our class with $K_0$ non-cyclic is contained in Jiang's
and Su's class if and only if $K_1$ is a torsion group. A $C^*$-algebra in our
class with cyclic $K_0$-group is contained in Jiang's and Su's class if
and only if the $K_1$-group is an inductive limit of a sequence of finite
cyclic groups, see \cite[Theorem 4.5]{JSSP}. Thus our classification theorem
can be applied to $C^*$-algebras that cannot be realized as inductive limits
of finite direct sums of the building blocks considered in \cite{TATD}, or in
\cite{JSSP}, namely those that have cyclic $K_0$-group and a $K_1$-group that
is not an inductive limit of a sequence of finite cyclic groups.

\bibliographystyle{amsplain}
\bibliography{mygind}

\providecommand{\bysame}{\leavevmode\hbox to3em{\hrulefill}\thinspace}
\begin{thebibliography}{10}

\bibitem{BTAF}
B.~Blackadar, \emph{Traces on simple {$AF$} {$C^*$}-algebras}, J. Funct. Anal.
  \textbf{38} (1980), 156--168.

\bibitem{BKOA}
\bysame, \emph{{$K$}-theory for operator algebras}, Cambridge University Press,
  1998.

\bibitem{BR}
B.~Blackadar and M.~R{\o}rdam, \emph{Extending states on preordered semigroups
  and the existence of quasitraces on {$C^*$}-algebras}, J. Algebra
  \textbf{152} (1992), 240--247.

\bibitem{BOLLO}
B.~Bollob{\'a}s, \emph{Combinatorics}, Cambridge University Press, 1986.

\bibitem{BAF}
O.~Bratteli, \emph{Inductive limits of finite dimensional {$C^*$}-algebras},
  Trans. Amer. Math. Soc. \textbf{171} (1972), 195--234.

\bibitem{DLUMCT}
M.~Dadarlat and T.~A. Loring, \emph{A universal multicoefficient theorem for
  the {K}asparov groups}, Duke Math. J. \textbf{84} (1996), 355--377.

\bibitem{ELPNCCW}
S.~Eilers, T.~A. Loring, and G.~K. Pedersen, \emph{Stability of anticommutation
  relations: {A}n application of noncommutative {$CW$} complexes}, J. Reine
  Angew. Math. \textbf{499} (1998), 101--143.

\bibitem{EAIS}
G.~A. Elliott, \emph{A classification of certain simple {$C^*$}-algebras},
  Quantum and non-commutative analysis, H. Araki et al. (eds.), Kluwer,
  Dordrecht (1993), 373--385.

\bibitem{ERR0}
\bysame, \emph{On the classification of {$C^*$}-algebras of real rank zero}, J.
  Reine Angew. Math. \textbf{443} (1993), 179--219.

\bibitem{EATS}
\bysame, \emph{A classification of certain simple {$C^*$}-algebras, {II}}, J.
  Ramanujan Math. Soc. \textbf{12} (1997), 97--134.

\bibitem{FIAG}
L.~Fuchs, \emph{Infinite abelian groups {I}}, Academic Press, 1970.

\bibitem{HQT}
U.~Haagerup, \emph{Quasitraces on exact {$C^*$}-algebras are traces},
  manuscript (1991).

\bibitem{HTQT}
U.~Haagerup and S.~Thorbj\o rnsen, \emph{Random matrices and {$K$}-theory for
  exact {$C^*$}-algebras}, preprint (1998).

\bibitem{HKAF}
D.~Handelman, \emph{{$K_0$} of von {N}eumann and {$AF$} {$C^*$}-algebras},
  Quart. J. Math. Oxford Ser. \textbf{29} (1978), 427--441.

\bibitem{JSSI}
X.~Jiang and H.~Su, \emph{A classification of splitting interval algebras}, J.
  Funct. Anal. \textbf{151} (1997), 50--76.

\bibitem{JSSP}
\bysame, \emph{On a simple unital projectionless {$C^*$}-algebra}, Amer. J.
  Math. \textbf{121} (1999), 359--413.

\bibitem{LIPHD}
L.~Li, \emph{Simple inductive limit {$C^*$}-algebras: {S}pectra and
  approximations by interval algebras}, J. Reine Angew. Math. \textbf{507}
  (1999), 57--79.

\bibitem{LLS}
T.~A. Loring, \emph{Lifting solutions to perturbing problems in
  {$C^*$}-algebras}, Fields Institute Monographs, 1997.

\bibitem{MURPHY}
G.~J. Murphy, \emph{{$C^*$}-algebras and operator theory}, Academic Press,
  1990.

\bibitem{NT}
K.~E. Nielsen and K.~Thomsen, \emph{Limits of circle algebras}, Exposition.
  Math. \textbf{14} (1996), 17--56.

\bibitem{RKL}
M.~R\o rdam, \emph{Classification of certain infinite simple {$C^*$}-algebras},
  J. Funct. Anal. \textbf{131} (1995), 415--458.

\bibitem{RSUCT}
J.~Rosenberg and C.~Schochet, \emph{The {K}\"unneth theorem and the universal
  coefficient theorem for {K}asparov's generalized {$K$}-functor}, Duke Math.
  J. \textbf{55} (1987), 431--474.

\bibitem{TNS}
K.~Thomsen, \emph{Nonstable {$K$}-theory for operator algebras}, {$K$}-Theory
  \textbf{4} (1991), 245--267.

\bibitem{TAI}
\bysame, \emph{On isomorphisms of inductive limit {$C^*$}-algebras}, Proc.
  Amer. Math. Soc. \textbf{113} (1991), 947--953.

\bibitem{TAIT}
\bysame, \emph{Inductive limits of interval algebras: The tracial state space},
  Amer. J. Math. \textbf{116} (1994), 605--620.

\bibitem{TATD}
\bysame, \emph{Limits of certain subhomogeneous {$C^*$}-algebras}, M{\'e}m.
  Soc. Math. Fr. \textbf{71} (1997).

\bibitem{VREI}
J.~Villadsen, \emph{The range of the {E}lliott invariant}, J. Reine Angew.
  Math. \textbf{462} (1995), 31--55.

\end{thebibliography}

\end{document}